\documentclass[12pt]{amsart}
\usepackage{graphicx, overpic, epsf,amssymb,amscd,amsfonts,diagrams,times}

\newtheorem{thm}{Theorem}[section]
\newtheorem*{claim}{Claim}

\newtheorem{prop}[thm]{Proposition}

\newtheorem{lemma}[thm]{Lemma}
\newtheorem{cor}[thm]{Corollary}
\newtheorem{conj}[thm]{Conjecture}
\newtheorem{q}[thm]{Question}
\newtheorem{rmk}[thm]{Remark}

\newcommand{\C}{\mathbb{C}}
\newcommand{\R}{\mathbb{R}}
\newcommand{\Z}{\mathbb{Z}}

\newcommand{\bdry}{\partial}

\newcommand{\s}{\vskip.1in}
\newcommand{\n}{\noindent}

\def\SFH{SFH}  
\def\CF{CF} 
\def\EH{EH} 
\def\HFhat{\widehat{\textit{HF}}\hspace{.03in}}

\newcommand{\be}{\begin{enumerate}}
\newcommand{\ee}{\end{enumerate}}

\numberwithin{equation}{subsection}

\begin{document}

\title{Contact structures, sutured Floer homology and TQFT}

\author{Ko Honda}
\address{University of Southern California, Los Angeles, CA 90089}
\email{khonda@usc.edu} \urladdr{http://rcf.usc.edu/\char126
khonda}

\author{William H. Kazez}
\address{University of Georgia, Athens, GA 30602} \email{will@math.uga.edu}
\urladdr{http://www.math.uga.edu/\char126 will}

\author{Gordana Mati\'c}
\address{University of Georgia, Athens, GA 30602} \email{gordana@math.uga.edu}
\urladdr{http://www.math.uga.edu/\char126 gordana}

\date{This version: July 13, 2008. (The pictures are in color.)}



\thanks{KH supported by an NSF
CAREER Award (DMS-0237386); GM supported by NSF grant DMS-0410066;
WHK supported by NSF grant DMS-0406158.}

\begin{abstract}
We describe the natural gluing map on sutured Floer homology which
is induced by the inclusion of one sutured manifold $(M',\Gamma')$
into a larger sutured manifold $(M,\Gamma)$, together with a contact
structure on $M-M'$.  As an application of this gluing map, we
produce a $(1+1)$-dimensional TQFT by dimensional reduction and
study its properties.
\end{abstract}

\maketitle

\section{Introduction}

Since its inception around 2001, Ozsv\'ath and Szab\'o's Heegaard
Floer homology~\cite{OS1,OS2} has been developing at a breakneck
pace. In one direction, Ozsv\'ath-Szab\'o~\cite{OS4} and,
independently, Rasmussen~\cite{Ra} defined knot invariants, called
knot Floer homology, which categorified the Alexander polynomial.
Although its initial definition was through Lagrangian Floer
homology, knot Floer homology was recently shown to admit a
completely combinatorial description by
Manolescu-Ozsv\'ath-Sarkar~\cite{MOS}. Knot Floer homology is a
powerful invariant which detects the genus of a knot by the work of
Ozsv\'ath-Szab\'o~\cite{OS6}, and detects fibered knots by the work
of Ghiggini~\cite{Gh} and Ni~\cite{Ni}. (The latter was formerly
called the ``fibered knot conjecture'' of Ozsv\'ath-Szab\'o).

One of the offshoots of the effort to prove this fibered knot
conjecture is the definition of a relative invariant for a
$3$-manifold with boundary. In a pair of important
papers~\cite{Ju1,Ju2}, Andr\'as Juh\'asz generalized the hat
versions of Ozsv\'ath and Szab\'o's Heegaard Floer
homology~\cite{OS1,OS2} and link Floer homology~\cite{OS4} theories,
and assigned a Floer homology group $\SFH(M,\Gamma)$ to a balanced
sutured manifold $(M,\Gamma)$. (A related theory is being worked out
by Lipshitz~\cite{Li1,Li2} and
Lipshitz-Ozsv\'ath-Thurston~\cite{LOT}.)

In \cite{HKM3}, the present authors defined an invariant
$EH(M,\Gamma,\xi)$ of $(M,\Gamma,\xi)$, a contact 3-manifold
$(M,\xi)$ with convex boundary and dividing set $\Gamma$ on $\bdry
M$, as an element in $\SFH(-M,-\Gamma)$. Our invariant generalized
the contact class in Heegaard Floer homology in the closed case, as
defined by Ozsv\'ath and Szab\'o~\cite{OS3} and reformulated by the
authors in~\cite{HKM2}. The definition of the contact invariant was
made possible by the work of Giroux~\cite{Gi2}, which provides
special Morse functions (called {\em convex Morse functions}) or,
equivalently, open book decompositions which are adapted to contact
structures.

Recall that a {\em sutured manifold} $(M,\Gamma)$, due to
Gabai~\cite{Ga}, is a compact, oriented, not necessarily connected
3-manifold $M$ with boundary, together with an oriented embedded
$1$-manifold $\Gamma\subset \bdry M$ which bounds a subsurface of
$\bdry M$.\footnote{This definition is slightly different from that
of Gabai~\cite{Ga}.} More precisely, there is an open subsurface
$R_+(\Gamma)\subset \bdry M$ (resp.\ $R_-(\Gamma)$) on which the
orientation agrees with (resp.\ is the opposite of) the orientation
on $\bdry M$ induced from $M$, and $\Gamma=\bdry R_+(\Gamma)=\bdry
R_-(\Gamma)$ as oriented $1$-manifolds.  A sutured manifold
$(M,\Gamma)$ is {\em balanced} if $M$ has no closed components,
$\pi_0(\Gamma)\rightarrow \pi_0(\bdry M)$ is surjective, and
$\chi(R_+(\Gamma))=\chi(R_-(\Gamma))$ on the boundary of every
component of $M$. In particular, every boundary component of $\bdry
M$ nontrivially intersects the suture $\Gamma$.

In this paper, we assume that all sutured manifolds are balanced and
all contact structures are cooriented. Although every connected
component of a balanced sutured manifold $(M,\Gamma)$ must have
nonempty boundary, our theorems are also applicable to closed,
oriented, connected $3$-manifolds $M$.  Following
Juh\'asz~\cite{Ju1}, a closed $M$ can be replaced by a balanced
sutured manifold as follows: Let $B^3$ be a $3$-ball inside $M$, and
consider $M-B^3$. On $\bdry (M-B^3)=S^2$, let $\Gamma=S^1$. Since
$\SFH(M-int(B^3),\Gamma=S^1)$ is naturally isomorphic to
$\HFhat(M)$, we can view a closed $M$ as $(M-int(B^3),S^1)$.

The goal of this paper is to understand the effect of cutting/gluing
of sutured manifolds. We first define a map which is induced from
the inclusion of one balanced sutured manifold $(M',\Gamma')$ into
another balanced sutured manifold $(M,\Gamma)$, in the presence of a
``compatible'' contact structure $\xi$ on $M-int(M')$. Here we say
that $(M',\Gamma')$ is a {\em sutured submanifold} of a sutured
manifold $(M,\Gamma)$ if $M'$ is a submanifold with boundary of $M$,
so that $M'\subset int(M)$.
If a connected component $N$ of $M-int(M')$ contains no components
of $\bdry M$  we say  that $N$ is {\em isolated}.

We will work with Floer homology groups over the ring $\Z$. With
$\Z$-coefficients, the contact invariant $\EH(M,\Gamma,\xi)$ is a
subset of $\SFH(-M,-\Gamma)$ of cardinality $1$ or $2$ of type
$\{\pm x\}$, where $x\in \SFH(-M,-\Gamma)$. (The cardinality is $1$
if and only if $x$ is a $2$-torsion element.) Over $\Z/2\Z$, the
$\pm 1$ ambiguity disappears, and $\EH(M,\Gamma,\xi)\in
\SFH(-M,-\Gamma)$.

The following theorem is the main technical result of our paper.

\begin{thm} \label{thm: gluing}
Let $(M',\Gamma')$ be a sutured submanifold of $(M,\Gamma)$, and let
$\xi$ be a contact structure on $M-int(M')$ with convex boundary and
dividing set $\Gamma$ on $\bdry M$ and $\Gamma'$ on $\bdry M'$. If
$M-int(M')$ has $m$ isolated components, then $\xi$ induces a
natural map:
$$\Phi_\xi: \SFH(-M',-\Gamma')\rightarrow \SFH(-M,-\Gamma)\otimes
V^{\otimes m},$$ which is well-defined only up to an overall $\pm$
sign.  Moreover,
$$\Phi_\xi(EH(M',\Gamma',\xi'))= EH(M,\Gamma,\xi'\cup \xi)\otimes
(x\otimes\dots\otimes x),$$ where $x$ is the contact class of the
standard tight contact structure on $S^1\times S^2$ and $\xi'$ is
any contact structure on $M'$ with boundary condition $\Gamma'$.
Here $V=\widehat{HF}(S^1\times S^2)\simeq \Z\oplus \Z$ is a
$\Z$-graded vector space where the two summands have grading which
differ by one, say $0$ and $1$.
\end{thm}

The choice of a contact structure $\xi$ on $M-M'$ plays a key role
as the ``glue'' or ``field'' which takes classes in
$\SFH(-M',-\Gamma')$ to classes in $\SFH(-M,-\Gamma)\otimes
V^{\otimes m}$. We emphasize that the gluing map $\Phi_\xi$ is
usually not injective.  The statement of the theorem, in particular
the ``naturality'' and the $V$ factor, will be explained in more
detail in Section~\ref{section:explanation}.

One immediate corollary of Theorem~\ref{thm: gluing} is the
following result, essentially proved in~\cite{HKM3}:

\begin{cor}\label{cor: inclusion}
Let $i:(M',\Gamma',\xi')\rightarrow (M,\Gamma,\xi)$ be an inclusion
such that $\xi|_{M'}=\xi'$. If $\EH(M,\Gamma,\xi)$ $\not=0$, then
$\EH(M',\Gamma', \xi')\not=0$.
\end{cor}

\s\n {\bf Gluing along convex surfaces.} Specifying a suture
$\Gamma$ on $\bdry M$ is equivalent to prescribing a
translation-invariant contact structure $\zeta_{\bdry M}$ in a
product neighborhood of $\bdry M$ with dividing set $\Gamma$. Let
$U$ be a properly embedded surface of $(M,\Gamma)$ satisfying the
following:
\begin{itemize}
\item There exists an invariant contact structure $\zeta_U$,
defined in a neighborhood of $U$, which agrees with $\zeta_{\bdry
M}$ near $\bdry U$;
\item $U$ is convex with possibly empty Legendrian boundary
and has a dividing set $\Gamma_U$ with respect to $\zeta_U$.
\end{itemize}
Let $(M',\Gamma')$ be the sutured manifold obtained by cutting
$(M,\Gamma)$ along $U$ and edge-rounding. (See \cite{H1} for a
description of edge-rounding.) By slightly shrinking $M'$, we obtain
the tight contact structure $\zeta=\zeta_{\bdry M}\cup \zeta_U$ on
$M-int(M')$.  The contact structure $\zeta$ induces the map
$$\Phi_\zeta: \SFH(-M',-\Gamma')\rightarrow \SFH(-M,-\Gamma)\otimes
V^{\otimes m},$$ for an appropriate $m$.

Summarizing, we have the following:

\begin{thm}[Gluing Map] \label{cor: gluing map}
Let $(M',\Gamma')$ be a sutured manifold and let $U_+$ and $U_-$ be
disjoint subsurfaces of $\bdry M'$ {\em(}with the orientation
induced from $\bdry M'${\em)} which satisfy the following:
\begin{enumerate}
\item Each component of $\bdry U_\pm$ transversely and nontrivially
intersects $\Gamma'$.
\item There is an orientation-reversing diffeomorphism $\phi:
U_+\rightarrow U_-$ which takes $\Gamma'|_{U_+}$ to $\Gamma'|_{U_-}$
and takes $R_\pm(U_+)$ to $R_\mp(U_-)$.
\end{enumerate}
Let $(M,\Gamma)$ be the sutured manifold obtained by gluing $U_+$
and $U_-$ via $\phi$, and smoothing.  Then there is a natural gluing
map
$$\Phi: \SFH(-M',-\Gamma')\rightarrow \SFH(-M,-\Gamma)\otimes
V^{\otimes m},$$ where $m$ equals the number of components of $U_+$
that are closed surfaces. Moreover, if $(M,\Gamma,\xi)$ is obtained
from $(M',\Gamma',\xi')$ by gluing, then
$$\Phi(\EH(M',\Gamma',\xi'))= \EH(M,\Gamma,\xi)\otimes
(x\otimes\dots\otimes x),$$ where $x$ is the contact class of the
standard tight contact structure on $S^1\times S^2$.
\end{thm}

In particular, when $\Gamma_U$ is {\em $\bdry$-parallel}, i.e., each
component of $\Gamma_U$ cuts off a half-disk which intersects no
other component of $\Gamma_U$, then the convex decomposition
$(M,\Gamma)\stackrel{(U,\Gamma_U)}\rightsquigarrow (M',\Gamma')$
corresponds to a sutured manifold decomposition by \cite{HKM1}. In
Section~\ref{section: 4.5} we indicate why our gluing map
$$\Phi:\SFH(-M',-\Gamma')\hookrightarrow \SFH(-M,-\Gamma)$$ is the
same as the direct summand map constructed in
\cite[Section~6]{HKM3}.\footnote{In \cite{Ju2}, Juh\'asz proves that
a sutured manifold gluing induces a direct summand map
$\SFH(-M',-\Gamma')\hookrightarrow \SFH(-M,-\Gamma)$. Although it is
expected that this map agrees with the natural gluing map, this has
not been proven.}

\s\n {\bf $(1+1)$-dimensional TQFT}. We now describe a
$(1+1)$-dimensional TQFT, which is obtained by dimensional reduction
of sutured Floer homology and gives an invariant of multicurves on
surfaces. (In this paper we loosely use the terminology ``TQFT''.
The precise properties satisfied by our ``TQFT'' are given in
Section~\ref{section: TQFT}.)

Let $\Sigma$ be a compact, oriented surface with nonempty boundary
$\bdry\Sigma$, and $F$ be a finite set of points of $\bdry \Sigma$,
where the restriction of $F$ to each component of $\bdry \Sigma$
consists of an even number $\geq 2$ of points. Moreover, the
connected components of $\bdry\Sigma-F$ are alternately labeled $+$
and $-$. Also let $K$ be a properly embedded, oriented
$1$-dimensional submanifold of $\Sigma$ whose boundary is $F$ and
which divides $\Sigma$ into $R_+$ and $R_-$ in a manner compatible
with the labeling of $\bdry \Sigma-F$. Let $\xi_K$ be the
$S^1$-invariant contact structure on $S^1\times\Sigma$ which traces
the dividing set $K$ on each $\{pt\}\times \Sigma$. Let $F_0\subset
\bdry \Sigma$ be obtained from $F$ by shifting slightly in the
direction of $\bdry \Sigma$. The corresponding contact invariant
$\EH(\xi_K)$ is a subset of $SFH(-(S^1\times\Sigma),-(S^1\times
F_0))$ of the form $\{\pm x\}$. The TQFT assigns to each
$(\Sigma,F)$ a graded $\Z$-module $V(\Sigma,F)=
\SFH(-(S^1\times\Sigma),-(S^1\times F_0))$ and to each $K$ the
subset $\EH(\xi_K)\subset V(\Sigma,F)$.

One application of the TQFT is the following:

\begin{thm} \label{thm: no single value}
The contact invariant in sutured Floer homology does not always
admit a single-valued representative with $\Z$-coefficients.
\end{thm}

Next, we say that $K$ is {\em isolating} if $\Sigma - K$ contains a
component that does not intersect $\bdry \Sigma$. Using the TQFT
properties we will prove:

\begin{thm} \label{thm: vanishing}
Over $\Z/2\Z$, $\EH(\xi_K) \neq 0$ if and only if $K$ is
nonisolating.
\end{thm}

Theorem~\ref{thm: vanishing}, combined with Corollary~\ref{cor:
inclusion}, expands our repertoire of universally tight contact
structures which are not embeddable in Stein fillable contact
$3$-manifolds.

\begin{cor}
Let $\xi_K$ be the $S^1$-invariant contact structure on
$S^1\times\Sigma$ corresponding to the dividing set $K \subset
\Sigma$. Then $\xi_K$ cannot be embedded in a Stein fillable $($or
strongly symplectically fillable$)$ closed contact $3$-manifold if
$K$ is isolating.
\end{cor}

Finally, we remark that $V(\Sigma,F)$ is the Grothendieck group of a
category $\mathcal{C}(\Sigma,F)$, called the {\em contact category},
whose objects are dividing sets on $(\Sigma,F)$ and whose morphisms
are contact structures on $\Sigma\times[0,1]$.  The contact category
will be treated in detail in \cite{H3}.

\s\n {\bf Organization of the paper.} In Section~\ref{section:
preliminaries} we review the notions of sutured Floer homology and
partial open book decompositions, which appeared in \cite{Ju1, Ju2,
HKM3}. Section~\ref{section:explanation} is devoted to explaining
Theorem~\ref{thm: gluing}, in particular the $V$ factor and the
naturality statement. Theorem~\ref{thm: gluing} will be proved in
Sections~\ref{section: definition} and \ref{section:
well-definition}.  The map $\Phi_\xi$ will be defined in
Section~\ref{section: definition} and the fact that $\Phi_\xi$ is a
natural map will be proved in Section~\ref{section:
well-definition}. We remark that, although the basic idea of the
definition of $\Phi_\xi$ is straightforward, the actual definition
and the proof of naturality are unfortunately rather involved. Basic
properties of the gluing map will be given in Section~\ref{section:
4.5}. Section~\ref{section: TQFT} is devoted to analyzing the
$(1+1)$-dimensional TQFT.

\section{Preliminaries} \label{section: preliminaries}

We first review some notions which appeared in \cite{Ju1,Ju2} and
\cite{HKM3}.

Let $(M,\Gamma)$ be a balanced sutured manifold. Then a {\em
Heegaard splitting $(\Sigma,\alpha,\beta)$ for $(M,\Gamma)$}
consists of a properly embedded oriented surface $\Sigma$ in $M$
with $\partial \Sigma = \Gamma$ and two sets of disjoint simple
closed curves $\alpha= \{\alpha_1, \dots, \alpha_r\}$ and $\beta=
\{\beta_1, \dots, \beta_r\}$. The Heegaard surface $\Sigma$
compresses to $R_-(\Gamma)$ along the collection $\alpha$ and to
$R_+(\Gamma)$ along the collection $\beta$. The number of $\alpha$
curves equals the number of $\beta$ curves since $(M,\Gamma)$ is
assumed to be balanced.

To define the sutured Floer Homology groups, as introduced by
Juh\'asz, we consider the Lagrangian tori $\mathbb{T}_\alpha=
\alpha_1\times\dots\times \alpha_r$ and $\mathbb{T}_{\beta}=
\beta_1\times\dots\times \beta_r$ in $Sym^r(\Sigma)$. Let
$\CF(\Sigma,\alpha,\beta)$ be the free $\Z$-module generated by the
points $\mathbf{x}=(x_1,\dots,x_r)$ in $\mathbb{T}_\alpha\cap
\mathbb{T}_\beta$.  In the definition of the boundary map for
sutured Floer homology, the suture $\Gamma$ plays the role of the
basepoint.  Denote by $\mathcal{M}_{\mathbf{x} ,\mathbf{y}}$ the
0-dimensional (after quotienting by the natural $\R$-action) moduli
space of holomorphic maps $u$ from the unit disk $D^2\subset \C$ to
$Sym^r(\Sigma)$ that (i) send $1\mapsto \mathbf{x}$, $-1\mapsto
\mathbf{y}$, $S^1 \cap \{\mbox{Im } z \geq 0 \}$ to
$\mathbb{T}_{\alpha}$ and $S^1 \cap \{\mbox{Im } z \leq 0\}$ to
$\mathbb{T}_{\beta}$, and (ii) avoid $\bdry\Sigma \times
Sym^{r-1}(\Sigma)\subset Sym^r(\Sigma)$. Then define
$$\bdry \mathbf{x} =  \sum_{\mu(\mathbf{x},\mathbf{y})=1} ~~
\#(\mathcal{M}_{\mathbf{x}, \mathbf{y} })~~  \mathbf{y},$$ where
$\mu(\mathbf{x},\mathbf{y})$ is the relative Maslov index of the
pair and $\#(\mathcal{M}_{\mathbf{x}, \mathbf{y} })$ is a signed
count of points in $\mathcal{M}_{\mathbf{x}, \mathbf{y} }$. The
homology of $\CF(\Sigma,\alpha,\beta)$ is the sutured Floer homology
group $SFH(\Sigma,\alpha,\beta)=SFH(M,\Gamma)$.

In \cite{HKM3}, the present authors defined an invariant
$EH(M,\Gamma,\xi)$ of $(M,\Gamma,\xi)$, a contact 3-manifold with
convex boundary and dividing set $\Gamma$ on $\bdry M$, as an
element in $\SFH(-M,-\Gamma)$. This invariant generalizes the
contact class in Heegaard Floer homology in the closed case, as
defined by Ozsv\'ath and Szab\'o~\cite{OS3}, and described from a
different point of view in ~\cite{HKM2}. For the sake of
completeness, we sketch the definition of the invariant
$EH(M,\Gamma,\xi)$.

First consider the case when $\xi$ is a contact structure on a
closed manifold $M$.  In~\cite{HKM2} we used an open book
decomposition compatible with $\xi$ to construct a convenient
Heegaard decomposition $(\Sigma,\alpha,\beta)$ for $M$ in which the
contact class was a distinguished element in
$\widehat{HF}(\Sigma,\alpha,\beta)$. Recall that an open book
decomposition for $M$ is a pair $(S,h)$ consisting of a surface $S$
with boundary and a homeomorphism $h:S\stackrel\sim\rightarrow S$
with $h|_{\bdry S}=id$, so that $M \simeq S\times[0,1]/\sim_{h}$,
where $(x,1) {\sim} _h (h(x),0)$ for $x \in S$ and
$(x,t){\sim}_h(x,t')$ for $x \in \bdry S$, $t,t'\in[0,1]$. A
Heegaard decomposition $(\Sigma,\beta,\alpha)$ for $-M$ (recall that
the contact class lives in the Heegaard Floer homology of $-M$) is
obtained from the two handlebodies $H_1=S\times[0,{1\over
2}]/\sim_h$ and $H_2=S\times[{1\over 2},1]/\sim_h$, which are glued
along the common boundary $\Sigma=(S \times \{{1\over 2} \}) \cup
-(S\times \{0 \})$ by $id \cup h$. Take a family of properly
embedded disjoint arcs $a_i$ that cuts the surface $S$ into a disk,
and small push-offs $b_i$ of $a_i$ (in the direction of the
boundary) such that $b_i$ intersects $a_i$ in exactly one point. The
compressing disks for $H_1$ and $H_2$, respectively, are
$D_{a_i}=a_i\times[0,{1\over 2}]$ and $D_{b_i}=b_i\times[{1\over
2},1]$; set $\alpha_i=\bdry D_{a_i}$ and $\beta_i=\bdry D_{b_i}$. We
call the family of arcs $a_i$ a {\em basis} for $S$, and show
in~\cite{HKM2} that the element of Heegaard Floer homology that
corresponds to the generator $\mathbf{x}=(x_1,\dots , x_n)$, where
$x_i$ is the unique intersection point of $a_i \times \{{1\over 2}
\}$ and $b_i \times \{{1\over 2} \}$,  is independent of the choice
of basis for $S$ and the compatible open book decomposition.
Moreover, it is the contact class defined by Ozsv\'ath and Szab\'o.

To define the contact class $EH(M,\Gamma,\xi)$ in the case of a
balanced sutured manifold, we generalize the notions of an ``open
book'' and a ``basis'', involved in the definition of the contact
invariant above.  Let $(A,B)$ be a pair consisting of a surface $A$
with nonempty boundary and a subsurface $B\subset A$. A collection
$\{a_1,\dots,a_k\}$ of properly embedded disjoint arcs in $A$ is
called a {\em basis for $(A,B)$} if each $a_i$ is disjoint from $B$
and $A-\cup_{i=1}^k a_i$ deformation retracts to $B$. A {\em partial
open book} $(S,R_+(\Gamma),h)$ consists of the following data: a
compact, oriented surface $S$ with nonempty boundary, a subsurface
$R_+(\Gamma)\subset S$, and a ``partial'' monodromy map $h:
P\rightarrow S,$ where $P\subset S$ is the closure of
$S-\overline{R_+(\Gamma)}$ and $h(x)=x$ for all $x\in (\bdry S)\cap
P$. We say that $(S,R_+(\Gamma),h)$ is a {\em partial open book
decomposition for} $(M,\Gamma)$ if $M \simeq S\times[0,1]/\sim_{h}$,
where the equivalence relation is $(x,1) {\sim} _h (h(x),0)$ for $x
\in P$ and $(x,t){\sim}_h(x,t')$ for $x \in \bdry S$,
$t,t'\in[0,1]$. Since the monodromy $h$ is defined only on $P$, the
space obtained after gluing has boundary consisting of  $R_+ \times
\{1\}$ and $R_-\times \{0\}$, where $R_- = S- h(P)$. The suture
$\Gamma$ is the common boundary of $R_+ \times\{1\}$ and $R_-\times
\{0\}$.

To see a handlebody decomposition of $M$ from this point of view,
let $H_1=S\times[0,{1\over 2}]/\sim$, where $(x,t)\sim (x,t')$ if $x
\in \bdry S$ and $t,t' \in [0,{1\over 2}]$ and let $H_2=
P\times[{1\over 2},1]/\sim$, where $(x,t)\sim (x,t')$ if $x\in \bdry
P$ and $t,t'\in [{1\over 2},1]$. It is clear that we can think of $M
\simeq S\times[0,1]/\sim_h$ as $M \simeq H_1\cup H_2/$gluing, where
the handlebodies are glued along portions of their boundary as
follows: $(x,{1\over 2}) \in H_1$ is identified to $(x,{1\over 2})
\in H_2$ and $(x,1) \in H_2$ is identified with $(h(x),0) \in H_1$
for $x\in P$. This leaves $R_+ \times \{{1\over 2} \}$ and $R_-
\times\{0\}$ as the boundary of the identification space. Now let
$\{a_1,\dots,a_k\}$ be a {\em basis} for $(S,R_+(\Gamma))$ in the
sense defined above. Let $b_i$, $i=1,\dots,k$, be pushoffs of $a_i$
in the direction of $\bdry S$ so that $a_i$ and $b_i$ intersect
exactly once. Then it is not hard to see that if we set
$\Sigma=(S\times\{0\})\cup (P\times\{{1\over 2}\})$, $\alpha_i=\bdry
(a_i\times[0,{1\over 2}])$ and $\beta_i= (b_i\times\{{1\over
2}\})\cup(h(b_i)\times\{0\})$, then $(\Sigma,\beta,\alpha)$ is a
Heegaard diagram for $(-M,-\Gamma)$.

The two handlebodies $H_1$ and $H_2$ defined above by the open book
decomposition $(S,R_+(\Gamma),h)$ carry unique product disk
decomposable contact structures. After gluing, they determine a
contact structure $\xi_{(S,R_+(\Gamma),h)}$ on $(M,\Gamma)$. We say
that a partial open book decomposition $(S,R_+(\Gamma),h)$ and a
contact structure $\xi$ are {\em compatible} if
$\xi=\xi_{(S,R_+(\Gamma),h)}$. On the other hand, as in the closed
manifold case, every contact structure $\xi$ with convex boundary on
a sutured manifold $(M,\Gamma)$ gives rise to a compatible partial
open book decomposition $(S,R_+(\Gamma),h)$.


\section{Explanation of Theorem~\ref{thm: gluing}}\label{section:explanation}

\subsection{Naturality}
We now explain what we mean by a ``natural map'' $\Phi_\xi$.  Recall
the following theorem from Ozsv\'ath-Szab\'o~\cite{OS7}:

\begin{thm}[Ozsv\'ath-Szab\'o] \label{thm: OS}
Given two Heegaard decompositions $(\Sigma,\alpha,\beta)$,
$(\overline\Sigma,\overline\alpha,\overline\beta)$ of a closed
3-manifold $M$, the isomorphism
$$\Psi: \HFhat(\Sigma,\alpha,\beta)\stackrel\sim\rightarrow
\HFhat(\overline\Sigma,\overline\alpha,\overline\beta),$$ given as
the composition of stabilization/destabilization, handleslide, and
isotopy maps, is well-defined up to an overall factor of $\pm 1$ and
does not depend on the particular sequence chosen from
$(\Sigma,\alpha,\beta)$ to
$(\overline\Sigma,\overline\alpha,\overline\beta)$.
\end{thm}

This lack of monodromy allows us to ``naturally'' identify the
isomorphic Heegaard Floer homology groups
$\HFhat(\Sigma,\alpha,\beta)$ and
$\HFhat(\overline\Sigma,\overline\alpha,\overline\beta)$, up to an
overall sign. Sutured Floer homology enjoys the same naturality
property, that is, the isomorphism
$$\Psi: \SFH(\Sigma,\alpha,\beta)\stackrel\sim\rightarrow \SFH(\overline\Sigma,
\overline\alpha,\overline\beta)$$ is also well-defined up to an
overall factor of $\pm 1$ and is independent of the same type of
choices if $(\Sigma,\alpha,\beta)$,
$(\overline\Sigma,\overline\alpha,\overline\beta)$ are two Heegaard
decompositions for $(M,\Gamma)$.

Next, suppose $(\Sigma',\beta',\alpha')$,
$(\overline{\Sigma}',\overline \beta',\overline\alpha')$ are
Heegaard splittings for $(-M',-\Gamma')$ and
$(\Sigma,\beta,\alpha)$,
$(\overline\Sigma,\overline\beta,\overline\alpha)$ are their
extensions to $(-M,-\Gamma)$. We will restrict ourselves to working
with a certain subclass of Heegaard splittings of $(-M',-\Gamma')$,
namely those that are contact-compatible on a neighborhood of $\bdry
M'$, with respect to an invariant contact structure $\zeta$ which
induces the dividing set $\Gamma'$ on $\bdry M'$.   The Heegaard
splittings for $(-M,-\Gamma)$ we will use extend those of
$(-M',-\Gamma')$ and are contact-compatible with respect to $\xi$ on
$M-M'$. Assume $M-int(M')$ has no isolated components.  Then, in the
statement of Theorem~\ref{thm: gluing}, we take the commutativity of
the following diagram to be the definition of the {\em naturality}
of $\Phi_\xi$:
\begin{equation} \label{commutative}
\begin{diagram}
\SFH(\Sigma',\beta',\alpha') & \rTo^{(\Phi_\xi)_1} &
\SFH(\Sigma,\beta,\alpha)  \\
\dTo^{\Psi_1} & &  \dTo^{\Psi_2} \\
\SFH(\overline{\Sigma}',\overline\beta',\overline\alpha')  &
\rTo^{(\Phi_\xi)_2} &
\SFH(\overline\Sigma,\overline\beta,\overline\alpha)  \\
\end{diagram}
\end{equation}
Here the vertical maps $\Psi_1,\Psi_2$ are the natural isomorphisms
of Theorem~\ref{thm: OS} and $(\Phi_\xi)_1, (\Phi_\xi)_2$ are the
maps induced by $\xi$, to be defined in Section~\ref{section:
definition}.

\subsection{Explanation of the $V$ factor.}
Consider $(M,\Gamma,\xi)$ and a compatible partial open book
decomposition $(S,R_+(\Gamma),h)$.  Let $\{a_1,\dots,a_k\}$ be a
basis for $(S,R_+(\Gamma))$.  Consider a larger collection $\{a_1,
\dots, a_k, a_{k+1},\dots,a_{k+l}\}$ of properly embedded disjoint
arcs in $S$ which satisfy $a_i\subset P$, so that
$S-\cup_{i=1}^{k+l} a_i$ is a disjoint union of disks $D_j$,
$j=1,\dots,l$, and a surface that deformation retracts to
$R_+(\Gamma)$. For each $j$, pick $z_j\in D_j$ and consider a small
neighborhood $N(z_j)\subset D_j$. Then $\{a_1,\dots,a_{k+l}\}$
becomes a basis for $(S,R_+(\Gamma)\cup (\cup_{j=1}^l N(z_j)))$. The
Heegaard surface for $(S,R_+(\Gamma))$ is $\Sigma=(P\times\{1\})\cup
(S\times\{0\})$, whereas the Heegaard surface for $
(S,R_+(\Gamma)\cup (\cup_{j=1}^{l} N(z_j)))$ is $\Sigma'=\Sigma-
\cup_{j=1}^l N(z_j)$.  As in Section~\ref{section: preliminaries},
the $a_i$ determine arcs $b_i$ as well as closed curves $\alpha_i$,
$\beta_i$. We refer to the procedure of adding extra arcs to a basis
for $(S,R_+(\Gamma))$ and extra $N(z_j)$'s to $R_+(\Gamma)$ as
``placing extra dots'' or ``placing extra $z_j$'s''.

\begin{claim}
The effect of placing an extra dot on $(S,R_+(\Gamma))$ on sutured
Floer homology is that of taking the tensor product with
$\HFhat(S^1\times S^2)\simeq V$.
\end{claim}

\begin{proof}
Consider the following situation: Suppose $\{a_1,\dots,a_k\}$ is a
basis for $(S,R_+(\Gamma))$. Then add an extra properly embedded arc
$a_{k+1}\subset P$ of $S$ which is disjoint from $a_1,\dots,a_k$ and
such that one component $D_1$ of $S-a_{k+1}$ is a half-disk which is
contained in $P$.  Also add an extra dot $z_1$ in the component
$D_1$. The $\alpha_{k+1}$ and $\beta_{k+1}$ corresponding to
$a_{k+1}$ intersect in exactly two points, and do not interact with
the other $\alpha_i$ and $\beta_i$.
By the placement of the extra dot,
$$\SFH(\Sigma',\{\beta_1,\dots,\beta_{k+1}\},\{\alpha_1,\dots,\alpha_{k+1}\})$$
$$\simeq
\SFH(\Sigma,\{\beta_1,\dots,\beta_{k}\},\{\alpha_1,\dots,\alpha_{k}\})\otimes
\HFhat(S^1\times S^2).$$

Next, after a sequence of arc slides as in Section~3.1 of
\cite{HKM2} (or just handleslides), we can pass between any two
bases of $(S,R_+(\Gamma)\cup N(z_1))$, where $N(z_1)$ is a small
disk in $S-R_+(\Gamma)$ about $z_1$.  Since $\SFH$ is invariant
under any sequence of handleslides, the claim follows.
\end{proof}

Now we explain the $V$ factors that appear in Theorem~\ref{thm:
gluing}. Let us consider a partial open book decomposition $(S',h')$
for any contact structure $\xi'$ which is compatible with
$(M',\Gamma')$, and let $(S,h)$ be a partial open book decomposition
for $\xi\cup\xi'$ which extends $(S',h')$. If no connected component
of $M-int(M')$ is isolated, then a basis $\{a_1',\dots,a_k'\}$ for
$(S',h')$ easily extends to a basis for $(S,h)$. If there are $m$
isolated components of $M-int(M')$, then $S-\cup_{i=1}^k a_i'$ has
$m$ connected components which do not intersect $R_+(\Gamma)$ (and
hence can never be completed to a basis for $(S,h)$).  Instead, by
adding $m$ extra dots $z_1,\dots,z_m$, we can extend
$\{a_1',\dots,a_k'\}$ to a basis for $(S, R_+(\Gamma)\cup
(\cup_{i=1}^m N(z_i)))$.

\section{Definition of the map $\Phi_\xi$} \label{section:
definition}

In this section we define the chain map:
$$\Phi_\xi:\CF(-M',-\Gamma')\rightarrow \CF(-M,-\Gamma),$$
which induces the map, also called $\Phi_\xi$ by slight abuse of
notation, on the level of homology. Let us assume that $M-int(M')$
has no isolated components. The general case follows without
additional effort, by putting extra dots.

\begin{proof}[Sketch of the construction.]
We start by giving a quick overview of the construction of
$\Phi_\xi$.  The actual definition needed to prove naturality is
considerably more complicated and occupies the remainder of the
section.

Let us first decompose $M=M'\cup M''$, where $M''=M-int(M')$. Let
$\Sigma'$ be a Heegaard surface for the sutured manifold
$(M',\Gamma')$.  By definition, $\Sigma' \cap
\partial M' =\Gamma'$. Next choose compressing disks $\alpha',
\beta'$ on $\Sigma'$.  Also let $\Sigma''$ be a Heegaard surface for
the sutured manifold $(M'',\Gamma''\cup -\Gamma')$.

Although it might appear natural to take the union of $\Sigma'$ and
$\Sigma''$ along their common boundary $\Gamma'$ to create a
Heegaard surface for $M$, we are presented with a problem. If we
glue $M'$ and $M''$ to obtain $M$, then, on the common boundary of
$M'$ and $M''$, $R_\pm(\Gamma')$ from $M'$ is glued to
$R_\mp(-\Gamma')$ from $M''$. As a result, the $\alpha$-curves for
$\Sigma'$ and $\beta$-curves for $\Sigma''$ will be paired, and the
$\beta$-curves for $\Sigma'$ and $\alpha$-curves for $\Sigma''$ will
be paired, and we will be mixing homology and cohomology. A way
around this problem is to insert the layer $N=T\times[0,1]$, where
$T=\bdry M'$, so that $M=M'\cup N\cup M''$, $M'\cap N=T\times\{0\}$,
and $M''\cap N=T\times\{1\}$. Let $\Sigma_N$ be a Heegaard surface
for $(N,(-\Gamma'\times\{0\})\cup (\Gamma'\times\{1\}))$. Then
$\Sigma = \Sigma' \cup \Sigma_N \cup \Sigma''$ is a Heegaard surface
for $M$.

A second issue which arises is that the union of compressing disks
for $M'$, $N$, and $M''$ is not sufficient to give a full set of
compressing disks for $M$.  Our remedy is to use the contact
invariant: First we take $(\Sigma',\beta',\alpha')$ to be {\em
contact-compatible} near $\bdry M'$.  Roughly speaking, this means
that $(\Sigma',\beta',\alpha')$, near $\bdry M'$, looks like a
Heegaard decomposition arising from a partial open book
decomposition of a contact structure $\zeta$ which is defined near
$\bdry M'$ and has dividing set $\Gamma'$ on $\bdry M'$.  Let
$\Sigma''$ be a Heegaard surface which is compatible with
$\xi|_{M''}$ and let $\Sigma_N$ be a Heegaard surface which is
compatible with the $[0,1]$-invariant contact structure $\xi|_N$. We
then extend $\alpha'$ and $\beta'$ by adding $\alpha''$ and
$\beta''$ which are compatible with $\xi\cup \zeta$, and then define
$$\Phi_\xi: \CF(\Sigma',\beta',\alpha')\rightarrow \CF(\Sigma,\beta'\cup \beta'',
\alpha'\cup\alpha''),$$
$$\mathbf{y}\mapsto (\mathbf{y}, \mathbf{x}''),$$
where $\bf x''$ is the contact class $EH(\xi\cup \zeta)$, consisting
of a point from each $\beta_i''\cap \alpha_i''$.
\end{proof}

We now give precise definitions.  Let $T=\bdry M'$ and let
$T\times[-1,1]$ be a neighborhood of $T=T\times\{0\}$ with a
$[-1,1]$-invariant contact structure $\zeta$ which satisfies the
following:
\begin{itemize}
\item $T_t=T\times\{t\}$, $t\in[-1,1]$, are convex surfaces with
dividing set $\Gamma'\times\{t\}$;
\item $T\times[-1,0]\subset M'$ and $T\times [0,1]\subset M-int(M')$;
\item $\xi|_{T\times[0,1]}=\zeta|_{T\times[0,1]}$.
\end{itemize}

In order to define $\Phi_\xi$, we need to construct a suitable
Heegaard splitting $(\Sigma',\beta',\alpha')$ for the sutured
manifold $(-M',-\Gamma')$ and a contact-compatible extension to
$(\Sigma,\beta,\alpha)$ for $(-M,-\Gamma)$.  This will be done in
several steps.

\s\n {\bf Step 1: Construction of $(\Sigma',\beta',\alpha')$.} In
this step we construct $(\Sigma',\beta',\alpha')$ which is {\em
contact-compatible with respect to $\zeta$ near $\bdry M'$}.
(Although a little unwieldy, we take the construction below as the
{\em definition} of a contact-compatible $(\Sigma',\beta',\alpha')$
with respect to $\zeta$ near $\bdry M'$.) The technique is similar
to the proof of Theorem~1.1 of \cite{HKM3}.

\begin{figure}[ht]
\begin{overpic}[width=10cm]{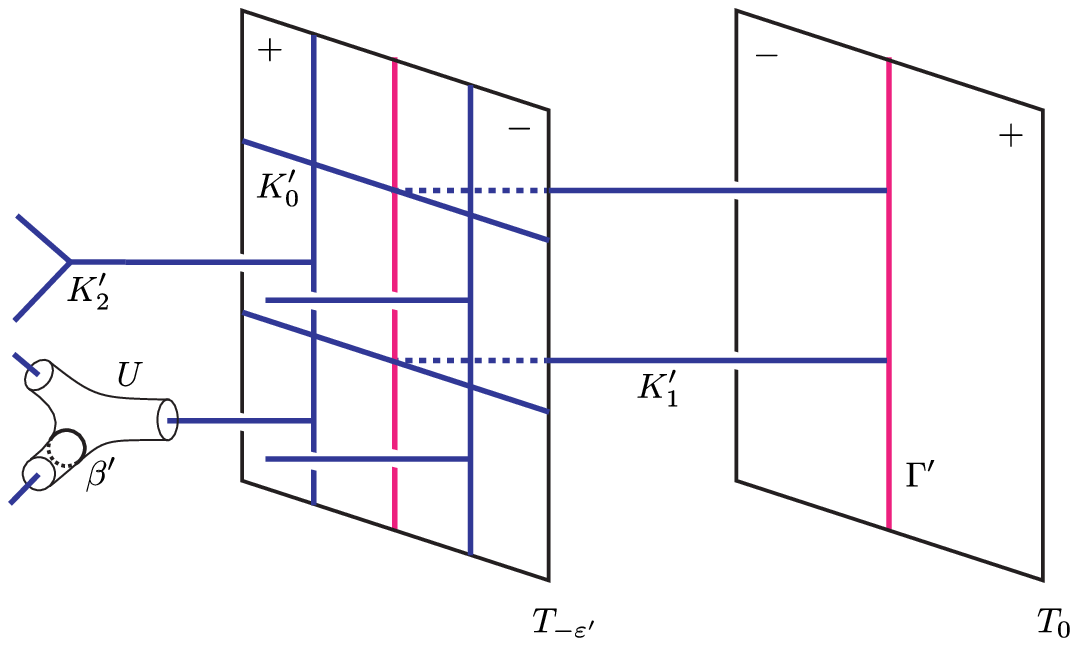}
\end{overpic}
\caption{} \label{step1}
\end{figure}

Let $0<\varepsilon' <1$. Start by choosing a cellular decomposition
of $T_{-\varepsilon'}$ so that the following hold:
\begin{itemize}
\item The $1$-skeleton $K'_0$ is Legendrian;
\item Each edge of the cellular decomposition lies on the boundary of
two distinct $2$-cells $\Delta$, $\Delta'$;
\item The boundary of each 2-cell $\Delta$ intersects the dividing set
$\Gamma_{T_{-\varepsilon'}}$ exactly twice.
\end{itemize}
Here we use the Legendrian realization principle and isotop
$T_{-\varepsilon'}$, if necessary.  Let $K'_1$ be a finite
collection of Legendrian segments $\{p\}\times[-\varepsilon',0]$, so
that every endpoint $(p,-\varepsilon')$ in $T_{-\varepsilon'}$ lies
in $K'_0\cap (\Gamma'\times\{-\varepsilon'\})$ and for each
connected component $\gamma$ of $\Gamma'$ there are at least two
$p$'s in $\gamma$. Now let $K'_2$ be a graph attached to $K'_0$ so
that $K'_0\cup K'_2$ is a $1$-skeleton of a cellular decomposition
of $M'-(T\times(-\varepsilon',0])$ and $int(K'_2)\subset M'-(T\times
[-\varepsilon',0])$.  The graph $K'_2$ is obtained without reference
to any contact structure. If we set $K'=K'_0\cup K'_1\cup K'_2$,
then $\bdry N(K')$ is the union of the tubular portion $U$ and small
disks $D_1,\dots, D_s\subset \bdry M'$. Here $N(G)$ denotes the
tubular neighborhood of a graph $G$. See Figure~\ref{step1}.

Define the Heegaard surface $\Sigma'$ to be (a surface isotopic to)
the union $(R_-(\Gamma')-\cup_i D_i) \cup U$.  Also modify
$R_+(\Gamma')$ slightly so that $R_+(\Gamma')-\cup_i D_i$ is the new
$R_+(\Gamma')$. The $\beta'$-curves are meridians ($=$ boundaries of
compressing disks) of $N(K')$, and the $\alpha'$-curves are
meridians of the complement, as chosen below.

After a contact isotopy, we may take the standard contact
neighborhood $N(K'_0)$ to be $T\times[-{3\varepsilon'\over 2},
-{\varepsilon'\over 2}]$ with standard neighborhoods of Legendrian
arcs of type $\{q\}\times[-{3\varepsilon'\over
2},-{\varepsilon'\over 2}]$, $q\in\Gamma'$, removed. Now define the
following decomposition of $T\times[-{3\varepsilon'\over 2},0]$ into
two handlebodies:
$$H_1'=(T\times[-3\varepsilon'/ 2,0])-N(K'_0\cup K'_1),$$
$$H_2'=N(K'_0\cup K'_1),$$
where $N(K'_0\cup K'_1)$ denotes the standard contact neighborhood.
Both $H_1'$ and $H_2'$ are product disk decomposable. (The product
disk decomposability of $H_2'$ is clear.  As for $H_1'$, observe
that $(T\times [-{\varepsilon'\over 2},0])-N(K_1')$ is product disk
decomposable.) Hence we may write $H_1'= S'\times[0,1]/\sim$,
$(x,t)\sim (x,t')$ if $x\in \bdry S'$ and $t,t'\in [0,1]$. Here
$\bdry S'\times[0,1]/\sim$ is the dividing set of $\bdry H_1'$, and
$R_+(\Gamma')\subset S'\times\{1\}$. Let
$P'=S'-\overline{R_+(\Gamma')}$.  Similarly we can write
$H_2'=S'_{(2)}\times[0,1]/\sim$.  If $-\Gamma'_{-3\varepsilon'/2}$
is the dividing set of $T_{-3\varepsilon'/2}$, with the outward
orientation induced from $T\times[-{3\varepsilon'\over 2},0]$, then
let $P'_{(2)}=S'_{(2)}-\overline{R_-(-\Gamma'_{-3\varepsilon'/2})}$.
Observe that $P'\times\{1\}$ is identified with
$P'_{(2)}\times\{0\}$; let $\psi: P'\stackrel\sim\rightarrow
P'_{(2)}$ be the corresponding identification map.  Also let $h':
Q'\rightarrow S'$ be the monodromy map for the partially defined
open book, where the domain of definition $Q'$ is a subset of $P'$
and contains arcs that correspond to compressing disks of $N(K_1')$.

Next let $\{a_1',\dots,a_k'\}$ be a maximal set of properly embedded
arcs on $S'$ such that the corresponding
$\alpha_i'=\bdry(a_i'\times[0,1])$, $i=1,\dots,k$, on $\bdry H_1'$
form a maximal collection of curves which can be extended to a full
$\alpha'$ set. (We will abuse notation and call such a maximal
collection of arcs a {\em basis} for $(S',R_+(\Gamma'))$.  Let
$b_i'$ be the usual pushoff of $a_i'$, and define
$\beta_i'=\bdry(\psi(b_i')\times[0,1])$, $i=1,\dots,k$, on $\bdry
H_2'$.

\begin{lemma}
$\{\alpha_1',\dots,\alpha_k'\}$ and $\{\beta_1',\dots,\beta_k'\}$
can be completed to full $\alpha'$ and $\beta'$ sets which are
weakly admissible.
\end{lemma}

\begin{proof}
The decomposition $T\times[-{3\varepsilon'\over 2},0]=H_1'\cup H_2'$
can be extended to a decomposition of $M'$ into two handlebodies. To
accomplish this, let $K'_2$ be the graph defined above, and choose
$K_3'$ to be the graph such that $N(K_2')\cup N(K_3')$ is a
decomposition of $M'-(T\times[-{3\varepsilon'\over 2},0])$ into two
handlebodies. Then $M'$ is the union of the two handlebodies
$H_1'\cup N(K'_3)=M'-N(K')$ and $H_2'\cup N(K'_2)=N(K')$.  The
collection $\{\alpha_1',\dots,\alpha_k'\}$ can be completed to a
full $\alpha'$ set by adding $\alpha_{k+1}',\dots,\alpha_{k+l}'$
which are meridians of $N(K'_3)$. On the other hand,
$\{\beta_1',\dots,\beta_k'\}$ can be completed by adding meridians
of $N(K'_2)$, in addition to $\bdry (c_i'\times[0,1])\subset \bdry
H_2'$, where $c_i'$ are properly embedded arcs of
$R_-(-\Gamma'_{-3\varepsilon'/2})$. (Add enough compressing disks of
$N(K'_2)$ so that $H_2'\cup N(K'_2)$ compresses to $H_2'$.  Then add
enough arcs $c_i$ so that $S'_{(2)}-\cup_i c_i'-\cup_i \psi(b_i')$
deformation retracts to the ``ends'' of $S'_{(2)}$, namely the arcs
of intersection with $T_0$.)

We now prove that the above extension can be made weakly admissible,
without modifying $\alpha_i'$ and $\beta_i'$, $i=1,\dots,k$. If a
periodic domain uses any $\alpha_i'$ or $\beta_i'$ with $1\leq i\leq
k$, then the position of $R_+(\Gamma')$ and the relative positions
of $a_i'$ and $b_i'$ imply that the periodic domain has both
positive and negative signs. Hence assume that we are not using
$\alpha_i'$ or $\beta_i'$ with $1\leq i\leq k$.  It is easy to find
disjoint closed curves $\gamma_{k+1}',\dots,\gamma_{k+l}'$ which are
duals of $\alpha_{k+1}',\dots,\alpha_{k+l}'$, i.e., $\gamma_i'$ and
$\alpha_j'$ have geometric intersection number $\delta_{ij}$, and
which do not enter $\bdry H_1'$. (Hence the $\gamma_i'$ do not
intersect $\alpha_j'$ with $j=1,\dots, k$.) If we wind the
$\alpha_i'$, $i=k+1,\dots,k+l$, about the curves $\gamma_i'$ as in
\cite[Section~5]{OS1}, then the result will be weakly admissible.
\end{proof}

\n {\bf Remark.} An alternate way of thinking of the contact
compatibility with respect to $\zeta$ near $\bdry M'$ is as follows:
Start with any Heegaard decomposition $(\Sigma',\beta',\alpha')$ for
$(-M',-\Gamma')$.  Take $T\times[0,1]$ with the invariant contact
structure $\zeta$, and form a partial open book decomposition for
$(T\times[0,1],\zeta)$ by choosing a Legendrian skeleton consisting
of sufficiently many arcs of type $\{p\}\times[0,1]$, where $p\in
\Gamma'$. Let $\Sigma'_\zeta$ be the corresponding
contact-compatible Heegaard surface.  Then the Heegaard surface for
$(-M',-\Gamma')$ which is contact-compatible with respect to $\zeta$
near $\bdry M'$ is obtained from $\Sigma'$ by attaching two Heegaard
surfaces of type $\Sigma'_\zeta$, one for
$T\times[-{3\varepsilon'\over 2},-{\varepsilon'\over 2}]$ and
another for $T\times [-{\varepsilon'\over 2},0]$.  (Note that the
choice of arcs of type $\{p\}\times[0,1]$ for the two
$\Sigma'_\zeta$'s may be different.) In other words, we are gluing
two copies of $(T\times[0,1],\zeta)$ to $(-M',-\Gamma')$.

\s\n{\bf Remark.} Another approach is to restrict attention to the
class of contact-compatible Heegaard splittings for an arbitrarily
chosen, tight or overtwisted, contact $(M',\Gamma',\xi')$ compatible
with the dividing set. Suppose we show that the definition of
$\Phi_\xi$ depends only on the partial open book $(S',h')$ for
$(M',\Gamma',\xi')$, up to positive and negative stabilizations.  By
the result of \cite{GG}, two open books become isotopic after a
sequence of positive and negative stabilizations, provided they
correspond to homologous contact structures.  This would show that
$\Phi_\xi$ is only dependent on the homology class of $\xi'$.
However, we would still need to remove the dependence on the
homology class.

\s\n {\bf Step 2: Extension of the Heegaard splitting to
$(-M,-\Gamma)$.} We extend the Heegaard splitting
$(\Sigma',\beta',\alpha')$ constructed in Step 1 to a Heegaard
splitting $(\Sigma,\beta,\alpha)$ for $(-M,-\Gamma)$ which is {\em
contact-compatible with respect to $\xi\cup\zeta$}. (Again, we take
the construction below as the definition of contact-compatibility.)

Let $\varepsilon''>0$.  Then we write $M= M'\cup N\cup M''$, where
$N=N_{\varepsilon''}=T\times[0,\varepsilon'']$ and
$M''=M''_{\varepsilon''}= M-int(M'\cup N)$.

The contact manifold $(M'',\xi|_{M''})$ admits a Legendrian graph
$K''$ with endpoints on $\Gamma_{\bdry M''}$ and a decomposition
into $N(K'')$ and $M''-N(K'')$, according to
\cite[Theorem~1.1]{HKM3}.  Assume that every connected component of
$K''$ intersects $\Gamma$ at least twice. (This will be useful in
Lemma~\ref{lemma:typeAorB}.) Similarly, $(N,\xi|_N)$ admits a
Legendrian graph $K'''$ consisting of Legendrian segments
$\{q\}\times[0,\varepsilon'']$, where there is at least one $q$ for
each component of the dividing set of $T_0$. We also assume that the
endpoints of $K'$, $K''$, and $K'''$ do not intersect.

We then decompose $M$ into $H_1=(M'-N(K'))\cup N(K''') \cup
(M''-N(K''))$ and $H_2=N(K')\cup (N-N(K'''))\cup N(K'')$,
respectively. Since $N-N(K''')$ and $N(K'')$ are product disk
decomposable with respect to $\xi$, their union is also product disk
decomposable.  Hence, we have:
\begin{itemize}
\item $H_2$ is a neighborhood of a graph $K$,
\item the restriction of $K$ to
$M''\cup (T\times[-\varepsilon',\varepsilon''])$ is Legendrian, and
\item restricted to $M''\cup (T\times[-{3\varepsilon'\over 2},
\varepsilon''])$, $H_2$ is a standard contact neighborhood of
$K\cap(M''\cup (T\times[-\varepsilon',\varepsilon'']))$.
\end{itemize}
Similarly, $((T\times[-{3\varepsilon'\over 2},0]) - N(K'_0\cup
K'_1)) \cup N(K''') \cup (M''-N(K''))$ is product disk decomposable
with respect to $\xi\cup\zeta$.  Therefore, $H_1$ extends $H_1'\cup
N(K'_3)= (S'\times[0,1]/\sim) \cup N(K'_3)$ so that $H_1=
(S\times[0,1]/\sim)\cup N(K'_3)$, where $(x,t)\sim(x,t')$ if $x\in
\bdry S$ and $t,t'\in[0,1]$.  Here, $R_+(\Gamma)\subset
S\times\{1\}$ and $S'$ is a subsurface of $S$.

Therefore, we may extend $(\Sigma',\beta',\alpha')$ to
$(\Sigma,\beta,\alpha)$ as follows: Consider a collection of arcs
$a_1'',\dots, a_m''$ which form a basis for $(S-P', R_+(\Gamma))$.
Then let $\alpha''_i=\bdry (a_i''\times[0,1])$, and $\beta''_i$ be
the corresponding closed curves derived from the pushoffs $b_i''$ of
$a_i''$.  The monodromy $h$ for $b_i''$ can be computed from the
partial open book decomposition on $M''\cup
(T\times[-{3\varepsilon'\over 2},\varepsilon''])$.  Then
$\alpha=\alpha'\cup\alpha''$ and $\beta=\beta'\cup \beta''$, where
$\alpha''$ (resp.\ $\beta''$) is the collection of the $\alpha_i''$
(resp.\ $\beta_i''$). The contact-compatibility on $\Sigma-\Sigma'$
immediately implies that the extension is weakly admissible.

\s We are now in a position to define the chain map $\Phi_\xi$.  Let
$(\Sigma',\beta',\alpha')$ be a Heegaard splitting for
$(-M',-\Gamma')$ which is contact-compatible near $\bdry M'$, and
let $(\Sigma,\beta,\alpha)$ be a contact-compatible extension of
$(\Sigma',\beta',\alpha')$ to $(-M,-\Gamma)$. Now let $x_i''$ be the
preferred intersection point (i.e., the only one on $S\times\{1\}$)
between $\alpha_i''$ and $\beta_i''$, and denote their collection by
$\mathbf{x}''$. Given $\mathbf{y}\in \CF(\Sigma',\beta',\alpha')$,
we define the map:
$$\Phi_\xi: \CF(\Sigma',\beta',\alpha')\rightarrow \CF(\Sigma,\beta'\cup \beta'',
\alpha'\cup\alpha''),$$
$$\mathbf{y}\mapsto (\mathbf{y}, \mathbf{x}'').$$
The fact that $\Phi_\xi$ is a chain map follows from observing that
every nonconstant holomorphic map which emanates from $x''_i$ must
nontrivially intersect $R_+(\Gamma)$.  Hence $\mathbf{x}''$ will be
used up, and the only holomorphic maps from
$(\mathbf{y},\mathbf{x}'')$ to $(\mathbf{y}',\mathbf{x}'')$ are
holomorphic maps from $\mathbf{y}$ to $\mathbf{y}'$ within
$\Sigma'$.  The tuple $\mathbf{x}''$ will be called {\em the $\EH$
class on $S-P'$}.  It is immediate from the definition of $\Phi_\xi$
that  when  $(M',\Gamma',\xi')$ is contact,
$\Phi_\xi(EH(M',\Gamma',\xi'))=EH(M,\Gamma,\xi'\cup\xi)$.

\s\n {\bf Remark.} Observe that the set $\{a''_1,\dots,a''_m\}$
contains arcs of $R_+(\Gamma')\subset S'$.  This is one of the
reasons $\Sigma'$ must be contact-compatible near $\bdry M'$.

\section{Naturality of $\Phi_\xi$} \label{section: well-definition}

In this section we prove that $\Phi_\xi$ does not depend on the
choices made in Section~\ref{section: definition}.  The proofs are
similar to the proofs of well-definition of the $\EH$ class in
\cite{HKM2,HKM3}, and we will only highlight the differences. The
proof of naturality under isotopy is identical to the proof of
\cite[Lemma 3.3]{HKM2}, and will be omitted.

\subsection{Handlesliding}

Consider the Heegaard surface $\Sigma$ and two sets of compressing
disks $(\beta,\alpha)$, $(\overline\beta,\overline\alpha)$ which are
contact-compatible with respect to $\xi\cup\zeta$.  In particular,
$(\beta'',\alpha'')$ and $(\overline\beta'',\overline\alpha'')$
correspond to bases $\{a''_1,\dots,a''_m\}$ and
$\{\overline{a}''_1,\dots, \overline{a}''_m\}$ for
$(S-P',R_+(\Gamma))$.

There are two types of operations to consider:

\be
\item[(A)] Arc slides in the $S-P'$ region, while fixing $\alpha'$
and $\beta'$.
\item[(B)] Handleslides within $\Sigma'$, while preserving the
contact-compatible $\alpha''$ and $\beta''$.
\ee

\begin{lemma}\label{lemma:typeAorB}
Suppose the closure of each component of
$S-P'-\overline{R_+(\Gamma)}$ intersects $\Gamma$ along at least two
arcs. Then one can take $(\beta,\alpha)$ to $(\overline\beta,
\overline\alpha)$ through a sequence of moves of type {\em (A)} or
{\em (B)}.
\end{lemma}

The required connectivity of $S-P'-\overline{R_+(\Gamma)}$ was
already incorporated in the definition in Section~\ref{section:
definition}.

\begin{figure}[ht]
\begin{overpic}[width=5cm]{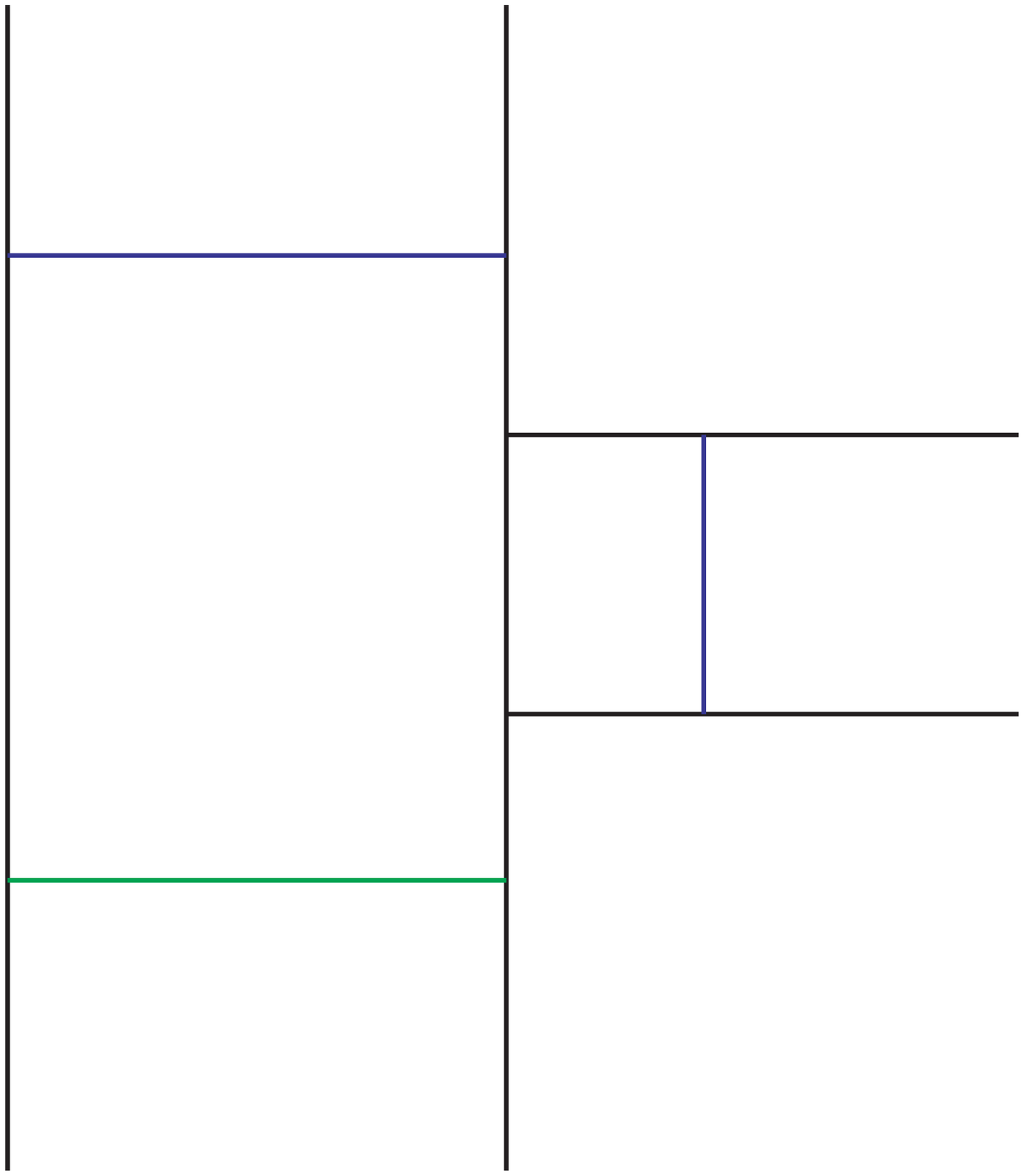}
\put(12,50){\tiny $S-P'$} \put(70,50){\tiny $P'$} \put(54,50){\tiny
$c'$} \put(17,80.5){\tiny $a_i''$} \put(17,19.2){\tiny
$\overline{a}_i''$}
\end{overpic}
\caption{} \label{connection}
\end{figure}

\begin{proof}
According to \cite[Lemma 3.3]{HKM2}, any basis $\{a_i''\}_{i=1}^m$
for $(S-P', R_+(\Gamma))$ can be taken to any other basis
$\{\overline{a}_i''\}_{i=1}^m$ for $(S-P',R_+(\Gamma))$ through a
sequence of arc slides within $S-P'-\overline{R_+(\Gamma)}$,
assuming sufficient connectivity of $S-P'-\overline{R_+(\Gamma)}$.
We must, however, not forget $\Sigma'$.  If $\Sigma'$ is taken into
account, the situation given in Figure~\ref{connection} must be
dealt with: Locally $P'$ is attached to $S-P'$ along an arc $c'$ (in
the diagram, we have pushed $c'$ into $P'$), and we would like to
arc slide $a_i''$ over $c'$ to obtain $\overline{a}_i''$. However,
this $c'$ may not be an $a_i'$. If this is the case, we must perform
a sequence of handleslides on $\alpha'$ and $\beta'$ first (while
fixing $\alpha''$ and $\beta''$), so that $\alpha_1'=\bdry
(c'\times[0,1])$ and $\beta_1'= (c'\times\{1\})\cup
(h(c')\times\{0\})$.  This is possible since we required at least
two arcs $\{p\}\times[-\varepsilon',0]$ in the definition of $K'_1$.
Then we may arc slide $a_i''$ over $c'$.
\end{proof}

We now discuss naturality under the moves (A) and (B).

(A).  Recall that an arc slide corresponds to a sequence of two
handleslides by \cite{HKM2}. For each handleslide of an arc slide in
the $S-P'$ region, the ``tensoring with $\Theta$'' map $\Psi$ sends
the $\EH$ class $\mathbf{x}''$ on $S-P'-\overline{R_+(\Gamma)}$ to
the $\EH$ class on $S-P'-\overline{R_+(\Gamma)}$, also called
$\mathbf{x}''$ by abuse of notation. Since the $\alpha''$ and
$\beta''$ are used up, the restriction of $\Psi$ to the remaining
$r$-tuple $\mathbf{y}\in \CF(\Sigma',\beta',\alpha')$ is the natural
``tensoring with $\Theta'$'' map $\Psi'$ from
$\CF(\Sigma',\beta',\alpha')$ to
$\CF(\overline{\Sigma}',\overline{\beta'},\overline{\alpha'})$.
Therefore,
$$\Psi(\mathbf{y}, \mathbf{x}'')=(\Psi'(\mathbf{y}),\mathbf{x}'').$$
The proof is identical to the proof of \cite[Lemma~5.2]{HKM2}.

(B). The ``tensoring with $\Theta$'' operation for a handleslide in
the $\Sigma'$ region clearly sends $\mathbf{x}''$ to $\mathbf{x}''$
as well. Therefore we have:
$$\Psi(\mathbf{y}, \mathbf{x}'')=(\Psi'(\mathbf{y}),\mathbf{x}'').$$

\subsection{Stabilization}

In this subsection we prove naturality under stabilization. For
this, we need to prove two things: (A) naturality under
stabilizations (contact or otherwise) inside $M'$, and (B)
naturality under positive (contact) stabilizations inside $M-M'$.

Let $A$ be a surface with nonempty boundary and $B\subset A$ be a
subsurface.  Let $c$ be a properly embedded arc in $A$; after
isotopy rel boundary, we assume $c$ intersects $\bdry B$
transversely and efficiently. Then we define the {\em complexity of
$c$ with respect to $(A,B)$} as the number of subarcs of $c$ which
are contained in $B$ and have both endpoints on the common boundary
of $A-B$ and $B$.

Given two Heegaard splittings $(\Sigma',\beta',\alpha')$ and
$(\overline\Sigma',\overline\beta',\overline\alpha')$ for
$(-M',-\Gamma')$ which are contact compatible with respect to
$\zeta$ near $\bdry M'$ (i.e., of the type constructed in Step 1 of
Section~\ref{section: definition}) and their extensions
$(\Sigma,\beta,\alpha)$ and
$(\overline\Sigma,\overline\beta,\overline\alpha)$ to $(-M,-\Gamma)$
which are contact compatible with respect to $\xi\cup\zeta$ (i.e.,
of the type constructed in Step 2 of Section~\ref{section:
definition}), we first find a common stabilization
$(\widetilde\Sigma,\widetilde\beta,\widetilde\alpha)$, which is also
contact compatible with respect to $\xi\cup\zeta$. If we place a
line (resp.\ tilde) over a symbol, then it stands for the
corresponding object for
$(\overline\Sigma,\overline\beta,\overline\alpha)$ (resp.\
$(\widetilde\Sigma,\widetilde\beta,\widetilde\alpha)$), e.g.,
$\overline{K}'$ is $K'$ for
$(\overline\Sigma,\overline\beta,\overline\alpha)$. (The exception
is $\overline{R_+(\Gamma')}$, which refers to the closure of
$R_+(\Gamma')$.)

\s\n (A) We will first discuss the subdivision on $(-M',-\Gamma')$.
Take $\widetilde\varepsilon'>0$ so that $\widetilde\varepsilon'\ll
\varepsilon',\overline\varepsilon'$. Given the Legendrian portion
$L'_0=K'_0\cup K'_1$ of the $1$-skeleton $K'$, we successively
attach Legendrian arcs $c'_i$ to $L'_i$ to obtain $L'_{i+1}$ in the
following order:
\begin{enumerate}
\item[($\alpha$)] First attach arcs to construct the Legendrian $1$-skeleton of a sufficiently
fine Legendrian cell decomposition of $T_{-\widetilde\varepsilon'}$,
after possibly applying Legendrian realization.
\item[($\beta$)] Then attach Legendrian arcs of the type $\{p\}\times [-\widetilde\varepsilon',0]$
with $p\in \Gamma'$.
\end{enumerate}
The arcs are attached so that in the end we obtain a Legendrian
graph containing
$\widetilde{L}'_0=\widetilde{K}'_0\cup\widetilde{K}'_1$, where
$\widetilde{K}'_0$ is a Legendrian skeleton of
$T_{-\widetilde\varepsilon'}$ and $\widetilde{K}'_1$ is the union of
arcs of type $\{p\}\times [-\widetilde\varepsilon',0]$, and so that
the restrictions of $L_0'$ and $\overline{L}'_0$ to
$T\times[-\widetilde\varepsilon',0]$ are subsets of
$\widetilde{L}_0'$. If we start with
$\overline{L}'_0=\overline{K}'_0\cup \overline{K}'_1$ instead, then
there is a sequence $\overline{L}'_i$ which eventually yields a
Legendrian graph containing $\widetilde{L}'_0$. The stabilizations
of contact type will be treated in (A$_1$). Next extend
$\widetilde{L}'_0$ to a common refinement $\widetilde{K}'$ of $K'$
and $\overline K'$ by subdividing on
$M'-(T\times[-\widetilde\varepsilon',0])$. These stabilizations will
be treated in (A$_2$).

\s\n (A$_1$) The above attachments of Legendrian arcs are done in
the same way as in \cite[Theorem~1.2]{HKM3} and in particular
\cite[Figures~1, 3, and 4]{HKM3}.

\begin{figure}[ht]
\begin{overpic}[width=8.4cm]{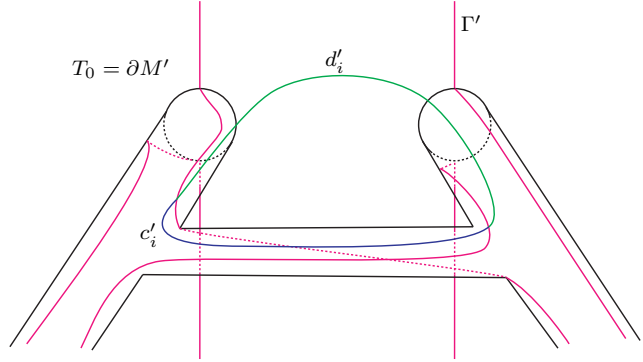}
\put(10,45){\tiny $T_0=\bdry M'$} \put(71.5,52){\tiny $\Gamma'$}
\put(50,46.5) {\tiny $d'_i$} \put(20.7,19){\tiny $c'_i$}
\end{overpic}
\caption{Arc of type ($\alpha_1$). The surface in the back is
$T_0=\bdry M'$, whose orientation as the boundary of $M'$ points
into the page. The cylinders on the left and right are thickenings
of arcs $\{p\}\times[-\varepsilon',0]$ and
$\{q\}\times[-\varepsilon',0]$ of $K'_1$, and the horizontal
cylinder is a thickening of $c'_i$. The blue arc is $c'_i$ and the
green arc is its isotopic copy $d'_i$. } \label{attacharc1}
\end{figure}

\begin{figure}[ht]
\begin{overpic}[width=8.8cm]{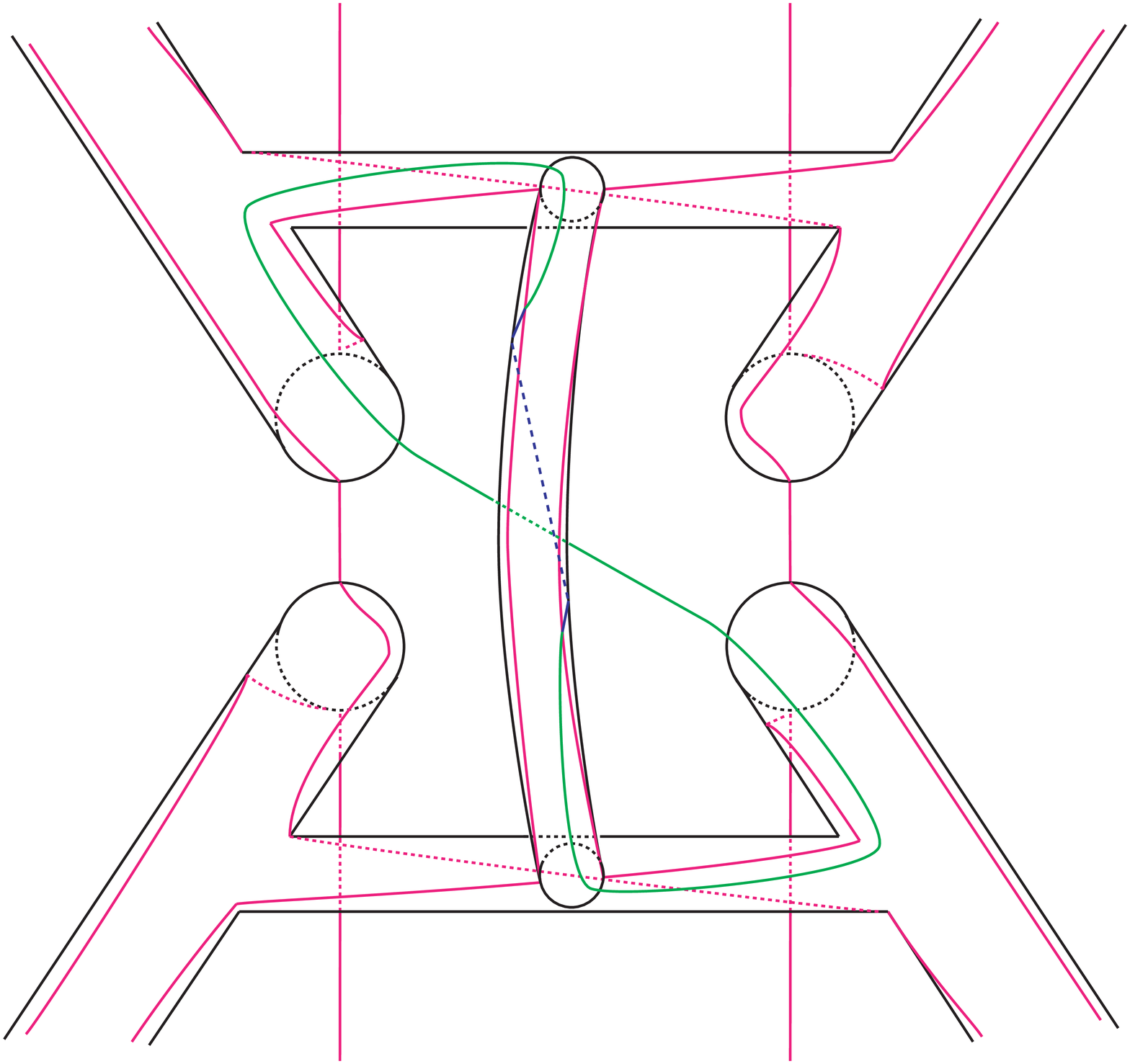}
\put(47,57){\tiny $c_i'$} \put(58,42){\tiny $d_i'$}
\end{overpic}
\caption{Arc of type ($\alpha_2$).} \label{attacharc3}
\end{figure}

\begin{figure}[ht]
\begin{overpic}[width=8.5cm]{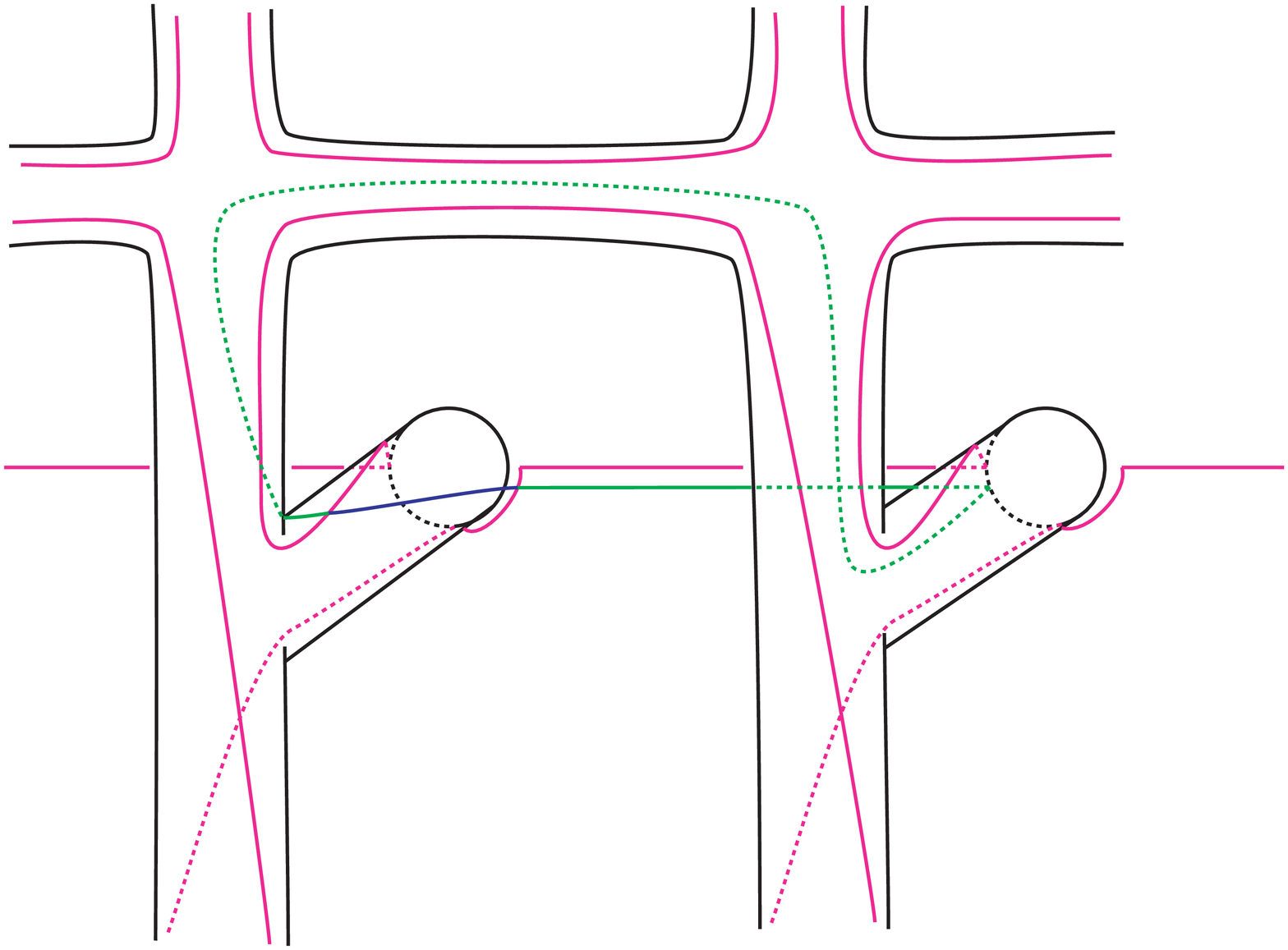}
\put(33,36){\tiny $c_i'$} \put(47,32.3){\tiny $d_i'$}
\end{overpic}
\caption{Arc of type ($\beta$).} \label{attacharc2}
\end{figure}

The attachment of arcs $c'_i\subset T_{-\widetilde\varepsilon'}$ of
type ($\alpha$) can be decomposed into three stages ($\alpha_1$),
($\alpha_2$) and ($\alpha_3$). Figure~\ref{attacharc1} depicts an
arc of type ($\alpha_1$). An arc $c'_i$ of type ($\alpha_1$)
connects between $\{p\}\times[-\varepsilon',0]$ and
$\{q\}\times[-\varepsilon',0]$, where $p,q\in \Gamma'$. After
attaching all the arcs of type ($\alpha_1$), we attach the arcs of
type ($\alpha_2$), depicted in Figure~\ref{attacharc3}. Here, the
arc $c'_i$ connects between two arcs of type ($\alpha_1$) and does
not cross the dividing set of $T_{-\widetilde\varepsilon'}$.
Finally, an arc of type ($\alpha_3$) is an arc that intersects the
dividing set of $T_{-\widetilde\varepsilon'}$ exactly once, and in
its interior. Arcs of type ($\alpha_1$), ($\alpha_2$), and
($\alpha_3$) are sufficient to construct the Legendrian skeleton of
$T_{-\widetilde\varepsilon'}$. Figure~\ref{attacharc2} depicts an
arc attachment of type ($\beta$).

In particular, we observe that the following holds:
\begin{itemize}
\item Each endpoint of $c'_i$ lies on $\Gamma_{\bdry (M'-N(L'_i))}$,
and $int(c'_i)\subset int(M'-N(L'_i))$.
\item $N(L'_{i+1})=N(c'_i)\cup N(L'_i)$, and $L'_{i+1}$ is a Legendrian
graph so that $N(L'_{i+1})$ is its standard neighborhood.
\item There is a Legendrian arc $d'_i$ on $\bdry (M'-N(L'_i))$ with the same
endpoints as $c'_i$, after possible application of the Legendrian
realization principle. The arc $d'_i$ intersects $\Gamma_{\bdry
(M'-N(L'_i))}$ only at its endpoints.
\item \label{unknot}
The Legendrian knot $\gamma'_i=c'_i\cup d'_i$ bounds a disk in
$M'-N(L'_i)$ and has $tb(\gamma'_i)=-1$ with respect to this disk.
This implies that $c'_i$ and $d'_i$ are isotopic relative to their
endpoints inside the closure of $M'-N(L'_i)$.
\end{itemize}

For simplicity, consider the situation of attaching a single arc
$c'_0$ to $L_0'$ to obtain $\widetilde{L}'_0$. Consider the (very)
partial open book decomposition on $T\times[-{3\varepsilon'\over
2},0]$, corresponding to the decomposition into $H_2'=N(K'_0\cup
K'_1)$ and $H_1'=(T\times[-{3\varepsilon'\over 2},0])-N(K'_0\cup
K'_1)= S'\times[0,1]/\sim$. The monodromy map is $h': Q'\rightarrow
S'$ as before. The arc $c'_0$ can be viewed as a Legendrian arc on
$S'\times\{{1\over 2}\}$ with endpoints on $\bdry Q'\times \{{1\over
2}\}$. Hence, removing a neighborhood of $c'_0$ from $H_1'$ and
adding it to $H_2'$ is equivalent to the following positive
(contact) stabilization: Let $e_0'$ be the Legendrian arc on
$S'=S'\times\{1\}$ which is Legendrian isotopic to $c_0'$ rel
endpoints, via an isotopy inside $H_1'$. Add a $1$-handle to $S'$
along the endpoints of $e'_0$ to obtain $\widetilde{S}'$, and
complete $e'_0$ to a closed curve $\gamma_0$ on $\widetilde{S}'$ by
attaching the core of the $1$-handle. Then the stabilization is the
data $(\widetilde{S}', R_+(\Gamma'), \widetilde{Q}',
\widetilde{h}'=R_{\gamma_0}\circ h')$, where $R_{\gamma_0}$ is a
positive Dehn twist about $\gamma_0$ and $\widetilde{Q}'$ is the
domain of $\widetilde{h}'$. Let $(S,h)$ (resp.\
$(\widetilde{S},\widetilde{h})$) be the partial open book which
extends $(S',h')$ (resp.\ $(\widetilde{S}',\widetilde{h}')$).

\begin{lemma} \label{lemma: choice of c0}
The arc $c'_0$ can be chosen so that the corresponding $e'_0\subset
S'$ has complexity $0$ with respect to $(S',P')$ and complexity at
most $1$ with respect to $(S',R_+(\Gamma'))$.
\end{lemma}

\begin{proof}
We treat the ($\alpha_1$) case, and leave the other cases to the
reader.  Refer to Figure~\ref{attacharc1}; in the figure replace
$c_i',d_i'$ by $c_0',d_0'$. If $d'_0$ intersects $R_+(\Gamma')$,
then $d'_0$, viewed on $S'\times\{1\}$, is the desired isotopic copy
$e'_0$ of $c'_0$. It is clear that $e'_0$ has complexity $0$ with
respect to $(S',P')$ and complexity $1$ with respect to
$(S',R_+(\Gamma'))$. On the other hand, if $d'_0$ intersects
$R_-(\Gamma')$, then we need to isotop $c'_0$ towards
$T_{-\varepsilon'/2}$ instead, in order to obtain $e'_0$. The
procedure is still the same --- in Figure~\ref{attacharc1} assume
that the surface in the back is $T_{-\varepsilon'/2}$ (instead of
$T_0$) so that $\bdry (M'-(T\times (-{\varepsilon'\over 2},0]))$
points out of the page. The resulting $e'_0$ has complexity $0$ with
respect to both $(S',P')$ and $(S',R_+(\Gamma'))$.
\end{proof}

In view of Lemma~\ref{lemma: choice of c0}, there exists a basis
$\{a_1',\dots,a_k'\}$ of $(S',R_+(\Gamma'))$ and an extension to a
basis $\{a_1',\dots,a_k',a_1'',\dots,a_m''\}$ of $(S,R_+(\Gamma))$,
so that $e'_0$ does not intersect any basis element. Let $a_0'$ be
the cocore of the $1$-handle of the stabilization along $e'_0$, and
let $b'_0$ be the pushoff of $a_0'$. Then let $\alpha_0'=\bdry
(a_0'\times[0,1])$ and
$\beta_0'=(b_0'\times\{1\})\cup(R_{\gamma_0}(b_0')\times\{0\})$,
where both are viewed on $\bdry H_1'$. Observe that $\beta_0'$ does
not intersect any of
$\{\alpha_1',\dots,\alpha_k',\alpha_1'',\dots,\alpha_m''\}$, where
$\alpha_i'=\bdry(a_i'\times[0,1])$ and
$\alpha_i''=\bdry(a_i''\times[0,1])$. Since $\beta_0'$ also does not
intersect the remaining $\alpha'$-curves
$\alpha_{k+1}',\dots,\alpha_{k+l}'$, the only intersection between
$\beta_0'$ and some $\alpha$-curve is the sole intersection with
$\alpha_0'$.

Let $\Psi: \CF(\Sigma,\beta,\alpha)\rightarrow \CF
(\widetilde\Sigma,\widetilde\beta,\widetilde\alpha)$ be the
composition of (Heegaard decomposition) stabilization and
handleslide maps corresponding to the stabilization along $e'_0$. We
have the following:

\begin{lemma}
The map $$\Psi: \CF(\Sigma,\beta,\alpha)\rightarrow \CF
(\widetilde\Sigma,\widetilde\beta,\widetilde\alpha)$$ is given by:
$$(\mathbf{y},\mathbf{x}'')\mapsto (\Psi'(\mathbf{y}),\widetilde{\mathbf{x}}''),$$
where $\mathbf{y}\in CF(\Sigma',\beta',\alpha')$,  $\mathbf{x}''$
$($resp.\ $\widetilde{\mathbf{x}}'')$ is the $\EH$ class in the
$S-P'$ region $($resp.\ $\widetilde{S}-\widetilde{P}'$ region$)$ for
$(\Sigma,\beta,\alpha)$ $($resp.\
$(\widetilde\Sigma,\widetilde\beta,\widetilde\alpha))$, and $\Psi'$
is the natural map from $(\Sigma',\beta',\alpha')$ to
$(\widetilde\Sigma',\widetilde\beta',\widetilde\alpha')$.
\end{lemma}

\begin{proof}
This follows from the technique in \cite[Lemma~3.5]{HKM3}. We use
the fact that the only intersection between $\beta_0'$ and an
$\alpha$-curve is the unique intersection with $\alpha_0'$.  We
decompose the positive stabilization along $e'_0$ into a trivial
stabilization, followed by a sequence of handleslides. (By a {\em
trivial stabilization} we mean the addition of a $1$-handle to
$\Sigma'$, together with curves $\alpha_0'$ and $\beta_0'$ that
intersect each other once, say at $x_0'$, and no other $\alpha_i'$,
$\beta_i'$, $i=1,\dots, s$, and where the regions of
$\Sigma'-\cup_{i=0}^s \alpha_i'- \cup_{i=0}^s\beta_i'$ adjacent to
$x_0'$ are path-connected to $R_+(\Gamma')$.) This is done exactly
as described in \cite[Lemma~3.5]{HKM3}: whenever $\beta_i'$ (could
be $\beta_i''$) intersects $e'_0\times\{0\}$, and
$\overline{\beta}_i'$ is the result of applying a positive Dehn
twist about $\gamma_0\times\{0\}$, then $\overline{\beta}_i'$ can be
obtained from $\beta_i'$ by applying a trivial stabilization,
followed by handleslides over $\beta_0'$ as in
\cite[Figure~7]{HKM3}. Here, the triple diagrams are weakly
admissible for the same reasons as \cite[Lemma~3.5]{HKM3}.

The slight complication that we need to keep in mind is that the
arcs $h(b_i'')$, where $b''_i$ is the usual pushoff of $a_i''$ in
$S-P'$, may enter the region $S'$ and intersect $e'_0$. If
$h(b_i'')$ intersects $e'_0$, then the ``tensoring with $\Theta$''
map corresponding to handlesliding $\beta_i''$ over $\beta_0'$ sends
$\EH$ to $\EH$ in the $S-P'$ region and restricts to the natural
``tensoring with $\Theta$'' map in the $\Sigma'$ region. (The proof
is the same as that of \cite[Lemma~3.5]{HKM3}. Also refer to
\cite[Figure~8]{HKM3}.) On the other hand, if $\beta_i'$ intersects
$e'_0\times\{0\}$, the $S-P'$ region is unaffected (hence $\EH$ is
mapped to $\EH$ in the $S-P'$ region), and we are doing a standard
handleslide map in the $\Sigma'$ region.
\end{proof}

\s\n (A$_2$) Next we discuss the effect of a stabilization, in the
handlebody sense, in the portion of $M'$ which is not
contact-compatible, i.e., away from
$T\times[-\widetilde\varepsilon',0]$.  Assume all the contact
stabilizations have already taken place on
$T\times[-\widetilde\varepsilon',0]$.  By abuse of notation, we
reset $\varepsilon'=\widetilde\varepsilon'$ and use the same
notation $S'$, $P'$, $Q'$, $h'$, $\Sigma'$, $K'$, $K'_0$, $K'_1$,
$K'_2$, $K'_3$, $H_1'$, $H_2'$, used in Step 1 of
Section~\ref{section: definition}, for the new (finer) Heegaard
decomposition which is contact-compatible on
$T\times[-\varepsilon',0]=T\times[-\widetilde\varepsilon',0]$.

\begin{claim}
The stabilization can be decomposed into a trivial stabilization,
followed by a sequence of handleslides which avoids $R_-(\Gamma')$.
\end{claim}

\begin{proof}
Observe that the arc of stabilization $c'_0$ is contained in
$N(K'_3)$. The meridian of the tubular neighborhood of $c'_0$ will
be called $\beta_0'$ and it is not difficult to see that there
exists a curve $\alpha_0'$ which intersects $\beta_0'$ once and lies
on $\bdry N(K'_3)-\bdry H_1'$. (Note that $\beta_0'$ only intersects
$\alpha_0'$.) After a sequence of handleslides that takes place away
from $\bdry H_1'$ (in fact the change takes place in a neighborhood
of $\alpha_0'\cup \beta_0'$), we may assume that $\alpha_0'$ and
$\beta_0'$ satisfy the conditions of a trivial stabilization.
\end{proof}

The claim implies that the handleslides and stabilization do not
interact with $h(a_i'')$ (or equivalently with $\beta_i''$). Hence
the $\EH$ class is mapped to the $\EH$ class in the $S-P'$ region,
and we are doing a standard sequence of handleslide maps plus one
stabilization in the $\Sigma'$ region.

\s\n (B)  Next we discuss the subdivision on $M-M'$.  Take
$\widetilde\varepsilon''>0$ so that $\widetilde\varepsilon''\ll
\varepsilon'', \overline\varepsilon''$.  Consider
$N_{\widetilde\varepsilon''}=T\times[0,\widetilde\varepsilon'']$. On
$N_{\widetilde\varepsilon''}$ attach the following arcs in the given
order to $K'''$ to obtain a common refinement of the restriction of
$K'''$ and $\overline{K}'''$ to $N_{\widetilde\varepsilon''}$:
\begin{enumerate}
\item First attach the Legendrian skeleton of a sufficiently fine Legendrian cell
decomposition of $T_{\widetilde\varepsilon''}$, after possibly
applying Legendrian realization.
\item Then attach Legendrian arcs of type
$\{q\}\times[0,\widetilde\varepsilon'']$ with $q\in \Gamma'$.
\end{enumerate}
Each of the above Legendrian arc attachments leads to a
stabilization --- however, since the arcs are contained in the
complement of $S\times[0,1]/\sim$, the stabilization is a {\em
precomposition} $h\mapsto h\circ R_{\gamma}$. More precisely, let
$c$ be the Legendrian attaching arc.  There is an isotopy of $c$ rel
endpoints, inside the complement of $S\times[0,1]/\sim$, to an arc
$e\subset S-P'$, viewed on $S\times\{1\}$, and also to $h(e)$,
viewed on $S\times\{0\}$. Observe that $h(e)$ may enter the
$R_-(\Gamma')$ region. Add a $1$-handle to $S$ along the endpoints
of $e$ to obtain $\widetilde{S}$, and complete $e$ to a closed curve
$\gamma$ by attaching the core of the handle.

In the following lemma, we identify $S=S\times\{0\}$ and determine
the complexity of the restriction of $h(e)$ to $S'$ and to $P'$.
Observe that $h(e)\cap S'=h(e)\cap \overline {R_-(\Gamma')}$.

\begin{lemma} \label{lemma:h(d)} $\mbox{}$
\begin{enumerate}
\item $h(e)$ has complexity at most one with respect to $(S,S')$.
\item $h(e)$ has complexity at most one with respect to $(S,P')$.
\end{enumerate}
\end{lemma}

\begin{proof}
We isotop $c$ rel endpoints in two stages: first through the product
structure given by the complement of $S\times[0,1]/\sim$, and then
through the product structure given by $S\times[0,1]/\sim$.

(1) follows from examining the proof of \cite[Theorem~1.2]{HKM3} as
in Lemma~\ref{lemma: choice of c0}. The three types of arc
attachments are ($\alpha_1$), ($\alpha_2$), and ($\alpha_3$).
Consider an arc of type ($\alpha_1$), given in
Figure~\ref{attacharc1}. In the current case, the surface in the
back is still $T_0=\bdry M'$, but the orientation is pointing out of
the page; also $c_i'$ and $d_i'$ should be changed to $c$ and $d$.
If the arc $d$ intersects $R_-(\Gamma')$ (where the orientation on
$\bdry M'$ is the orientation induced from $M'$), then $d$ is an arc
on $S\times\{0\}$, which means that $d=h(e)$. Hence $h(e)$ has
complexity $1$ with respect to $(S,R_-(\Gamma'))$, and also
complexity at most $1$ with respect to $(S,S')$. On the other hand,
if $d$ intersects $R_+(\Gamma')$, then $d=e$. Hence $h(e)$ is
contained in $(S-S')\times\{0\}$ and has complexity zero with
respect to $(S,S')$. It follows that $h(e)$ also has complexity zero
with respect to $(S,P')$.  The arcs of type ($\alpha_2$) and
($\alpha_3$) are treated similarly.

(2) follows from considerations similar to \cite[Section 5, Example
2]{HKM3}. Suppose $d=h(e)$, i.e., $d$ intersects $R_-(\Gamma')$.
(The situation of $d=e$ is easier, and is left to the reader.) Then
Figure~\ref{pusharc} depicts what happens when we push $h(e)$,
viewed as an arc on $S\times\{0\}$, to $S\times\{1\}$. The surface
in the front is $T_0$ and the surface in the back is
$T_{-\varepsilon'/2}$. The blue arc $(h(e)|_{S'})_0$ is the isotopic
copy of $h(e)|_{S'}$ on $T_0$ or $S'\times\{0\}$, and the green arc
$(h(e)|_{S'})_{-\varepsilon'/2}$ is the copy on $S'\times\{1\}$
which intersects $T_{-\varepsilon'/2}$.
\begin{figure}[ht]
\s
\begin{overpic}[width=8.5cm]{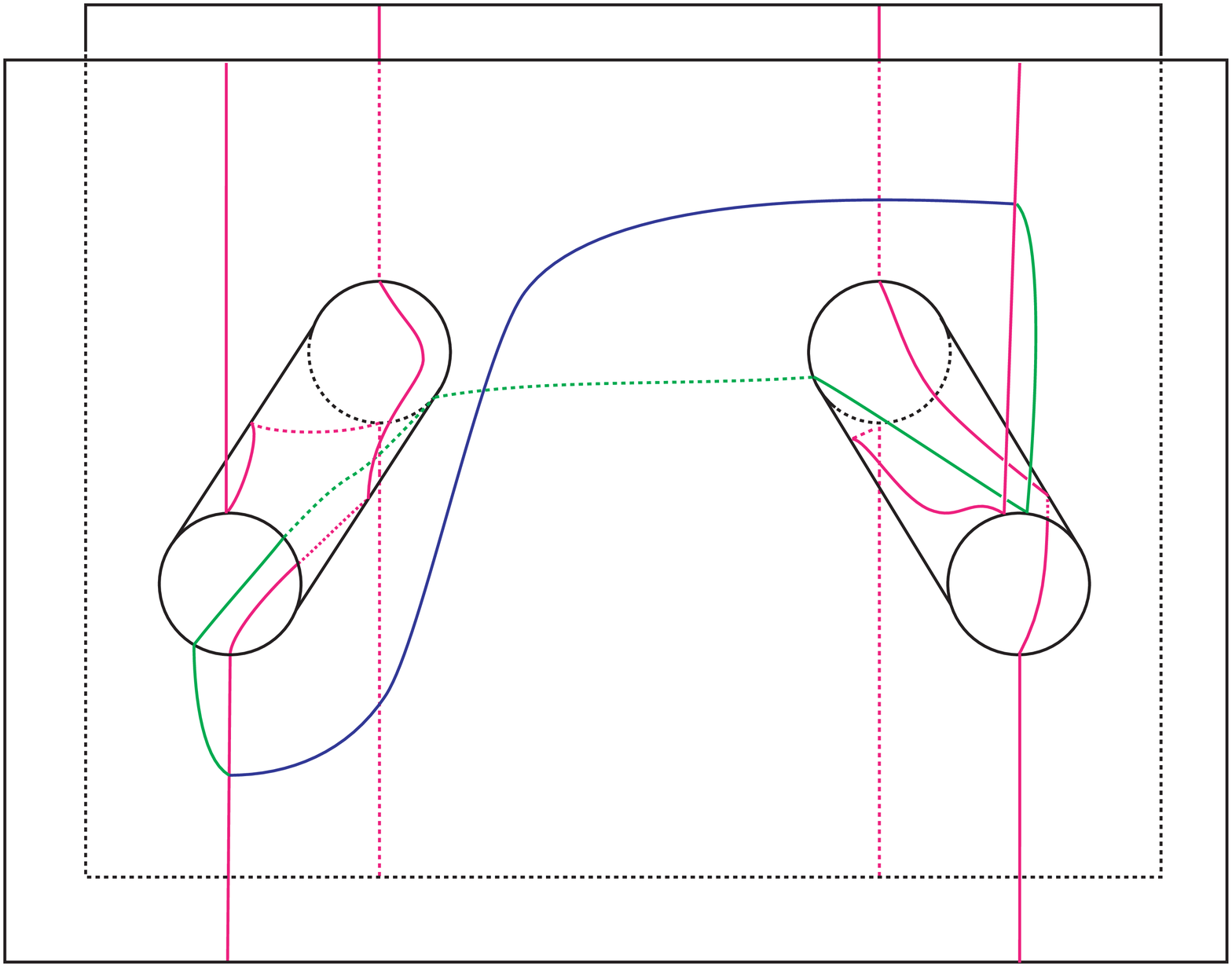}
\put(50,65){\tiny $(h(e)|_{S'})_0$} \put(50,3){\tiny $-$}
\put(10,3){\tiny $+$} \put(42,50){\tiny
$(h(e)|_{S'})_{-\varepsilon'/2}$} \put(-6.5,55){\tiny $T_0$}
\put(50,82){\tiny $T_{-\varepsilon'/2}$}
\end{overpic}
\caption{} \label{pusharc}
\end{figure}
We easily see that $h(e)$ has complexity $1$ with respect to
$(S,P')$. The arc corresponding to \cite[Figure~4]{HKM3} is simpler,
and does not enter $P'$.
\end{proof}

The following lemma follows from Lemma~\ref{lemma:h(d)}.

\begin{lemma}\label{lemma: basis}
There exists a basis $\{a_1',\dots,a_k'\}$ for $(S',R_+(\Gamma'))$
and basis $\{a_1'',\dots,a_m''\}$ for $(S-P',R_+(\Gamma))$ such that
the following hold:
\begin{enumerate}
\item $a_i''$ does not intersect $e$ for all $i$;
\item all but one of the $a_i'$ or $a_i''$ are disjoint from $h(e)$;
\item one of $a_1'$ or $a_1''$ intersects $h(e)$.
\end{enumerate}
\end{lemma}

The basis $\{a_1',\dots,a_k'\}$ can be used to construct $\alpha'$,
$\beta'$ for $\Sigma'$, and the basis $\{a_1'',\dots,a_m''\}$ gives
an extension to $\alpha$, $\beta$ on $\Sigma$.

\begin{proof}
Consider the ($\alpha_1$) case where $d$ intersects $R_-(\Gamma')$.
By (2) of Lemma~\ref{lemma:h(d)}, there exists a basis
$\{a_1',\dots,a_k'\}$ for $(S',R_+(\Gamma'))$ so that $a_1'$
intersects $h(e)$ once, and the remaining $a_i'$, $i=2,\dots, k$, do
not intersect $h(e)$.  Next observe that $e\subset S-P'$ and does
not intersect the $R_+(\Gamma)$ region. It is possible to choose a
basis $\{a_1'',\dots,a_m''\}$ for $(S-P',R_+(\Gamma))$ which does
not intersect $e$, as well as $h(e)$.  The other cases are similar.
\end{proof}

\begin{figure}[ht]
\s
\begin{overpic}[width=17cm]{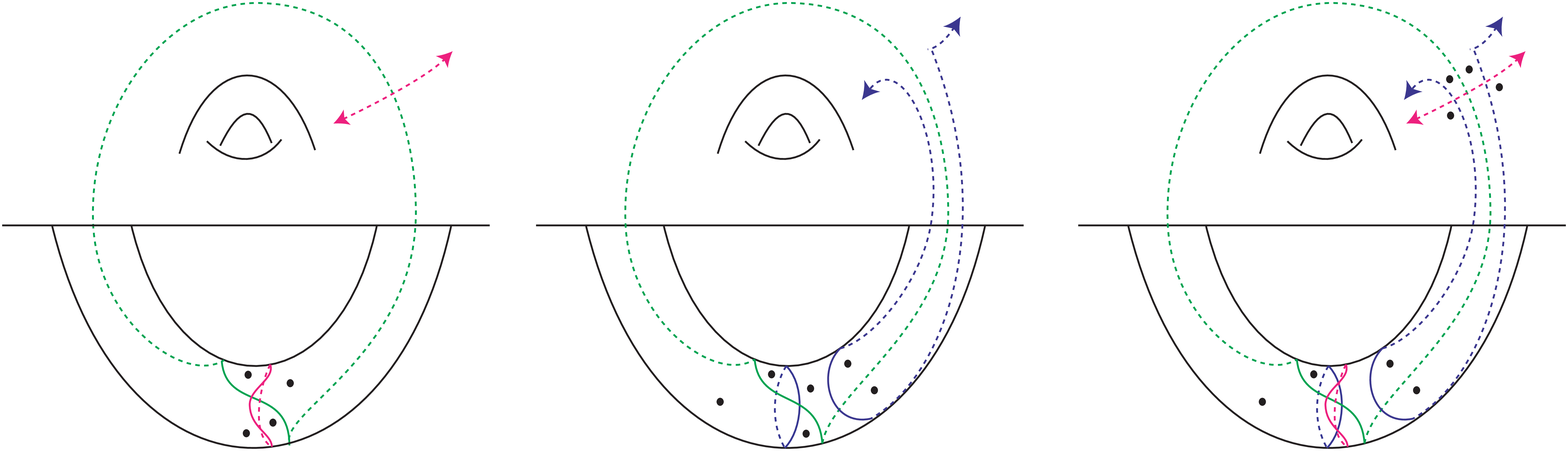}
\put(15,-1.7){\tiny $\alpha_0''$} \put(12,3.1){\tiny $\beta_0''$}
\put(27.5,26){\tiny $\alpha_1'$} \put(49,-1.7){\tiny $\gamma_0''$}
\put(61.5,26){\tiny $\gamma_1'$}
\end{overpic}
\s \caption{} \label{handleslide2}
\end{figure}

Let us consider the case where $a_1'$ intersects $h(e)$.  (The other
case is similar.)  When we stabilize $S$ along $e$, we add the
cocore $a_0''$ of the $1$-handle and obtain the corresponding
$\alpha_0''$ and $\beta_0''$.  The only intersection point of
$\alpha_0''$ with any $\beta$ arc is with $\beta_0''$, which we call
$x_0''$.  Hence we expect the following diagram to commute:
$$\begin{diagram}
\SFH(\beta',\alpha') & \rTo^{\hspace{.3in}\Phi_\xi}  & \SFH(\beta,\alpha)  \\
\dTo^{\Psi} & &  \dTo^{\Psi} \\
\SFH(\beta',\alpha') & \rTo^{\Phi_\xi} &
\SFH(\beta\cup\{\beta_0''\},\alpha\cup\{\alpha_0''\})  \\
\end{diagram}   $$
However, our stabilization is not a trivial stabilization, as
$\alpha_1'$ intersects $\beta_0''$ in one point.  Therefore we need
to decompose the stabilization into a trivial stabilization,
followed by a handleslide.  This will be done in a manner similar to
\cite[Lemma~3.5]{HKM3}.  Let $\gamma_i''$ be pushoffs of
$\alpha_i''$ for all $i$, $\gamma_j'$ be pushoffs of $\alpha_j'$ for
all $j\not=1$, and $\gamma_1'$ be obtained by pushing $\alpha_1'$
over $\alpha_0''$, as depicted in Figure~\ref{handleslide2}. In
Figure~\ref{handleslide2}, we place black dots in regions that are
path-connected to $\Gamma$; in other words, holomorphic curves are
not allowed to enter such regions.

Now consider the following diagram:
$$\begin{diagram}
\SFH(\beta',\gamma') &\rTo^{\Phi_\xi}  & \SFH(\beta,\gamma)  \\
& &\dTo_{\Psi_2}\\
\dTo^{\Psi_1}& & \SFH(\beta\cup \{\beta_0''\}, \gamma\cup \{\gamma_0''\}) \\
& & \dTo_{\Psi_3}\\
\SFH(\beta',\alpha') & \rTo^{\hspace{.3in}\Phi_\xi\hspace{.3in}} &
\SFH(\beta\cup\{\beta_0''\},\alpha\cup\{\alpha_0''\})  \\
\end{diagram}   $$
\s\n For the term $\SFH(\beta',\gamma')$ in the upper left-hand
corner, $\gamma'$ is the set consisting of all the $\gamma_i'$; for
$\SFH(\beta,\gamma)$ in the upper right-hand corner, $\gamma$ is the
set consisting of all the $\gamma_i''$ and $\gamma_i'$, with the
exception of $\gamma_0''$. This means that $(\beta,\gamma)$ is
obtained from the middle diagram of Figure~\ref{handleslide2} by a
destabilization; hence $\gamma$ effectively consists of pushoffs of
$\alpha$.  The map $\Psi_2$ is the map which corresponds to the
trivial stabilization, and $\Psi_3$ is the handleslide map which is
the ``tensoring with $\Theta$'' map, where $\Theta$ is the top
generator of $CF(\gamma\cup\{\gamma_0''\},\alpha\cup
\{\alpha_0''\})$.  The slightly tricky feature of this diagram is
that at $\SFH(\beta\cup \{\beta_0''\}, \gamma\cup \{\gamma_0''\})$
we leave the category of diagrams which nicely decompose into the
$M'$ part and the $M-M'$ part. The map $\Psi_1$ is the ``tensoring
with $\Theta'$'' map, where $\Theta'$ is the top generator of
$CF(\gamma',\alpha')$.  The maps $\Phi_\xi$ are the ``tensoring with
the $EH$ class" maps. By the placement of the dots in the right-hand
diagram of Figure~\ref{handleslide2}, it is not difficult to see the
following:

\begin{lemma} \label{lemma:fig3}
The $\EH$ class on $(\beta''\cup \{\beta_0''\}, \gamma''\cup
\{\gamma_0''\})$ is mapped to the $EH$ class on $(\beta''\cup
\{\beta_0''\}, \alpha''\cup\{\alpha_0''\})$ via $\Psi_3$.
\end{lemma}

\begin{proof}
The Heegaard triple diagram is weakly admissible for the same reason
as Lemma~3.5 of \cite{HKM3}, and the details are left to the reader.
In the right-hand diagram of Figure~\ref{handleslide2}, consider the
largest closed connected component $R$ which is bounded by the
$\alpha\cup\{\alpha_0''\},\beta\cup\{\beta_0''\},\gamma\cup
\{\gamma_0''\}$ curves, does not contain a dot (i.e., does not
intersect $\Gamma$), and contains the unique intersection point of
$\beta''_0$ and $\gamma''_0$. The set $R$ is an annulus which is
bounded by $\alpha_0''$ and $\gamma_0''$ on one side, and by
$\alpha_1'$ and $\gamma_1'$ on the other. There are two points of
$\Theta$ in $R$, but only one intersection point of
$\beta\cup\{\beta_0''\}$ and $\gamma\cup\{\gamma''_0\}$. Hence one
of the $\Theta$ points cannot be used towards $R$, namely the
intersection point between $\alpha_1'$ and $\gamma_1'$.  This allows
us to ``erase'' the boundary component of $R$ consisting of
$\alpha_1'$ and $\gamma_1'$, and conclude that $\alpha_0''\cap
\gamma_0''$ is mapped to $\alpha_0''\cap \beta_0''$. The rest of the
tuples of the $EH$ class are straightforward.
\end{proof}

Once the $\EH$ portion is used up by Lemma~\ref{lemma:fig3},
$\Psi_3$ restricts to $\Psi_1$ on the rest of the tuples, i.e.,
those that lie on $\Sigma'$. The commutativity of the above diagram
follows.

Now, inside $M''_{\widetilde\varepsilon''}=M-int(M'\cup
N_{\widetilde\varepsilon''})$, we attach Legendrian arcs to the
Legendrian graph which plays the role of $K''$ so that we have a
common refinement of $K''$ and $\overline{K}''$.  An arc attachment
in this region corresponds to a straightforward stabilization along
$c$ which lies in $S-S'$. The map on Floer homology induced by such
a stabilization clearly sends $\EH$ to $\EH$ and has a natural
restriction to the $\Sigma'$ region.

\section{Properties of the gluing map} \label{section: 4.5}

In this section we collect some standard properties of the gluing
map.

\begin{thm}[Identity] \label{thm: identity}
Let $(M,\Gamma)$ be a sutured manifold and $\xi$ be a
$[0,1]$-invariant contact structure on $\bdry M\times[0,1]$ with
dividing set $\Gamma\times\{t\}$ on $\bdry M\times\{t\}$.  The
gluing map $$\Phi_\xi: \SFH(-M,-\Gamma)\rightarrow
\SFH(-M,-\Gamma),$$ obtained by attaching $(\bdry M\times[0,1],\xi)$
onto $(M,\Gamma)$ along $\bdry M\times \{0\}$, is the identity map
$($up to an overall $\pm$ sign if over $\Z)$.
\end{thm}

The proof of Theorem~\ref{thm: identity} will be given in
Subsection~\ref{proof of id}, after some preliminaries.

\begin{prop}[Composition]
Consider the inclusions $(M_1,\Gamma_1)\subset (M_2,\Gamma_2)\subset
(M_3,\Gamma_3)$ of sutured manifolds, and let $\xi_{12}$ be a
contact structure on $M_2-int(M_1)$ which has convex boundary and
dividing sets $\Gamma_i$ on $\bdry M_i$, $i=1,2$. Similarly define
$\xi_{23}$.  If
$$\Phi_{12}: SFH(-M_1,-\Gamma_1)\rightarrow SFH(-M_2,-\Gamma_2),$$
$$\Phi_{23}: SFH(-M_2,-\Gamma_2)\rightarrow SFH(-M_3,-\Gamma_3),$$
$$\Phi_{13}: SFH(-M_1,-\Gamma_1)\rightarrow SFH(-M_3,-\Gamma_3),$$
are natural maps induced by $\xi_{12}$, $\xi_{23}$, and
$\xi_{12}\cup\xi_{23}$, respectively, then $\Phi_{23}\circ
\Phi_{12}=\Phi_{13}$, $($up to an overall $\pm$ sign if over $\Z)$.
\end{prop}

\begin{proof}
This is immediate, once we unwind the definitions. Let
$(S_1,R_+(\Gamma_1),h_1)$ be a partial open book decomposition for
$(M_1,\Gamma_1,\xi_1)$.  Here $\xi_1$ is arbitrary and may be tight
or overtwisted. Let $(\Sigma_1,\beta_1,\alpha_1)$ be the
corresponding contact-compatible Heegaard splitting. We assume that
the partial open book for $\xi_1$ is sufficiently fine and the
Heegaard splitting is of the type given in Step~1 of
Section~\ref{section: definition}. Extend $(S_1,R_+(\Gamma_1),h_1)$
to $(S_2,R_+(\Gamma_2),h_2)$ via $\xi_{12}$ (of the type given in
Step~2 of Section~\ref{section: definition}), and let
$\mathbf{x}_{12}$ be the $EH$ class for the arcs which complete a
basis for $(S_1,R_+(\Gamma_1),h_1)$ to a basis for
$(S_2,R_+(\Gamma_2),h_2)$. Similarly define $\mathbf{x}_{23}$. Then
the chain map $\Phi_{12}$ maps:
$$\mathbf{y}\mapsto (\mathbf{y},\mathbf{x}_{12}),$$
and $\Phi_{23}$ maps:
$$(\mathbf{y},\mathbf{x}_{12})\mapsto
(\mathbf{y},\mathbf{x}_{12},\mathbf{x}_{23}).$$ This is the same as
$\Phi_{13}(\mathbf{y})$, since $(\mathbf{x}_{12},\mathbf{x}_{23})$
is the $EH$ class for the arcs which complete a basis for
$(S_1,R_+(\Gamma_1),h_1)$ to a basis for $(S_3,R_+(\Gamma_3),h_3)$.
Moreover the extension is of the type given in Step~2 of
Section~\ref{section: definition}.
\end{proof}

\begin{prop}[Associativity]
Let $(M_1,\Gamma_1)$, $(M_2,\Gamma_2)$, and $(M_3,\Gamma_3)$ be
pairwise disjoint sutured submanifolds of $(M,\Gamma)$.  Let $\xi$
be a contact structure defined on $M-int(M_1\cup M_2\cup M_3)$ which
has convex boundary and dividing sets $\Gamma$ on $\bdry M$ and
$\Gamma_i$ on $\bdry M_i$.  Let $(M_{12},\Gamma_{12})$ be a sutured
submanifold of $(M,\Gamma)$ which is disjoint from $M_3$, contains
$M_1$ and $M_2$, and has dividing set $\Gamma_{12}$ on $\bdry
M_{12}$ with respect to $\xi$.  Similarly define
$(M_{23},\Gamma_{23})$.  Then the maps
\begin{equation}
\SFH(-M_1,-\Gamma_1)\otimes\SFH(-M_2,-\Gamma_2)\otimes
\SFH(-M_3,-\Gamma_3) \end{equation} $$\begin{diagram}
& \rTo^{\Phi_{\xi|_{M_{12}-M_1-M_2}}\otimes id} &
\SFH(-M_{12},-\Gamma_{12})\otimes \SFH(-M_3,-\Gamma_3) &
\rTo^{\Phi_{\xi|_{M-M_{12}-M_3}}} & \SFH(-M,-\Gamma)
\end{diagram}$$
and
\begin{equation}
\SFH(-M_1,-\Gamma_1)\otimes\SFH(-M_2,-\Gamma_2)\otimes
\SFH(-M_3,-\Gamma_3) \end{equation}
$$\begin{diagram} & \rTo^{id\otimes
\Phi_{\xi|_{M_{23}-M_2-M_3}}} & \SFH(-M_{1},-\Gamma_{1})\otimes
\SFH(-M_{23},-\Gamma_{23}) & \rTo^{\Phi_{\xi|_{M-M_{1}-M_{23}}}} &
\SFH(-M,-\Gamma)
\end{diagram}$$
are identical $($up to an overall $\pm$ sign if over $\Z)$.
\end{prop}

\begin{proof}
Let $(S_i, R_+(\Gamma_i),h_i)$ be a partial open book decomposition
for $(M_i,\Gamma_i,\xi_i)$, $i=1,2,3$, where $\xi_i$ is arbitrary.
Let $(\Sigma_i,\beta_i,\alpha_i)$ be the corresponding
contact-compatible Heegaard splitting. To define the chain map
$\Phi_{12}=\Phi_{\xi|_{M_{12}-M_1-M_2}}$, we extend
$(S_i,R_+(\Gamma_i),h_i)$ to a partial open book decomposition
$(S_{12},R_+(\Gamma_{12}),h_{12})$ for
$(M_{12},\Gamma_{12},\xi|_{M_{12}-M_1-M_2}\cup\xi_1\cup \xi_2)$.
Then
$$\Phi_{12}:(\mathbf{y}_1,\mathbf{y}_2)\mapsto (\mathbf{y}_1,\mathbf{y}_2,
\mathbf{x}_{12}),$$ where $\mathbf{x}_{12}$ is the $EH$ class for
the arcs which complete a basis for
$$(S_1,R_+(\Gamma_1),h_1)\cup(S_2,R_+(\Gamma_2),h_2)$$ to a basis
for $(S_{12},R_+(\Gamma_{12}),h_{12})$. Next we complete a basis for
$$(S_{12},R_+(\Gamma_{12}),h_{12})\cup (S_3,R_+(\Gamma_3),h_3)$$ to a
basis for the open book $(S_{123},R_+(\Gamma_{123}),h_{123})$
corresponding to $(M,\xi\cup \xi_1\cup\xi_2\cup \xi_3)$, and let
$\mathbf{x}_{(12)3}$ be the $EH$ class for the completing arcs.
Hence, $\Phi_{(12)3}=\Phi_{\xi|_{M-M_{12}-M_3}}$ maps:
$$(\mathbf{y}_1,\mathbf{y}_2, \mathbf{x}_{12})\otimes \mathbf{y}_3 \mapsto
(\mathbf{y}_1,\mathbf{y}_2,\mathbf{y}_3,\mathbf{x}_{12},\mathbf{x}_{(12)3}).$$
Similarly, $\Phi_{1(23)}\circ (id\otimes \Phi_{23})$ sends
$$(\mathbf{y}_1,\mathbf{y}_2,\mathbf{y}_3)\mapsto (\mathbf{y}_1,
\mathbf{y}_2, \mathbf{y}_3,\mathbf{x}_{23},\mathbf{x}_{1(23)}).$$ By
applying a sequence of positive stabilization and basis change moves
in the $M-int(M_1\cup M_2\cup M_3)$ region, as proven in
Section~\ref{section: well-definition}, we see that
$(\mathbf{x}_{12},\mathbf{x}_{(12)3})$ is taken to
$(\mathbf{x}_{23},\mathbf{x}_{1(23)}).$
\end{proof}

\begin{prop} \label{prop: bdry-parallel}
Let $(M',\Gamma')$ be obtained from $(M,\Gamma)$ by decomposing
along a properly embedded surface $T$ with $\bdry$-parallel dividing
set $\Gamma_T$.  The inclusion/direct summand map
$$\SFH(-M',-\Gamma')\rightarrow \SFH(-M,-\Gamma)$$ given in
\cite[Section~6]{HKM3} is the same as the gluing map of
Theorem~\ref{cor: gluing map}.
\end{prop}

Proposition~\ref{prop: bdry-parallel} can be proved using techniques
as that are similar to the proof of Theorem~\ref{thm: identity}
below, and is left to the reader.

\subsection{Proof of Theorem~\ref{thm: identity}} \label{proof of
id} In this subsection we prove Theorem~\ref{thm: identity}.

\subsubsection{Attaching a trivial bypass}  Let
$(S',R_+(\Gamma'),h')$ be a partial open book decomposition for the
triple $(M',\Gamma',\xi')$, where $\xi'$ is any contact structure.
We determine the effect of attaching a trivial bypass on the partial
open book $(S',R_+(\Gamma'),h')$. Let $(M,\Gamma)$ be the result of
attaching a bypass $D$ to $(M',\Gamma')$ along a trivial arc of
attachment $c$, and thickening. (Of course $(M,\Gamma)$ and
$(M',\Gamma')$ are isotopic, but we keep the distinction.) The
boundary $\bdry D$ decomposes into two arcs which intersect at their
common endpoints: the arc of attachment $c\subset \bdry M'$ and the
bypass arc $d$. As described in \cite[Section~5, Example~5]{HKM3},
attaching a neighborhood $N(D)$ of $D$ is equivalent to attaching a
tubular neighborhood of $d$ (a $1$-handle), followed by a
neighborhood $D_0\times[0,1]$ of a disk $D_0\subset D$ which is a
slight retraction of $D$ (a $2$-handle). Now, let $K'$ be the
Legendrian graph in $(M',\Gamma')$ with endpoints on $\Gamma'$,
which gives rise to the partial open book decomposition
$(S',R_+(\Gamma'),h')$. Then the Legendrian graph $K$ for
$(S,R_+(\Gamma),h)$ is obtained from $K'$ by taking the union with a
Legendrian arc $\{pt\}\times [0,1]\subset D_0\times[0,1]$, which is
the cocore of the $2$-handle. The complement of $N(K)$ in $M$ is
product disk decomposable. This decomposition gives rise to an
extension $(S,R_+(\Gamma),h)$ of $(S',R_+(\Gamma'),h')$ to
$(M,\Gamma)$, obtained by attaching a $1$-handle to $S'$. Let $a_0$
be the cocore of the $1$-handle.  The monodromy $h'$ on the
$S'$-portion remains unchanged.

\begin{figure}[ht]
\begin{overpic}[width=12cm]{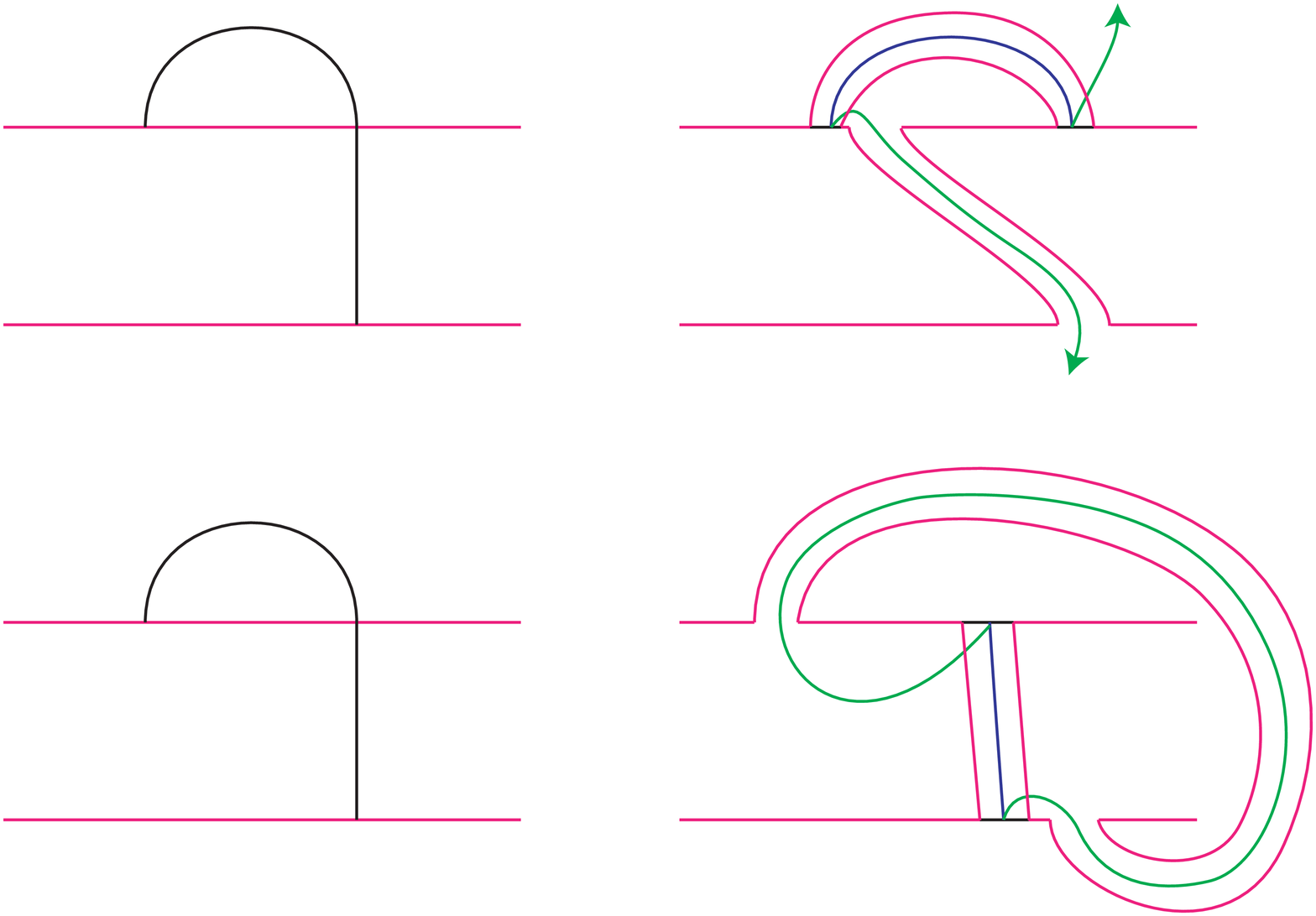}
\put(17.7,62){\tiny $+$} \put(17.7,52){\tiny $-$}
\put(17.7,40.8){\tiny $+$} \put(17.7,24){\tiny $-$}
\put(17.7,13){\tiny $+$} \put(17.7,3.5){\tiny $-$}
\end{overpic}
\caption{The top row is the $R_+$ case and the bottom row is the
$R_-$ case. The diagrams on the right-hand side depict the
$1$-handle attached to $S'$ to obtain $S$. The blue arc $a_0$
completes a basis for $(S',R_+(\Gamma'),h')$ to a basis
$(S,R_+(\Gamma),h)$, and the green arc is its image $h(a_0)$.}
\label{trivial-bypasses}
\end{figure}

We now apply the calculations done in \cite[Section~5,
Example~5]{HKM3} to obtain a description of $(S,R_+(\Gamma),h)$.
There are two cases of trivial bypasses: $c$ cuts off a half-disk
$D_1$ of $\bdry M'-\Gamma'$ which is either in $R_+(\Gamma')$ or in
$R_-(\Gamma')$. (If $c$ cuts off two half-disks $D_1$, $D_2$ and
$\bdry D_1$, $\bdry D_2$ intersect along an arc of $\Gamma'$, then
we take $D_1$ to be the ``smaller'' half-disk, i.e., $\bdry D_1\cap
\Gamma' \subset \bdry D_2\cap \Gamma'$.) The two cases will be
called the $R_+$ and $R_-$ cases, respectively. See
Figure~\ref{trivial-bypasses} for the determination of the monodromy
corresponding to the portion that is attached.

\subsubsection{Effect of a trivial bypass attachment on sutured Floer
homology} \label{subsub: triv}

Let $(\Sigma',\beta',\alpha')$ be the contact-compatible Heegaard
splitting for a basis $\{a_1',\dots,a_k'\}$ for
$(S',R_+(\Gamma'),h')$ and
$(\Sigma,\beta=\beta'\cup\{\beta_0\},\alpha=\alpha'\cup\{\alpha_0\})$
be its extension with respect to $\{a_1',\dots,a_k',a_0\}$ for
$(S,R_+(\Gamma),h)$. Here $\alpha_0=\bdry (a_0\times[0,1])$ and
$\beta_0=(b_0\times\{1\})\cup (h(b_0)\times\{1\})$, where $b_0$ is
the usual pushoff of $a_0$. Let $x_0$ be the $EH$ class
corresponding to $a_0$.

Let $c$ be the trivial arc of attachment along $\bdry M'$ and let
$D_1\subset \bdry M'$ be the half-disk cobounded by a subarc of $c$
and an arc of $\Gamma'$, as described above. Assume, without loss of
generality, that no endpoint of $K'$ lies on $\bdry D_1$. If
$D_1\subset R_+(\Gamma')$, then the only intersection of $\alpha_0$
with any $\beta_i$ is $x_0$.  On the other hand, if $D_1\subset
R_-(\Gamma')$, then the only intersection of $\beta_0$ with any
$\alpha_i$ is $x_0$. Therefore, for both $R_+$ and $R_-$, the
inclusion map
$$CF(\Sigma',\beta',\alpha')\rightarrow CF(\Sigma,\beta,\alpha),$$
$$\mathbf{y}\mapsto (\mathbf{y},x_0)$$
is an isomorphism of chain complexes. Therefore, tensoring with
$x_0$ gives an isomorphism
$$\Phi:\SFH(-M',-\Gamma')\stackrel\sim\rightarrow \SFH(-M,-\Gamma).$$

However, in order to show that the map is an identity morphism, we
need to decompose the stabilization (i.e., attaching a handle to
$\Sigma'$ and adding $\alpha_0$, $\beta_0$ to $\alpha'$, $\beta'$)
into a trivial stabilization and a sequence of handleslides.  Let us
consider the $R_-$ case.  (The $R_+$ case is left to the reader.) In
this case, $\beta_0$ only intersects $\alpha_0$, but $\alpha_0$ can
intersect $\beta_i'$. If there are no other intersections, then we
are done, since we have a trivial stabilization.  Otherwise,
consider the pushoff $\overline{a}_0$ of $a_0$, obtained by
isotoping the endpoints of $a_0$ along $\bdry S'$, against the
orientation of $\Gamma'$. If we stabilize $(S',R_+(\Gamma'),h')$
along $\overline{a}_0$, then all the intersections between
$\alpha_0$ and $\beta_i'$ will be eliminated, since the composition
with the positive Dehn twist forces the arcs to go around the core
of the attached $1$-handle. In its place, if $a_{k+1}'$ is the
cocore of the $1$-handle, then its image under the monodromy will
intersect $a_0$ exactly once. Let us rename open books and assume
$(S',R_+(\Gamma'),h')$ already has this property, namely we may
assume that there is only one intersection between $\alpha_0$ and
$\cup_i\beta_i'$.  The rest of the argument is identical to that of
Lemma~\ref{lemma:fig3}, and will be omitted.

\subsubsection{Reduction to a sequence of trivial bypasses}
Suppose now that $(M',\Gamma')$ is a sutured submanifold of
$(M,\Gamma)$, $M-M'=\bdry M'\times[0,1]$, $\bdry M'=\bdry
M'\times\{0\}$, and the contact structure $\xi$ on $\bdry
M'\times[0,1]$ with convex boundary condition $\Gamma\cup\Gamma'$ is
$[0,1]$-invariant. We now prove that there is an extension of
$(S',R_+(\Gamma'),h')$ for $(M',\Gamma')$ to $(S,R_+(\Gamma),h)$ for
$(M,\Gamma)$, of the type constructed in Step~2 of
Section~\ref{section: definition}, which can be decomposed into a
sequence of trivial bypass attachments. The nature of this extension
is such that it is obtained by adding ``horizontal'' Legendrian arcs
of type $\delta\times\{t\}\subset \bdry M'\times[0,1]$ and
``vertical'' Legendrian arcs of type $\{p\}\times[0,t]$. We will see
how the extension can be thought of as a sequence of trivial bypass
additions.

\begin{figure}[ht]
\begin{overpic}[width=6cm]{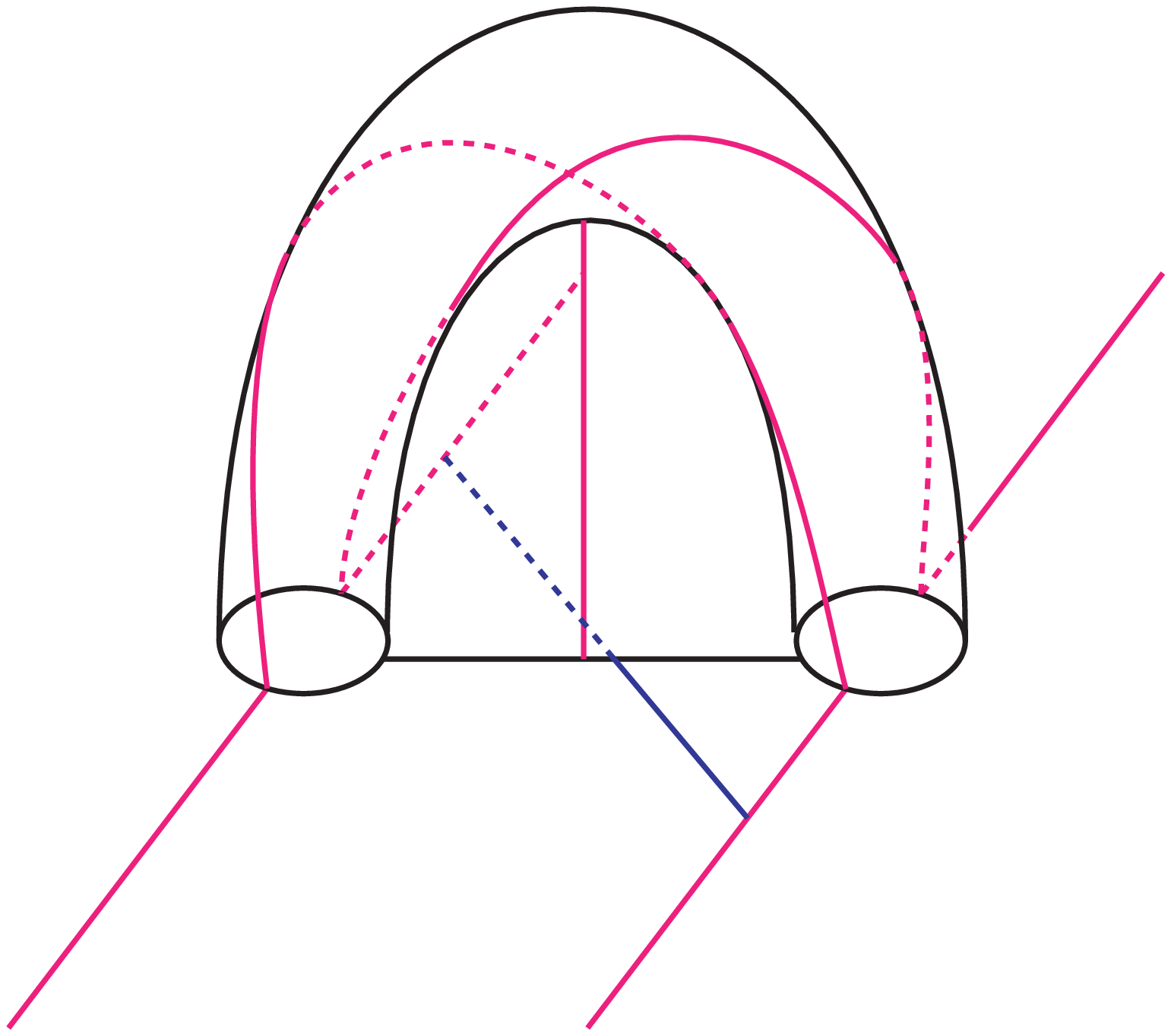}
\put(17,75){\tiny $N(d)$} \put(55,45){\tiny $D_0$}
\put(52.7,23){\tiny $\overline{a}_0$}
\end{overpic}
\caption{Attaching a trivial bypass arc $d$.  We stabilize along
$\overline{a}_0$ before attaching the bypass.} \label{attacharc}
\end{figure}

Observe that, when we attach a neighborhood $N(d)$ of a trivial
bypass arc $d$, then the result can be viewed more symmetrically as
in Figure~\ref{attacharc}. (This we leave as an exercise for the
reader.) This means that $d$ can be viewed as
the concatenation of three Legendrian arcs: two ``vertical'' arcs
$\{p_1,p_2\}\times[0,t]$ and a ``horizontal'' arc
$\delta\times\{t\}$, where $\delta$ connects $p_1$ and $p_2$. In
this subsection we make the assumption that all $\delta$'s, possibly
with subscripts, do not intersect $\Gamma'$ in the interior of
$\delta$, and all $p_i$'s are in $\Gamma'$.  Let $c'$ be the
component of $c-\Gamma'$ which is not part of $\bdry D_1$.  Slide
the endpoints of $c'$ in the direction of $\Gamma'$ if $c'\subset
R_-(\Gamma')$ and in the direction of $-\Gamma'$ if $c'\subset
R_+(\Gamma')$.  We will call the resulting Legendrian arc
$\overline{a}_0$; this notation agrees with the notation for the
stabilizing arc, used in Section~\ref{subsub: triv}.  If we take a
Legendrian-isotopic copy of $\overline{a}_0$ inside $M'$ via an
isotopy which fixes the endpoints, then we perform a stabilization
as in Section~\ref{subsub: triv} along the copy of $\overline{a}_0$
before attaching the bypass.

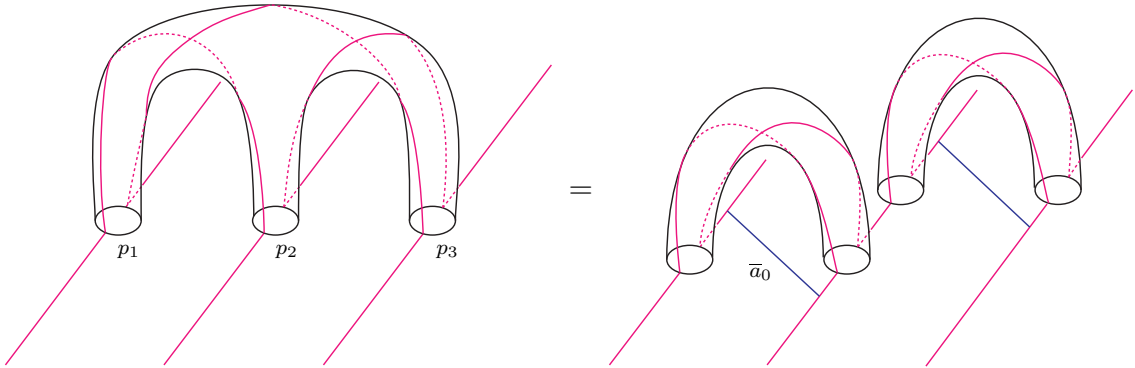
\begin{figure}[ht]
\s
\begin{overpic}[width=15cm]{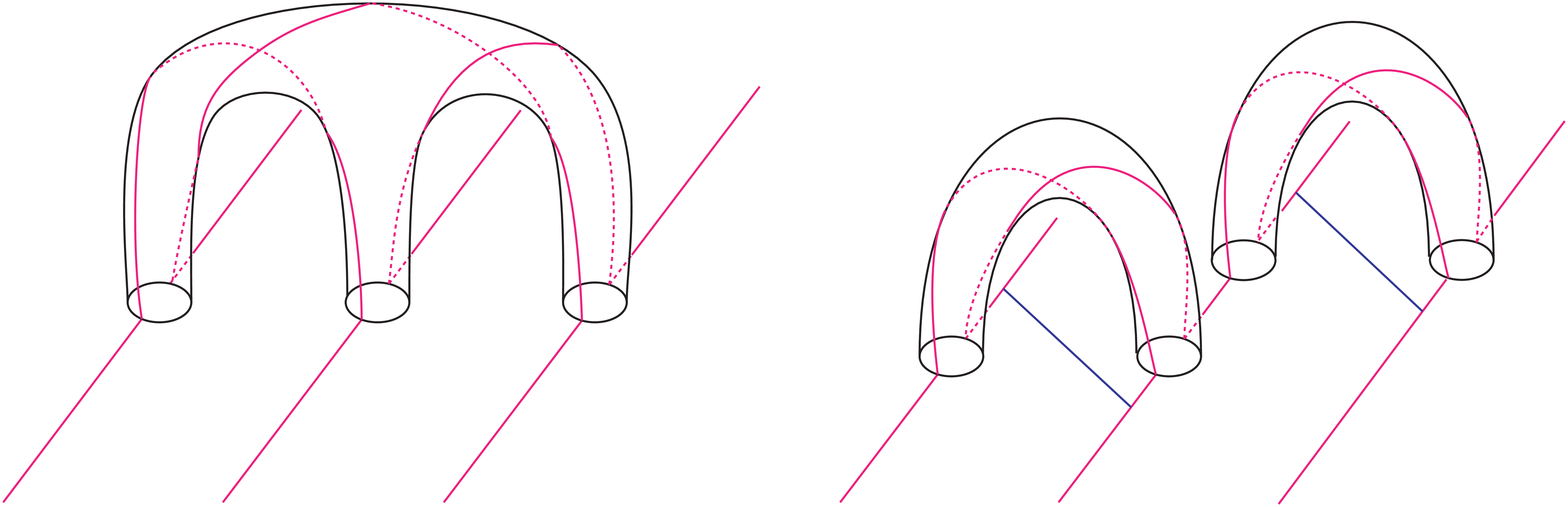}
\put(10,10){\tiny $p_1$} \put(38.25,10){\tiny $p_3$}
\put(24,10){\tiny $p_2$} \put(50,15){$=$} \put(66,8){\tiny
$\overline{a}_0$}
\end{overpic}
\caption{Sliding the bypass arc.} \label{twoarcs}
\end{figure}

Now, if we have a Legendrian graph consisting of
$\{p_1,p_2,p_3\}\times[0,t]$, together with $\delta_i\times\{t\}$,
$i=1,2$, with $\bdry \delta_i=\{p_i, p_{i+1}\}$, then attaching its
standard Legendrian neighborhood is equivalent to attaching two
bypass arcs as given in Figure~\ref{twoarcs}; this is readily seen
by sliding an endpoint of $\delta_2\times\{t\}$ along the dividing
set on the boundary of the union of $M'$ and the neighborhood of the
Legendrian arc $(\{p_1\}\times[0,t]) \cup (\delta_1\times\{t\}) \cup
(\{p_2\}\times[0,t])$.

\begin{figure}[ht]
\begin{overpic}[width=5.8cm]{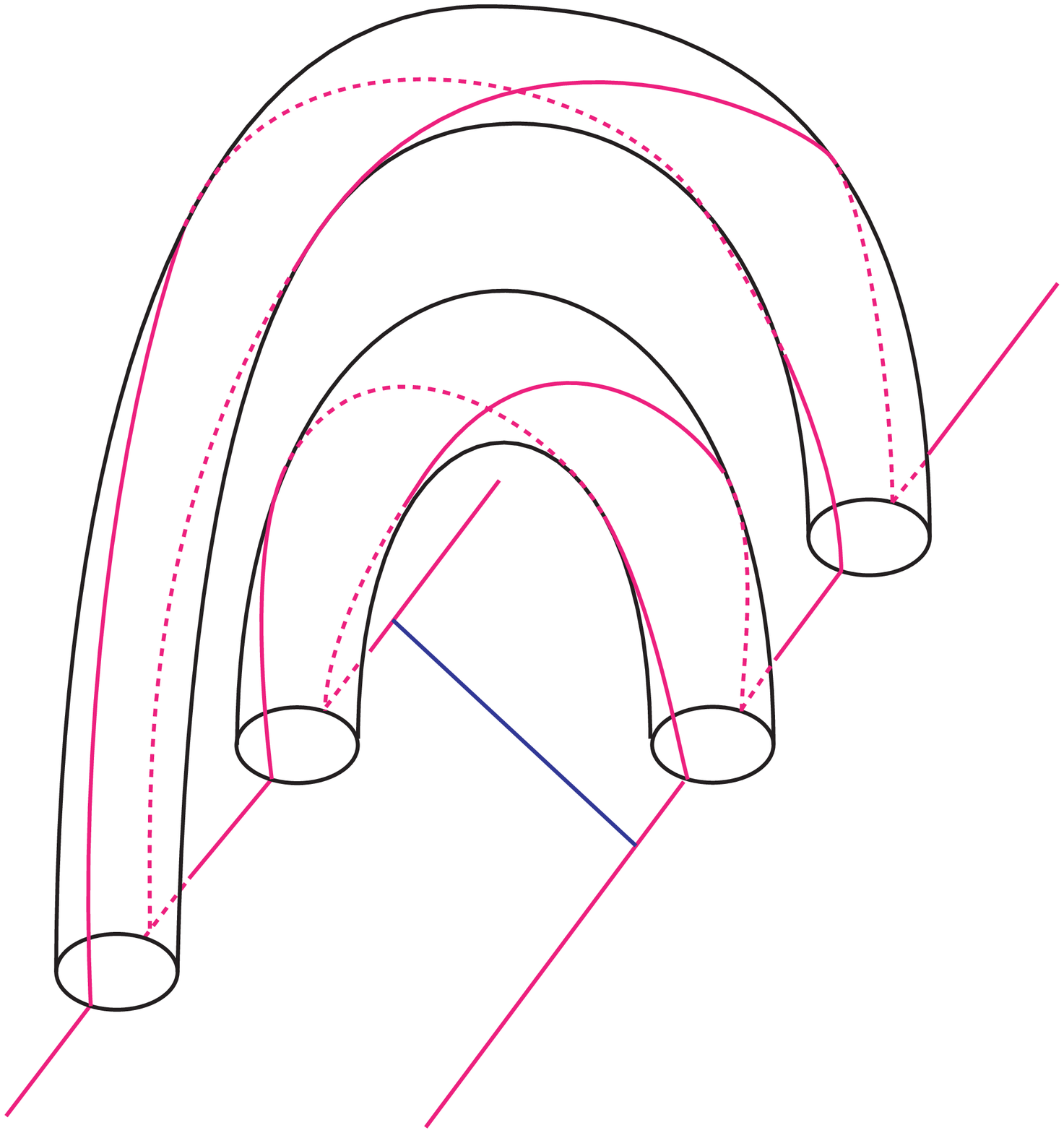}
\put(33,34){\tiny $p_1$} \put(70,34){\tiny $p_2$} \put(83,52){\tiny
$p_3$} \put(45,28){\tiny $\overline{a}_0$}
\end{overpic}
\caption{} \label{sliding}
\end{figure}

Finally, let $\gamma_1$ be a Legendrian arc given by the
concatenation of $\{p_1,p_2\}\times [0,t]$ and $\delta_1\times\{t\}$
with $\bdry \delta_1=\{p_1,p_2\}$, and we attach a Legendrian arc
$\gamma_2$ consisting of $\{p_3\}\times [0,t]$ and
$\delta_2\times\{t\}$ with $\bdry \delta_2=\{p_3,q\}$, where $q$ is
an interior point of $\delta_1$. By sliding the endpoint of
$\delta_2$, we see that attaching $\gamma_1$ and $\gamma_2$ is
equivalent to attaching the two Legendrian arcs given in
Figure~\ref{sliding}.  When attaching the first arc $\gamma_1$, we
first stabilize along $\overline{a}_0$; for the second arc
$\gamma_2$, the attachment of the first arc has the same effect as a
stabilization.

Therefore, using the above trivial bypass arcs, we can construct a
Legendrian graph $L$ in $\bdry M'\times[0,1-\varepsilon]$, which is
the union of arcs of type $\{p\}\times [0,1-\varepsilon]$ and the
$1$-skeleton $L_{1-\varepsilon}$ of a cell decomposition of $\bdry
M'\times \{1-\varepsilon\}$, each of whose cells have boundary with
$tb=-1$. (Here $(p,1-\varepsilon)$ must lie in $L_{1-\varepsilon}$.)
If we take the standard Legendrian neighborhood $N(L)$, then its
complement $(\bdry M'\times[0,1-\varepsilon]) - N(L)$ is also a
standard neighborhood of a Legendrian graph $K''$. The Legendrian
neighborhood $N(L_{1-\varepsilon})$ can be enlarged via a contact
isotopy so that $N(L_{1-\varepsilon})$ is $\bdry
M'\times[1-2\varepsilon,1]$, with neighborhoods of Legendrian arcs
of type $\{q\}\times[1-2\varepsilon,1]$, $q\in \Gamma'$, removed. On
the other hand, $N(\{p\}\times [0,1-\varepsilon])$ is viewed as a
sufficiently thin/small neighborhood of the Legendrian arc
$\{p\}\times [0,1-\varepsilon]$. The above description clearly
indicates that the extension of $K'$ to $K'\cup K''$ is an extension
of the partial open book decomposition $(S',R_+(\Gamma'),h')$ to
$(S,R_+(\Gamma),h)$ of the type described in Step 2 of
Section~\ref{section: definition}. This completes the proof of
Theorem~\ref{thm: identity}.

\section{A $(1+1)$-dimensional TQFT} \label{section: TQFT}

In this section we describe a $(1+1)$-dimensional TQFT, obtained by
dimensional reduction. (Strictly speaking, the theory does not quite
satisfy the TQFT axioms but has similar composition rules.)

\subsection{Invariants of multicurves on surfaces}

In this subsection we describe a TQFT-type invariant of a multicurve
on a surface. Let $\Sigma$ be a compact, oriented surface with
nonempty boundary $\bdry\Sigma$, and $F$ be a finite set of points
of $\bdry \Sigma$, where the restriction of $F$ to each component of
$\bdry \Sigma$ consists of an even number $\geq 2$ of points. Part
of the structure of a pair $(\Sigma, F)$ is a labeling of each
component of $\partial\Sigma - F$ by $+$ or $-$ so that crossing a
point of $F$ while moving along $\partial \Sigma$ reverses signs.
Let $\#F=2n$ be the cardinality of $F$. Also let $K$ be a properly
embedded, oriented $1$-dimensional submanifold of $\Sigma$ whose
boundary is $F$ and which divides $\Sigma$ into $R_+$ and $R_-$ in a
manner compatible with the labeling of $\bdry \Sigma-F$. As on
$\partial\Sigma - F$, the sign changes every time $K$ is crossed.
Such a $K$ will be called a {\em dividing set for $(\Sigma,F)$}.

We now list the properties satisfied by the TQFT.

\s\s\s\s\s\s\n {\bf TQFT Properties.}

\s\s
\begin{enumerate}
\item[I.]  \label{condition1} It assigns to each $(\Sigma,F)$ a
graded $\Z$-module $V(\Sigma,F)$. If $\Sigma$ is connected, then
$$V(\Sigma,F)=\Z^2\otimes\dots\otimes \Z^2,$$
where the number of copies of $\Z^2$ is $r=n-\chi(\Sigma)$, and
$\Z^2=\Z\oplus \Z$ is a graded $\Z$-module whose first summand has
grading $1$ and the second summand has grading $-1$. We will refer
to this grading as the {\em Spin$^c$-grading}. Moreover, if
$(\Sigma,F)$ is the disjoint union $(\Sigma_1,F_1)\sqcup
(\Sigma_2,F_2)$, then $$V(\Sigma_1\sqcup \Sigma_2,F_1\sqcup
F_2)\simeq V(\Sigma_1,F_1)\otimes V(\Sigma_2,F_2).$$

\s
\item[II.] To each $K$ it assigns a subset $c(K)\subset V(\Sigma,F)$
of cardinality $1$ or $2$ of type $\{\pm x\}$, where $x\in
V(\Sigma,F)$. If $K$ has a homotopically trivial closed component,
then $c(K)=\{0\}$.

\s
\item[III.] \label{condition3}
Given $(\Sigma,F)$, let $\gamma,\gamma'\subset \bdry\Sigma$ be
mutually disjoint $1$-dimensional submanifolds of $\bdry\Sigma$, so
that their endpoints do not lie in $F$. Suppose there is a
diffeomorphism $\tau:\gamma\stackrel\sim\rightarrow \gamma'$ which
sends $\gamma\cap F\stackrel\sim\rightarrow\gamma'\cap F$  and
preserves the orientations on $\gamma \cap \partial\Sigma$ and
$\gamma' \cap \partial \Sigma'$. If we glue $(\Sigma,F)$ by
identifying $\gamma$ and $\gamma'$ via $\tau$, then the result will
be denoted by $(\Sigma',F')$. Then there exists a map
$$\Phi_{\tau}:V(\Sigma,F)\rightarrow V(\Sigma',F'),$$
which satisfies $$c(K)\mapsto c(\overline{K}),$$ where
$\overline{K}$ is obtained from $K$ by gluing $K|_{\gamma}$ and
$K|_{\gamma'}$. Here $\Phi_{\tau}$ is well-defined up to an overall
$\pm 1$ multiplication.
\end{enumerate}

\begin{figure}[ht]
\begin{overpic}[width=7.5cm]{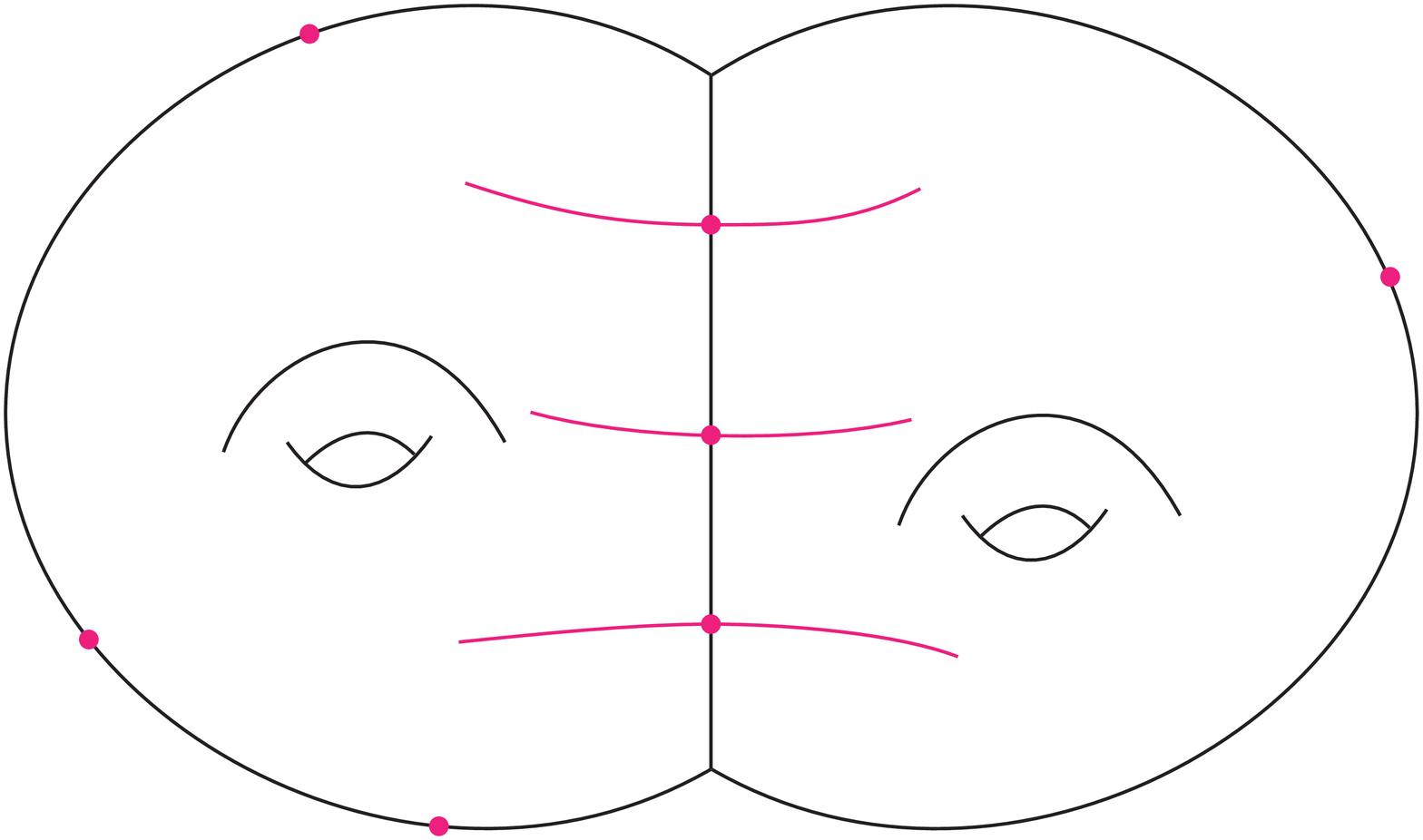}
\put(20,44){\small $\Sigma''$} \put(75,44){\small $\Sigma'''$}
\put(40,16.5){\tiny $K''$} \put(60,15.8){\tiny $K'''$}
\end{overpic}
\caption{Gluing $(\Sigma'',K'')$ and $(\Sigma''',K''')$. The red
dots are $F''$ and $F'''$.} \label{tensor}
\end{figure}
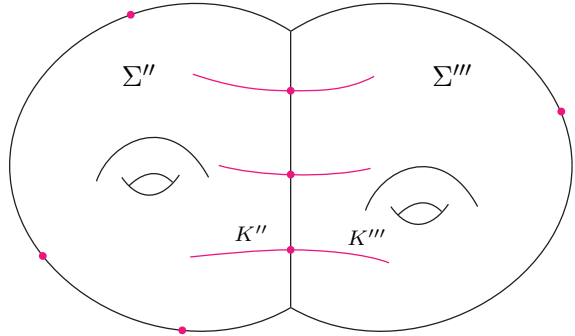

\s\s See Figure~\ref{tensor} for an illustration of the gluing in
Property III, when $(\Sigma,F)=(\Sigma'',F'')\sqcup
(\Sigma''',F''')$, $K=K''\sqcup K'''$, and $\gamma,\gamma'$ are
arcs. In this case, the gluing map is:
$$\Phi_\tau: V(\Sigma'',F'')\otimes V(\Sigma''',F''')\rightarrow
V(\Sigma',F').$$

We will use the subscripts $(i)$ to denote the Spin$^c$-grading:
$V(\Sigma,F)_{(i)}$ is the graded piece with grading $i$ and
$\Z^n_{(i)}$ is the $\Z^n$-summand representing the $i$th graded
piece.

\begin{thm}\label{thm: existence}
There exists a nontrivial TQFT satisfying Properties I-III above.
\end{thm}

\begin{proof}
This TQFT arises by dimensional reduction of sutured Floer homology.

\s\n I. ~~Given $(\Sigma,F)$, let $F_0\subset \bdry \Sigma$ be
obtained from $F$ by shifting slightly in the direction of
$\bdry\Sigma$.  (We may think of points of $F_0$ as being situated
halfway between successive points of $F$ on $\bdry \Sigma$.) We
consider $S^1$-invariant balanced sutured manifolds $(S^1\times
\Sigma, S^1\times F_0)$, and let
$$V(\Sigma,F)= SFH(-(S^1\times \Sigma),-(S^1\times F_0)).$$
The reason for using $F_0$ instead of $F$ in the definition is
explained below in II when the role of contact structures is
explained.  The Spin$^c$-grading for $V(\Sigma,F)$ corresponds to
the different relative Spin$^c$-structures on $(S^1\times
\Sigma,S^1\times F_0)$.

The next lemma determines $V(\Sigma,F)$, up to isomorphism.

\begin{lemma}
If $\Sigma$ is connected, then $$\SFH(-(S^1\times
\Sigma),-(S^1\times F_0))\simeq (\Z_{(-1)}\oplus \Z_{(1)})^{\otimes
r},$$ where $r= n-\chi(\Sigma)$.
\end{lemma}

\begin{proof}
This follows from Juh\'asz' tensor product
formula~\cite[Proposition~8.10]{Ju2} for splitting sutured manifolds
along product annuli, together with the observation that when $n=2$
and $\Sigma=D^2$, we have $SFH(-(S^1\times D^2),-(S^1\times
F_0))\simeq \Z_{(-1)}\oplus \Z_{(1)}$, split according to the
relative Spin$^{c}$-structure.  (See \cite[Section~5,
Example~2]{HKM3}.)
\end{proof}

Finally, the property
$$V(\Sigma_1\sqcup \Sigma_2,F_1\sqcup
F_2)\simeq V(\Sigma_1,F_1)\otimes V(\Sigma_2,F_2)$$ is immediate
from the definition of the sutured Floer homology groups.

\s\n II. ~~Next, there is a $1-1$ correspondence between dividing
sets $K$ of $(\Sigma,F)$ without homotopically trivial closed curves
and tight contact structures with boundary condition
$(S^1\times\Sigma,S^1\times F_0)$, up to isotopy rel boundary. For
the correspondence to hold we require that
$\bdry\Sigma\not=\emptyset$. The map from dividing sets to contact
structures is easy: simply consider the $S^1$-invariant contact
structure $\xi_K$ on $S^1\times \Sigma$ so that the dividing set on
each $\{pt\}\times \Sigma$ is $\{pt\}\times K$.  It was shown in
\cite{Gi3,H2} that the map, when restricted to the subset of
dividing sets $K$ without homotopically trivial curves, gives a
bijection with the set of isotopy classes of tight contact
structures on $(S^1\times \Sigma,S^1\times F_0)$. Now, to each $K$
we assign $\EH(\xi_K)\subset SFH(-(S^1\times \Sigma),-(S^1\times
F_0))$. If $K$ has a homotopically trivial curve, then $\xi_K$ is
overtwisted, and $\EH(\xi_K)=\{0\}$.

Finally we explain why we use $F_0$ instead of $F$ in $(S^1\times
\Sigma,S^1\times F_0)$. The dividing set $S^1\times F_0$ of
$\bdry(S^1\times \Sigma)$ does not intersect the dividing set of
$\{pt\}\times\Sigma$, since the two surfaces are transverse. This
means that $F_0$ must lie between the endpoints $F$ of $K$.

\s\n III. ~~This is a corollary of Theorem~\ref{thm: gluing}: In
order to apply Theorem~\ref{thm: gluing}, slightly shrink $\Sigma$
to $\Sigma_0$ inside the glued-up surface $\Sigma'$. See
Figure~\ref{tensor2}. If we write
$\Sigma-\Sigma_0=\bdry\Sigma\times[0,1]$ with
$\bdry\Sigma\times\{0\}=\bdry\Sigma_0$ and
$\bdry\Sigma\times\{1\}=\bdry\Sigma$, then the dividing set $K_0$ on
$\Sigma'-\Sigma_0$ is obtained from $F\times[0,1]$ by identifying
$F|_{\gamma}\times\{1\}$ with $F|_{\gamma'}\times\{1\}$ via $\phi$.
Let $\xi_{K_0}$ be the $S^1$-invariant contact structure on
$S^1\times (\Sigma'-\Sigma_0)$ corresponding to the dividing set
$K_0$.
\begin{figure}[ht]
\begin{overpic}[width=6cm]{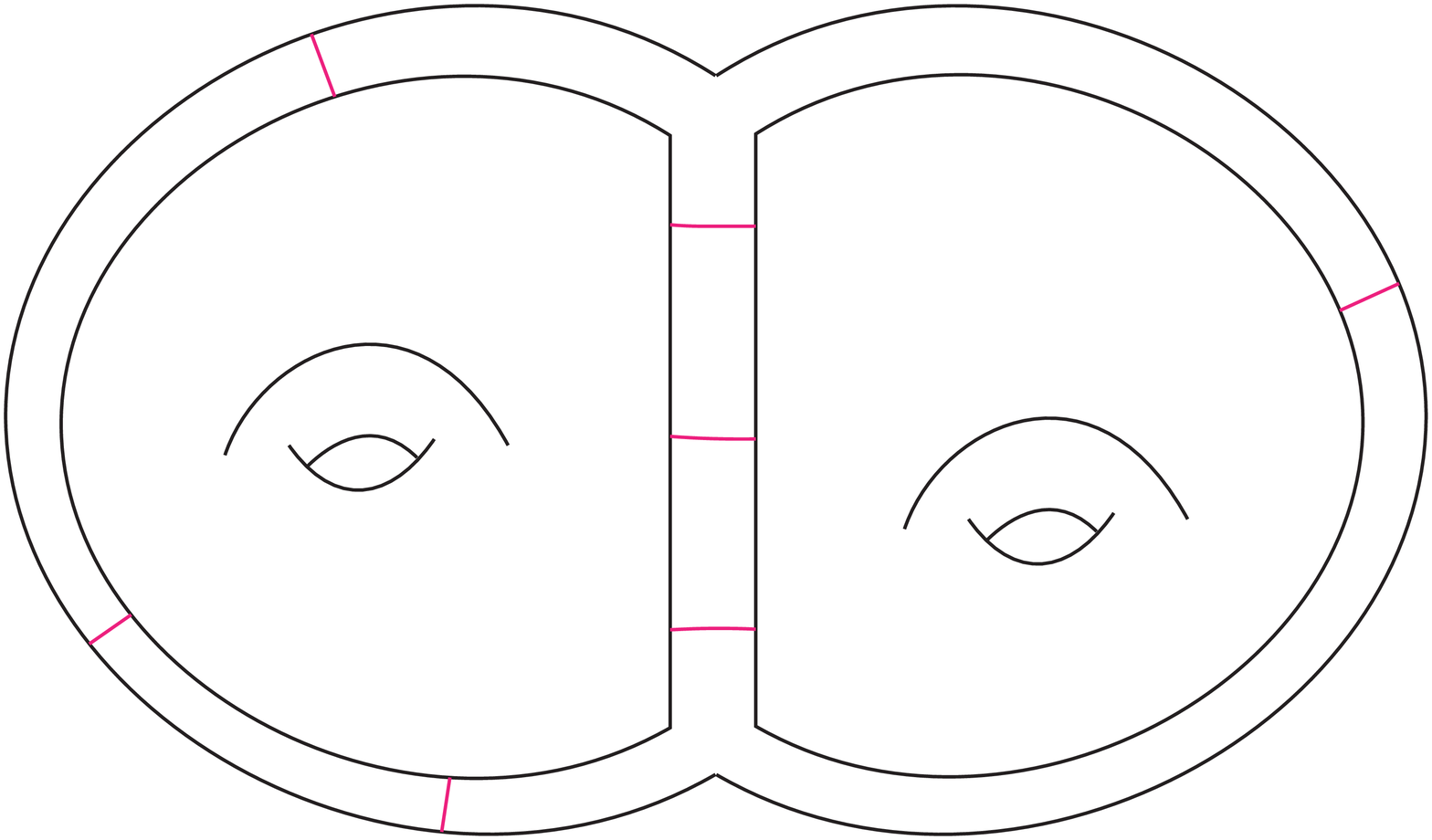}
\put(25,40){\small $\Sigma_0$} \put(70,40){\small $\Sigma_0$}
\end{overpic}
\caption{The dividing set $K_0$ on $\Sigma'-\Sigma_0$ is given in
red.} \label{tensor2}
\end{figure}
The contact structure $\xi_{K_0}$ induces the map
$\Phi_{\xi_{K_0}}=\Phi_\tau$ from $V(\Sigma,F)$ to $V(\Sigma',F')$.
This completes the proof of Theorem~\ref{thm: existence}.
\end{proof}

\begin{rmk}
There is another grading for $V(\Sigma,F)$, a relative grading
called the {\em Maslov grading}, which is largely invisible for the
time being since all the generators have the same Maslov grading.
\end{rmk}

\subsection{Analysis of $\Sigma=D^2$}
Suppose $\Sigma=D^2$ and $F$ consists of $2n$ points on $\bdry D^2$.
In this case, the set of dividing sets $K$ without closed components
corresponds to the set of crossingless matchings of $F$. A {\em
crossingless matching of $F$} is a collection of $n$ properly
embedded arcs in $D^2$ with endpoints on $F$ so that each endpoint
is used once and no two arcs intersect in $D^2$. The orientation
condition is trivially satisfied for a crossingless matching. If $K$
has a closed component, then the component must be homotopically
trivial. Thus the corresponding contact structure $\xi_K$ is
overtwisted, and $c(K)=\{0\}$.

\s\n {\bf n=1.} When $n=1$, $V(\Sigma,F)=\Z_{(0)}$, which is
generated by the unique $K$ which connects the two points. (By this
we mean $\Z$ is generated by either element of $c(K)$.)

\begin{figure}[ht]
\begin{overpic}[width=7.5cm]{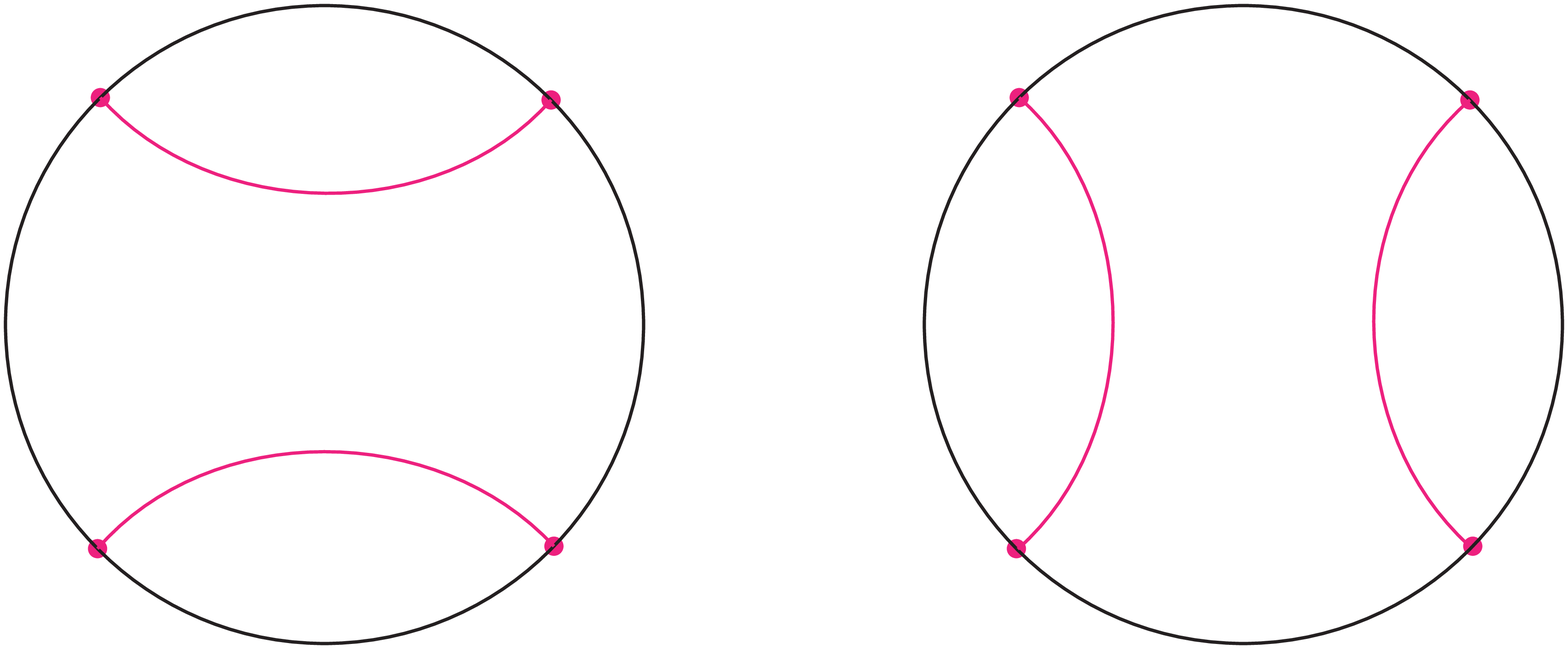}
\put(18.7,20){\tiny $-$} \put(18.2,-5.7){\tiny $K_+$}
\put(78.2,-5.7){\tiny $K_-$} \put(18.7,5){\tiny $+$}
\put(18.7,35){\tiny $+$} \put(78.2,20){\tiny $+$} \put(64,20){\tiny
$-$} \put(93,20){\tiny $-$}
\end{overpic}
\s \caption{} \label{baseballs}
\end{figure}

\s\n {\bf n=2.} When $n=2$, $V(\Sigma,F)=\Z_{(1)}\oplus \Z_{(-1)}$.
We claim that $V(\Sigma,F)$ is generated by $c(K_+)$ and $c(K_-)$,
given as in Figure~\ref{baseballs}.
Here $K_+$ and $K_-$ are the two dividing sets, both
$\bdry$-parallel. The grading for $c(K)$ can be calculated by taking
$\chi(R_+)-\chi(R_-)$, where $R_+$ (resp.\ $R_-$) is the positive
(resp.\ negative) region of $\Sigma-K$. Hence the degrees are $1$
and $-1$ for $K_+$ and $K_-$, respectively. As calculated in
\cite[Section~5, Example~3]{HKM3}, there is a Heegaard diagram for
which the $EH$ class for $\xi_{K_+}$ is the unique tuple
representing its Spin$^c$-structure (and similarly for $K_-$). Hence
$c(K_+)$ generates the first summand and $c(K_-)$ generates the
second summand, with respect to any coefficient system.

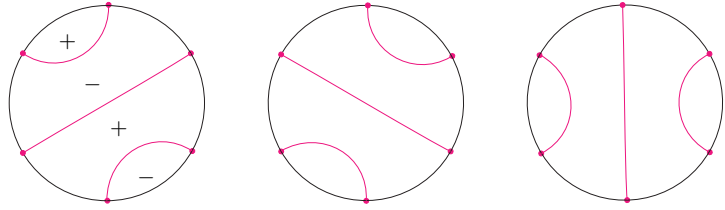
\begin{figure}[ht]
\begin{overpic}[width=9.5cm]{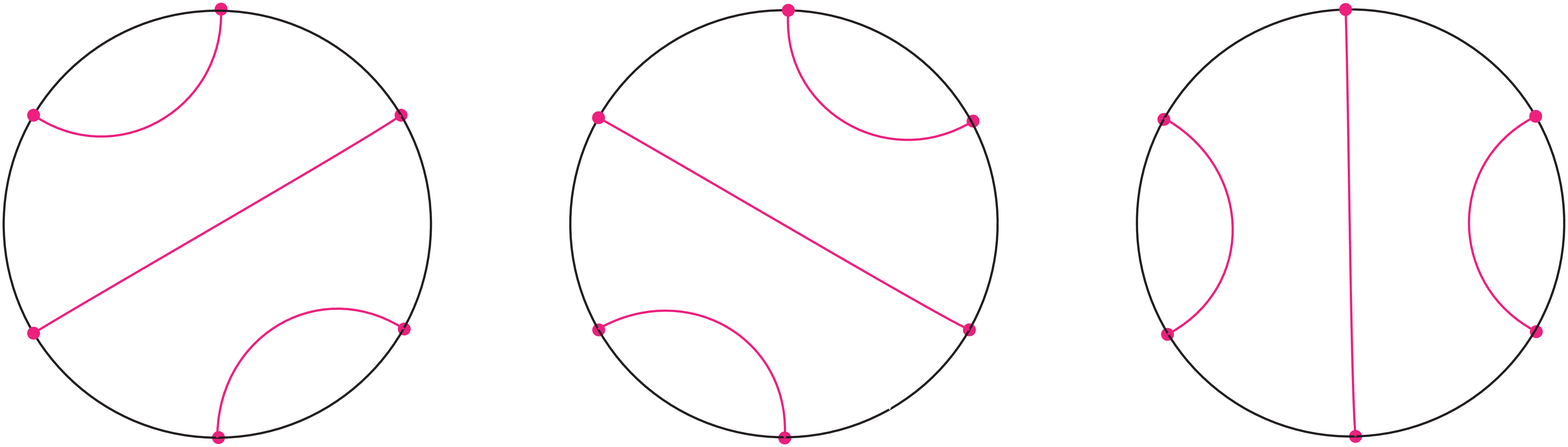}
\put(7,22){\tiny $+$} \put(10.5,16){\tiny $-$} \put(14.2,9.8){\tiny
$+$} \put(18,3){\tiny $-$}
\end{overpic}
\caption{$K_1, K_2, K_3$, from left to right.} \label{nthree}
\end{figure}

\s\n {\bf n=3.} When $n=3$, $V(\Sigma,F)$ decomposes into
$\Z_{(2)}\oplus \Z^2_{(0)} \oplus \Z_{(-2)}$. The first and last
summands are generated by $c(K)$ for $\bdry$-parallel $K$. The
middle $\Z^2_{(0)}$ must support three configurations $K_1, K_2,
K_3$. See Figure~\ref{nthree}.

We have the following:

\begin{lemma}\label{lemma: n is three}
The sets $c(K_1), c(K_2), c(K_3)$ are nonzero and distinct.
Moreover, their elements are primitive.
\end{lemma}

Over $\Z/2\Z$, the lemma implies that $c(K_1)+c(K_2)=c(K_3)$, i.e.,
$c(K_3)$ is a superposition of the other two states $c(K_1)$ and
$c(K_2)$.

\begin{proof}
Consider an arc $\gamma\subset \bdry \Sigma$ with
$\#(F\cap\gamma)=2$. Take a disk $\Sigma'=D^2$ with $\#F'=2$, and
pick an arc $\gamma'\subset \bdry\Sigma'$ with $\#(F'\cap
\gamma')=2$.  Then attach $\Sigma'$ onto $\Sigma$ so that $\gamma$
and $\gamma'$ are identified and $F''=F\cup F'-\gamma$ satisfies
$\#F''=4$. Observe that the $\Z$-module $V(\Sigma',F')\simeq \Z$ is
generated by a unique element $K'$, which is a $\bdry$-parallel arc.
Label the points of $F$ in clockwise order from $1$ to $6$, so that
$1$ is $2$pm, $2$ is $4$pm, etc., and let $\Phi_j$, $j=1,2,3$, be
the gluing map
$$V(\Sigma,F)\rightarrow V(\Sigma\cup\Sigma',F''),$$
obtained by attaching the $\bdry$-parallel arc $K'$ from $j$ to
$j+1$. It sends $c(K_i)\mapsto c(K_i\cup K')$. See
Figure~\ref{attach-bdry-parallel}. Restricted to $\Z^2_{(0)}$, the
image of $\Phi_j$ is one of the two summands $\Z_{(1)}$ or
$\Z_{(-1)}$.  Hence we view $\Phi_j$ as a map $\Z^2_{(0)}\rightarrow
\Z_{(\pm 1)}$.

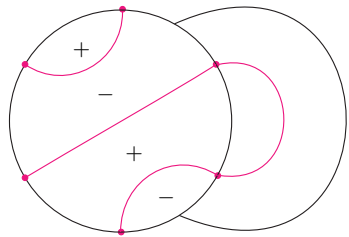
\begin{figure}[ht]
\begin{overpic}[width=4.5cm]{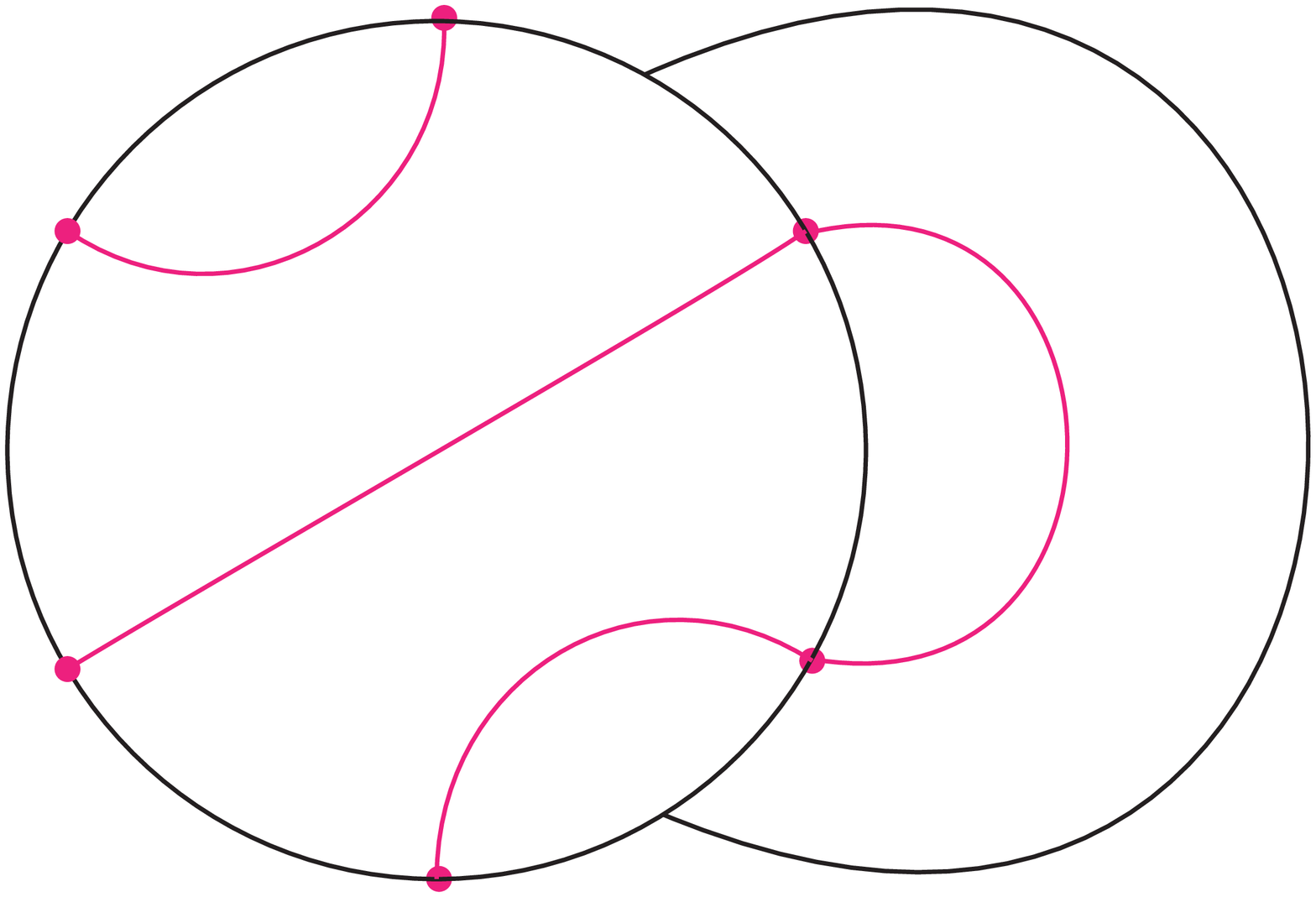}
\put(19,53.5){\tiny $+$} \put(26,40){\tiny $-$}
\put(34.5,22.7){\tiny $+$} \put(44,9.5){\tiny $-$}
\end{overpic}
\caption{The diagram represents $\Phi_1(c(K_1))$. $K_1$ is the
crossingless matching inside the circle, and the $\bdry$-parallel
arc, representing a generator of $V(\Sigma',F')=\Z$ with $\#F'=2$,
is attached from $1$ to $2$.} \label{attach-bdry-parallel}
\end{figure}

Suppose $K_i\cup K'$ does not have a closed component; there is
always some $\Phi_j$ for which this is true.  Then we have reduced
to the case $n=2$, where we already know that each representative of
$c(K_i\cup K')$ is nonzero and primitive. Since
$\Phi_j:\Z^2\rightarrow \Z$ maps $c(K_i)\mapsto c(K_i\cup K')$ and
the latter is primitive, it follows that $c(K_i)$ must also be
primitive.

Next, $c(K_i\cup K')=\EH(\xi_{K_i\cup K'})=\{0\}$ if $K_i\cup K'$
has a closed (and necessarily homotopically trivial) component.
Hence, by attaching $\Sigma'$ at the appropriate locations (i.e.,
checking which $\Phi_1$, $\Phi_2$ or $\Phi_3$ annihilates $c(K_i)$),
we can determine the locations of all the $\bdry$-parallel (or
outermost) arcs of $K_i\subset \Sigma$. Since the location of the
$\bdry$-parallel arcs determines $K_i$, it follows that the $c(K_i)$
must be distinct.
\end{proof}

By inductively applying the above procedure, we obtain the
following:

\begin{prop}\label{prop: crossingless}
All crossingless matchings $K$ of $(\Sigma=D^2,F)$ with $\#F=2n$ are
distinguished by $c(K)\subset V(\Sigma,F)$ and are primitive.
Equivalently, all the tight contact structures on $S^1\times D^2$,
$\#F=2n$, are distinguished by their contact invariant in
$SFH(-(S^1\times D^2),-(S^1\times F_0))$.
\end{prop}

The proof is left to the reader. Lemma~\ref{lemma: n is three} and
its generalization Proposition~\ref{prop: crossingless} are rather
surprising, since the dimension of $V(D^2,F)$ with $\#F=2n$ is
$2^{n-1}$, whereas the number of crossingless matchings on $(D^2,F)$
is the {\em Catalan number} $\displaystyle C_n={1\over n+1}\begin{pmatrix} 2n \\
n\end{pmatrix}$, which is greater than or equal to $2^{n-1}$, and
grows roughly twice as fast as a function of $n$. This means that
all the $c(K)$'s are ``tightly packed'' inside $V(D^2,F)$,
especially when the coefficient ring is $\Z/2\Z$.

Also recall that the dimensions of our $V(D^2,F)$ with $\#F=2n$ are
the same as that of $(1+1)$-dimensional, level $k=2$, $sl(2,\C)$
TQFT. It would be interesting to compare the two TQFT's.

\subsection{The $\pm 1$ ambiguity over $\Z$-coefficients} \label{subsection: Z}

In this subsection we prove Theorem~\ref{thm: unique lifting} and
deduce from it that the $\pm 1$ ambiguity of the contact invariant
$\EH(M,\Gamma,\xi)$ in $\SFH(-M,-\Gamma)$ over $\Z$ {\em cannot be
removed} and that the gluing map $\Phi_\xi$ of Theorem~\ref{thm:
gluing} is well-defined {\em only up to} an overall $\pm$ sign over
$\Z$.  This proves Theorem~\ref{thm: no single value}, stated in the
Introduction.

\begin{thm} \label{thm: unique lifting}
There is no single-valued lift of $c(K)\subset V(\Sigma,F)$ for all
$K, \Sigma, F$, with $\Z$-coefficients.
\end{thm}

\begin{proof}
Assume the invariants of dividing curves are single-valued. Consider
the example $\Sigma=D^2$ and $\#F=6$. Recall $K_1, K_2, K_3$ from
Figure~\ref{nthree}. By Lemma~\ref{lemma: n is three}, each element
of $c(K_i)$, $i=1,2,3$, is primitive in $\Z^2_{(0)}$. We also use
the same maps $\Phi_j: \Z^2_{(0)}\rightarrow \Z^2$, $j=1,2,3$.

We compute the following:
$$\Phi_1: c(K_1)\mapsto c(K_+), ~ c(K_2)\mapsto c(K_+), ~ c(K_3)\mapsto
0,$$
$$\Phi_2: c(K_1)\mapsto 0, ~ c(K_2)\mapsto c(K_-), ~ c(K_3)\mapsto
c(K_-),$$
$$\Phi_3: c(K_1)\mapsto c(K_+), ~ c(K_2)\mapsto 0, ~
c(K_3)\mapsto c(K_+).$$ Here $c(K_+)$ and $c(K_-)$ are generators of
$V(D^2,F')$ with $\#F'=4$. Since the image of each $\Phi_j(\Z^2)$ is
$\Z$, generated by either $c(K_+)$ or $c(K_-)$, we view $\Phi_j$ as
a map $\Z^2\rightarrow \Z$.

Let us analyze $\Phi_1$ in more detail.  Write $c(K_1)$ as $(1,0)\in
\Z^2$, since it is primitive. Then $\Phi_1: \Z^2\rightarrow \Z$ maps
$(1,0)\mapsto 1$. We can then decompose $\Z^2$ into $\Z\oplus \Z$ so
that $(0,1)$ generates $\ker \Phi_1$, possibly after an appropriate
isomorphism of $\Z^2$. Without loss of generality, $c(K_3)=(0,1)$.
Since $\Phi_1:c(K_2)\mapsto 1$, it follows that $c(K_2)=(1,a)$,
$a\in \Z$.

Next consider $\Phi_2$.  Since $(1,0)\mapsto 0$, $(0,1)\mapsto 1$,
and $(1,a)\mapsto 1$, it follows that $a=1$.

Finally, $\Phi_3$ maps $(1,0)\mapsto 1$, $(0,1)\mapsto 1$, and
should map $(1,1)\mapsto 2$, but instead sends it to $0$, a
contradiction.
\end{proof}

\begin{thm} \label{thm: ambiguity}
$\pm 1$ monodromy exists in sutured Floer homology. That is, there
is a sequence of stabilization, destabilization, handleslide, and
isotopy maps which begins and ends at the same configuration, so
that their composition is $-id$.
\end{thm}

\begin{proof}
When working over $\Z$,  $\EH(M,\Gamma,\xi)$ and $\Phi_\xi$ are
defined up to a factor of $\pm 1$.  The only reason for the
introduction of this factor was the {\em possibility} of the
existence of $\pm 1$ monodromy.  Since single-valued lifts do not
always exist by Theorem~\ref{thm: unique lifting}, it follows that
$\pm 1$ monodromy must exist.
\end{proof}

This proof is unsatisfying in the sense that it does not explain the
root cause of the existence of monodromy, nor does it give a
specific sequence of maps which exhibits nontrivial monodromy.

\begin{q}
Is there monodromy in Heegaard Floer homology, i.e., when the
$3$-manifold is closed? In particular, is there an explicit sequence
of stabilization, destabilization, handleslide, and isotopy maps
which begins and ends at the same configuration, so that their
composition is $-id$ for $\widehat{HF}(S^1\times S^2)\simeq \Z\oplus
\Z$?
\end{q}

\subsection{A useful gluing isomorphism}
In this subsection we give a useful gluing map and explore some
consequences.

Let $\gamma$ be a properly embedded arc on $\Sigma$ which is
transverse to $K$ and intersects $K$ exactly once. Suppose we cut
$(\Sigma,F)$ and $K$ along $\gamma$ to obtain $(\Sigma',F')$ and
$K'$. This is the reverse procedure of gluing $(\Sigma',F')$ and
$K'$ along disjoint subarcs $\gamma', \gamma''\subset \bdry
\Sigma'$, where each arc intersects $F'$ exactly once.  We then
have:

\begin{lemma} \label{lemma: cutting}
The gluing map $\Phi: V(\Sigma',F')\rightarrow V(\Sigma,F)$ is an
isomorphism.
\end{lemma}

If $\gamma$ decomposes $\Sigma$ into two components $(\Sigma'',F'')$
and $(\Sigma''',F''')$, then the gluing map is: $$\Phi:
V(\Sigma'',F'')\otimes V(\Sigma''',F''')\stackrel\sim\rightarrow
V(\Sigma,F).$$

\begin{proof}
We interpret the gluing map $\Phi$ as a gluing map $\Phi_0:
V(\Sigma',F')\rightarrow V(\Sigma,F)$, where the gluing occurs along
a $\bdry$-parallel convex annulus $A$ as given in
Figure~\ref{annulus-gluing}.
\begin{figure}[ht]
\begin{overpic}[width=4.5cm]{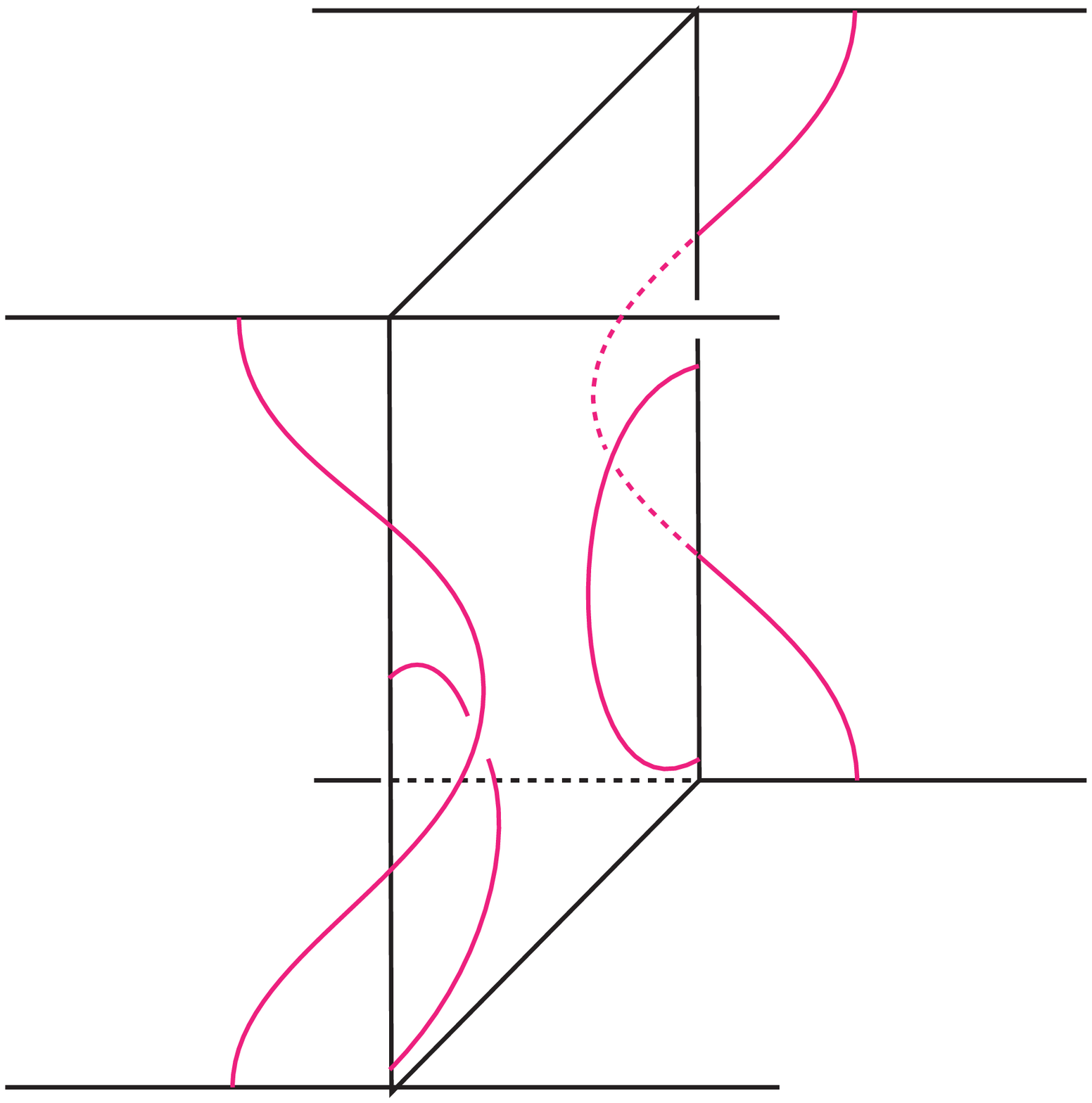}
\put(15,38){\tiny $+$} \put(37,16){\tiny $+$} \put(56,45){\tiny $+$}
\put(50,78.5){\tiny $A$}
\end{overpic}
\caption{The top and bottom of the annulus are identified.}
\label{annulus-gluing}
\end{figure}
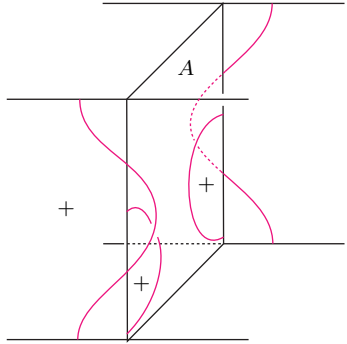

First we prove that $\Phi=\Phi_0$. Let
$(M,\Gamma)=(S^1\times\Sigma,S^1\times F_0)$, where $F_0$ is the
pushoff of $F$ in the direction of $\bdry\Sigma$.  Also let
$(M',\Gamma')$ be the sutured manifold obtained from $(S^1\times
\Sigma',S^1\times F_0')$ by slightly retracting $\Sigma'$ to
$\Sigma'_0$; here $F_0'$ is the pushoff of $F'$. Let $\xi_0$ be the
contact structure on $M-int(M')$, given as the union of the
invariant contact structures on a neighborhood of $\bdry M'$ with
dividing set $S^1\times F_0'$ and on a neighborhood of $A$ with
$\bdry$-parallel dividing set. Since the dividing set on $\bdry
(M-int(M'))$ is of the type $S^1\times \{\mbox{finite set}\}$,
$\xi_0$ is an $S^1$-invariant contact structure by \cite{Gi3,H2},
and is encoded by the ``minimal'' dividing set $K_0$ on
$\Sigma-\Sigma'_0$.

We now briefly sketch why $\xi_{K'}\cup \xi_0$ is isotopic to
$\xi_{K}$.  Let $\gamma_1,\gamma_2$ be the components of $\Gamma$
which intersect $\bdry A$ and $\delta_1,\dots,\delta_m$ be the
components of $\Gamma$ which do not intersect $\bdry A$.  For each
$\delta_i$, there is a parallel copy $\delta_i'$ on $\bdry M'$.
Moreover, there is a Legendrian arc from $\delta_i$ to $\delta_i'$
which has zero twisting number with respect to a surface parallel to
$\Sigma-\Sigma'_0$.  For each $\gamma_i$, there are two components
$\gamma_i'$ and $\gamma_i''$ on $\bdry M'$ which share a parallel
arc with $\gamma_i$.  Hence there are Legendrian arcs from
$\gamma_i$ to $\gamma_i'$ and from $\gamma_i$ to $\gamma_i''$ which
have zero twisting number as well.  The above Legendrian arcs
constrain $K_0$ so that $K'\cup K_0$ is isotopic to $K$. This proves
$\Phi=\Phi_0$.

We next prove that $\Phi_0$ is an isomorphism. According to
Juh\'asz~\cite{Ju2}, gluing along a product annulus gives an
isomorphism of sutured Floer homology groups.  Although our
situation is slightly different, the result is the same. By
\cite[Theorem~6.2]{HKM3}, $V(\Sigma',F')$ is a direct summand of
$V(\Sigma,F)$ since the dividing set on $A$ is $\bdry$-parallel.
Now, according to Proposition~\ref{prop: bdry-parallel}, $\Phi_0$ is
indeed the direct summand map of \cite[Theorem~6.2]{HKM3}.  More
precisely, $\Phi_0$ induces an isomorphism onto the Spin$^c$-direct
summand corresponding to the $\bdry$-parallel dividing set with
relative half-Euler class $\chi(R_+)-\chi(R_-)=2-0=2$.  To see that
$V(\Sigma',F')\simeq V(\Sigma,F)$ under the map $\Phi_0$, we use a
rank argument. Both $V(\Sigma',F')$ and $V(\Sigma,F)$ are isomorphic
to $\Z^r$, where $r={1\over 2}(\# F)-\chi(\Sigma)$.  Since
$V(\Sigma',F')$ is a direct summand of $V(\Sigma,F)$ and they are
both free with the same rank, it follows that $V(\Sigma',F')\simeq
V(\Sigma,F)$.
\end{proof}



As an application of Lemma~\ref{lemma: cutting}, we give a
sufficient condition for a dividing set $K$ for $(\Sigma,F)$ to have
$c(K)$ which is nonzero and primitive in $V(\Sigma,F)$, when
$\Z$-coefficients are used. A connected component of $\Sigma-K$
which is not connected to $\bdry \Sigma$ is called an {\em isolated
region of $K$ in $\Sigma$}.  We say that $K$ is {\em isolating} if
there is an isolated region of $K$ in $\Sigma$, and {\em
nonisolating} if there is no isolated region. For example, if $K$
has a homotopically trivial closed curve, then it is isolating.

We then have the following:

\begin{prop} \label{prop: nonisolating}
With $\Z$-coefficients, the dividing set $K$ has nonzero and
primitive $c(K)$ if $K$ is nonisolating.
\end{prop}

\begin{proof}
Suppose $\Sigma$ is connected. (If $\Sigma$ is not, we consider each
component of $\Sigma$ separately.) If $(\Sigma,F)=(D^2,F)$, then we
are done by Proposition~\ref{prop: crossingless}. Therefore, suppose
$\Sigma\not=D^2$. In view of Lemma~\ref{lemma: cutting}, it suffices
to find a properly embedded arc $\gamma\subset \Sigma$ which
intersects $K$ exactly once, so that cutting along it increases the
Euler characteristic of $\Sigma$ by one.  Let $\Sigma_0$ be a
connected component of $\Sigma-K$ which has Euler characteristic
$\not=1$.  Since $K$ is nonisolating, $\Sigma_0$ must nontrivially
intersect $\bdry\Sigma$.  It is then easy to find a properly
embedded arc $\gamma_0\subset \Sigma$ which lies in $\Sigma_0$, and
which is not $\bdry$-parallel in $\Sigma_0$.  We can isotop the
endpoints of $\gamma_0$ along $\bdry \Sigma$ so the resulting
$\gamma$ intersects $K$ exactly once.
\end{proof}

We also have the following corollary of Lemma~\ref{lemma: cutting}:

\begin{prop}\label{prop: generators}
With $\Z$-coefficients, $V(\Sigma,F)$ is generated $c(K)$, where $K$
ranges over all dividing sets for which $\bdry K=F$.
\end{prop}

\begin{proof}
The assertion is clearly true when $\Sigma=D^2$ and $\#F=2$ or $4$.
Now, any $(\Sigma,F)$ can be split along an arc $\gamma$ so that the
resulting $(\Sigma',F')$ satisfies $\chi(\Sigma')=\chi(\Sigma)+1$
and $\#F'=\#F+2$, and so that $V(\Sigma',F')\simeq V(\Sigma,F)$.
Once we reach $\Sigma'=D^2$, a good choice of splitting will
decrease $\#F'$ of each component, until each component is $(D^2,F)$
with $\#F=2$ or $4$. The proposition follows by gluing.
\end{proof}

\subsection{Analysis when $\Sigma$ is an annulus} \label{subsection:
annulus} Suppose $\Sigma$ is an annulus.  We consider the situation
where $F$ consists of two points on each boundary component.  The
calculations will be done in $\Z$-coefficients, but calculations in
a twisted coefficient system will certainly yield more information.
See for example \cite{GH}.

By Juh\'asz' formula, $V(\Sigma,F)=\Z^2\otimes \Z^2=\Z_{(2)}\oplus
\Z^2_{(0)}\oplus \Z_{(-2)}$. One can easily see that $\Z_{(2)}$ is
generated by a $\bdry$-parallel $K_+$ with two positive
$\bdry$-parallel arcs, and $\Z_{(-2)}$ is generated by a
$\bdry$-parallel $K_-$ with two negative $\bdry$-parallel arcs.

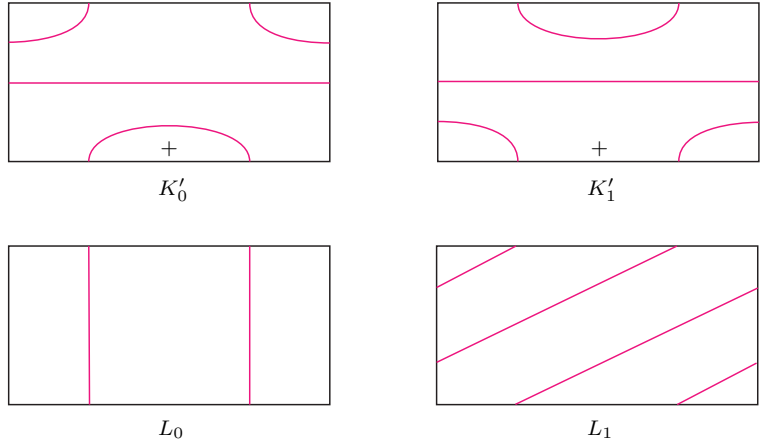
\begin{figure}[ht]
\begin{overpic}[width=10cm]{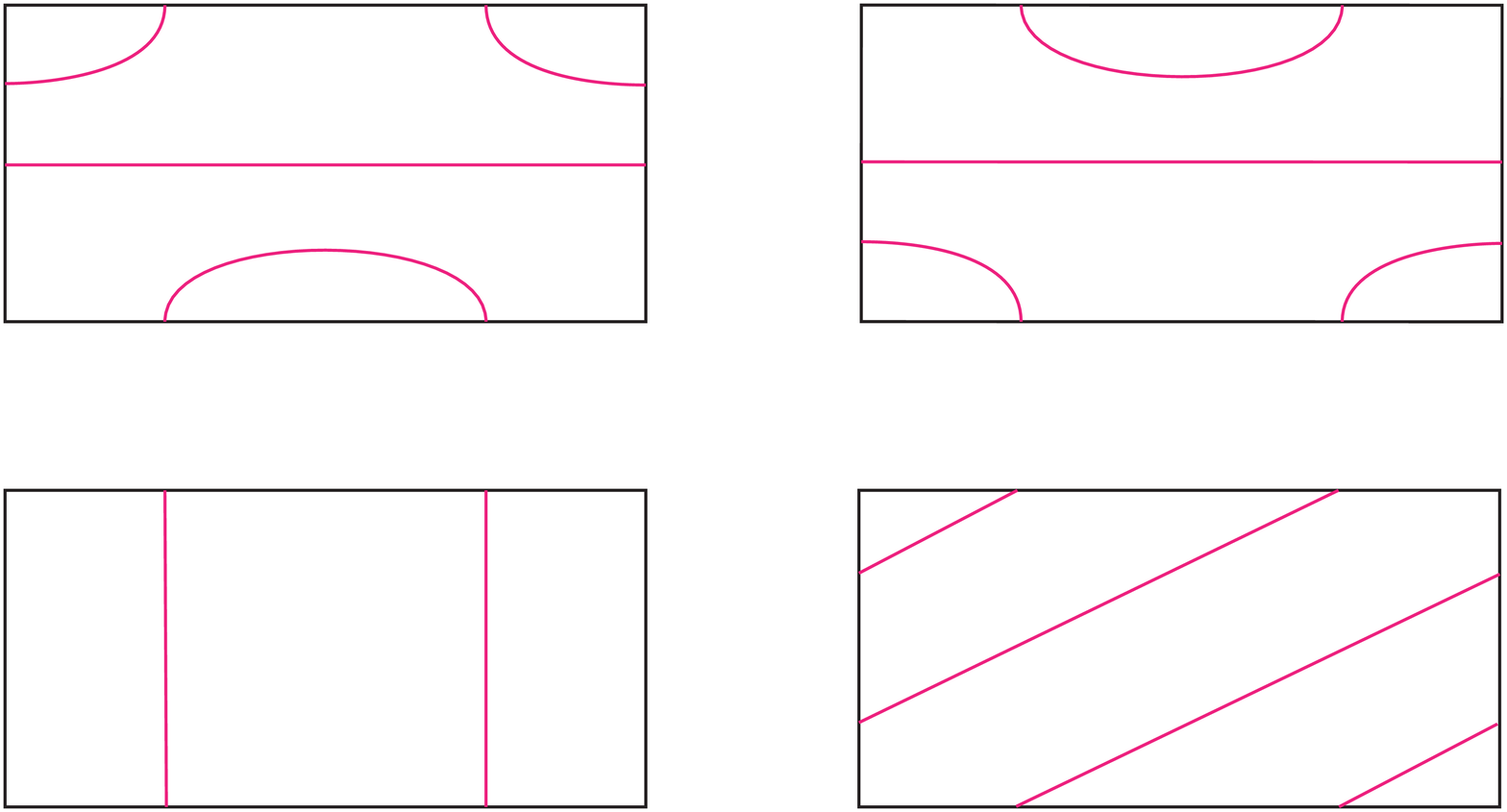}
\put(20,28){\tiny $K_0'$} \put(20.3,33.5){\tiny $+$}
\put(77,28){\tiny $K_1'$} \put(20,-4){\tiny $L_0$} \put(77,-4){\tiny
$L_1$} \put(77.5,33.5){\tiny $+$}
\end{overpic}
\s
\caption{The sides of each annulus are identified.} \label{annulus}
\end{figure}

It remains to analyze $\Z^2_{(0)}$.  The  nonisolating dividing sets
$K$ with nontrivial $c(K)\subset \Z^2_{(0)}$ are the following:
$K_0'$ and $K_1'$, which have two $\bdry$-parallel arcs of opposite
sign and one closed curve, $L_0$ consisting of two parallel arcs
from one boundary component to the other, as well as $L_j$, obtained
from $L_0$ by performing $j$ positive Dehn twists about the core
curve of the annulus. See Figure~\ref{annulus}. The other possible
dividing sets $K$, besides those with homotopically trivial
components, have at least two parallel closed curves. The
corresponding contact structure will necessarily have at least
$2\pi$-torsion. It was proved in \cite{GHV} that any contact
structure with $2\pi$-torsion has vanishing contact invariant over
$\Z$.

First consider the map $\Phi: \Z^2_{(0)}\rightarrow \Z_{(2)}$,
obtained by attaching an annulus with configuration $K_+$ from
below.  Since a homotopically trivial curve is created, we have
$\Phi(c(K_0'))=\{0\}$.  Also, $\Phi(c(K_1'))=\{0\}$, since the
resulting dividing set will have two parallel closed curves. On the
other hand, the $c(L_i)$ all map to the generator $c(K_+)$ of
$\Z_{(2)}$. Hence the map $\Phi$ is surjective, and must have $\ker
\Phi\simeq \Z$. Next, since $K_0'$ and $K_1'$ are nonisolating,
$c(K_0')$ and $c(K_1')$ must be primitive; this implies that
$c(K_0')=c(K_1')$ and generate $\ker \Phi$. Now, one can make a
coordinate change if necessary so that $c(L_0)=\{\pm (1,0)\}$ and
$c(K_i')=\{\pm (0,1)\}$.

Next we compute $c(L_1)$.  For this, we use Lemma~\ref{lemma: n is
three} and the following fact which follows from the proof of
Theorem~\ref{thm: unique lifting}: for $\Sigma=D^2$ and $n=3$, a
representative $\overline{c}(K_3)$ of $c(K_3)$ is a superposition of
representatives $\overline{c}(K_1)$ and $\overline{c}(K_2)$ of
$c(K_1)$ and $c(K_2)$ with $\pm 1$ coefficients. Observe that $K_1$
is obtained from $K_3$ by applying a bypass attachment from the
front, and $K_2$ is obtained from $K_3$ by a bypass attachment to
the back.  It is easy to see from the $\Phi$ in the previous
paragraph that $c(L_1)=\{\pm (1,n)\}$ for some integer $n$. Given
the configuration $K_0'$, take a bypass arc of attachment $\delta$
with endpoints on the two $\bdry$-parallel arcs and one other
intersection point with $K_0'$, namely along the closed component.
Take a small disk $D^2$ about $\delta$. Consider the gluing map
$$\Psi:V(D^2,K_0'|_{\bdry D^2})\otimes V(\Sigma-D^2,K_0'|_{\bdry
\Sigma-D^2}\cup F)\rightarrow V(\Sigma,F).$$ By tensoring the
$\overline{c}(K_i)$ with $\overline{c}(K_0'|_{\Sigma-D^2})$, the
equation $\overline{c}(K_3)=\pm \overline{c}(K_1)\pm
\overline{c}(K_2)$ becomes
$$\overline{c}(K_0')= \pm \overline{c}(L_0)\pm \overline{c}(L_1).$$
This means $(0,1)= \pm (1,0)\pm (1,n)$.  The only possible solutions
are $(0,1)= (1,0)-(1,-1)$ or $(0,1)=-(1,0)+(1,1)$.  (The two
possibilities are equivalent after a basis change.) Hence
$c(L_1)=\{\pm (1,1)\}$, for example.

In general, we conjecture that $c(L_n)=\{\pm (1,n)\}$. A proof of
this conjecture requires a more careful sign analysis than we are
willing to do for the moment.

\subsection{Determination of nonzero elements $c(K)$ in $V(\Sigma,F)$}
In this section we prove Theorem~\ref{thm: vanishing}, i.e., we
determine exactly which elements $K$ have nonzero invariants $c(K)$
in $V(\Sigma,F)$ with $\Z/2\Z$-coefficients.

\begin{prop} \label{prop: isolating}
If $K$ is isolating, then $c(K)=0$ with $\Z/2\Z$-coefficients.
\end{prop}

\begin{proof}
Suppose first that there is an isolated region $\Sigma_0$ which is
an annulus.  In that case, take an arc of attachment $\delta$ of a
bypass which intersects the two boundary components of $\Sigma_0$,
and some other component of $K$, in that order. By the TQFT property
applied to a small neighborhood $D$ of $\delta$ and $\Sigma-D$, we
see that if $K'$ (resp.\ $K''$) is obtained from $K$ by applying a
bypass from the front (resp.\ bypass to the back), then
$c(K)=c(K')+c(K'')$, since the corresponding fact is true on $D$.
One easily sees that $K'$ and $K''$ are isotopic, and is $K$ with
$\bdry \Sigma_0$ removed. With $\Z/2\Z$-coefficients, then,
$c(K)=2c(K')=0$.

Next suppose that $\Sigma_0$ has more than one boundary component,
and is outermost among all isolated regions, in the sense that one
boundary component $\gamma$ of $\Sigma_0$ is adjacent to a component
$\Sigma_1$ whose boundary intersects $\bdry\Sigma$. Also suppose
that $\Sigma_0$ is not an annulus. Take an arc of attachment
$\delta$ which begins on $\gamma$, intersects $\gamma$ after
traveling inside $\Sigma_0$, and ends on an arc component of $K$ on
$\bdry \Sigma_1$. Choose $\delta$ so that $\Sigma_0-\delta$ has two
components, one which is an annulus and the other which has Euler
characteristic $>\chi(\Sigma_0)$.  Then apply the bypass attachments
from the front and to the back to obtain $K', K''$ as in the
previous paragraph. Now, $c(K)=c(K')+c(K'')$, and one of $c(K')$ or
$c(K'')$ is zero, since it possesses an annular isolated region.
This reduces the number of components of $\bdry\Sigma_0$.

Finally suppose that $\bdry\Sigma_0$ is connected. If $\Sigma_0$
bounds a surface of genus $g>1$, then the above procedure can split
$c(K)=c(K')+c(K'')$, where both $c(K')$ and $c(K'')$ have isolated
regions with connected $\bdry\Sigma_0$ and strictly smaller genus.
Hence suppose that $\Sigma_0$ bounds a once-punctured torus.  Also,
by cutting along arcs as in Proposition~\ref{prop: nonisolating}, we
may assume that $\Sigma$ itself is a once-punctured torus with one
$\bdry$-parallel arc and one closed curve parallel to the boundary.
Choose $\delta$ as given in Figure~\ref{punctured-torus}, namely,
$\delta$ begins on the $\bdry$-parallel arc and intersects
$\bdry\Sigma_0$ twice, and restricts to a nontrivial arc on
$\Sigma_0$.
\begin{figure}[ht]
\begin{overpic}[width=11.5cm]{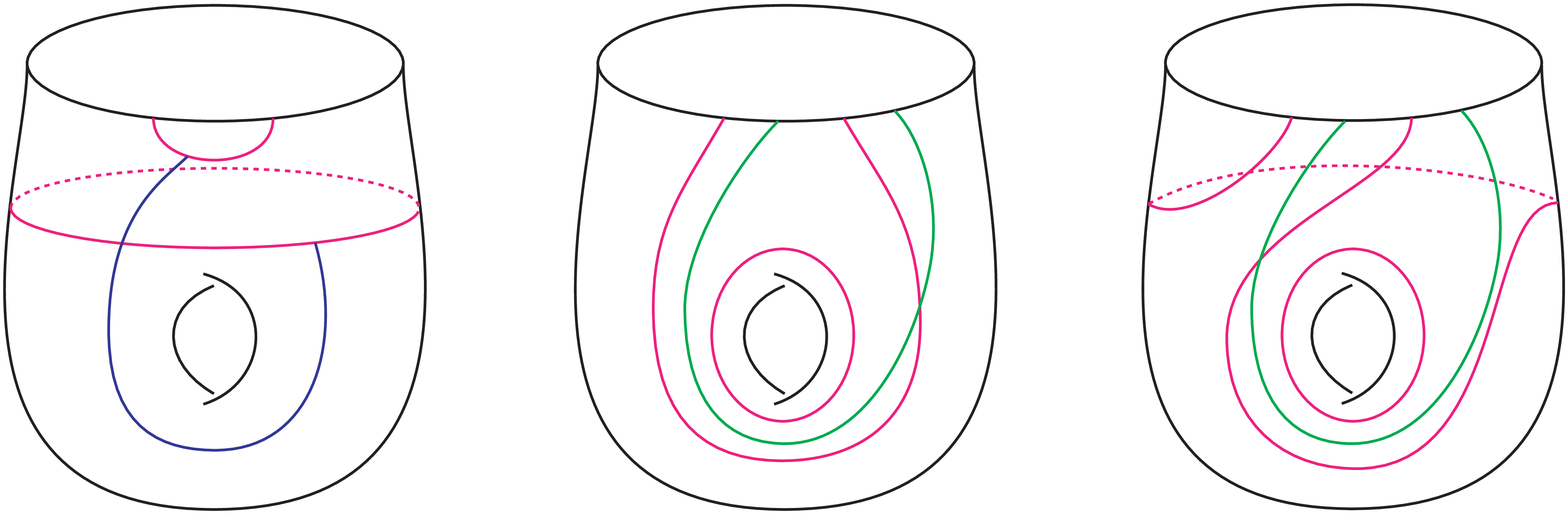}
\put(5.7,5){\tiny $\delta$} \put(59.5,22){\tiny $\tau$}
\put(95.5,22){\tiny $\tau$}
\end{overpic}
\caption{} \label{punctured-torus}
\end{figure}
The resulting $K'$ and $K''$ are the center and right-hand diagrams.
Now cut along the properly embedded, non-boundary-parallel arc
$\tau$ which intersects each of $K'$ and $K''$ exactly once.
Applying Lemma~\ref{lemma: cutting}, we see that $c(K')=c(K'')$ if
and only if the cut-open dividing curves $K'_0$ and $K''_0$ have
equal invariants in the cut-open surface. Finally, observe that, on
the cut-open surface (an annulus), $c(K'_0)=c(K''_0)$ since they
correspond to $K_0$ and $K_1$, discussed in
Subsection~\ref{subsection: annulus}.
\end{proof}

Propositions~\ref{prop: isolating} and ~\ref{prop: nonisolating},
together give Theorem~\ref{thm: vanishing}.

In the case of  $\Z$-coefficients we expect the following to hold:

\begin{conj}
Over $\Z$-coefficients, the following are equivalent: \be
\item $c(K)\not=0$;
\item $c(K)$ is primitive;
\item $K$ is nonisolating.
\ee
\end{conj}

The difficulty comes from not being able to determine whether $c(K)$
is divisible by $2$ with $\Z$-coefficients, which in turn stems from
our $\pm 1$ difficulty in Subsection~\ref{subsection: Z}. When
twisted coefficients are used, the result is quite different, and
will yield substantially more information \cite{GH}.

\maketitle \s\n {\em Acknowledgements.} We thank John Etnyre,
Andr\'as Juh\'asz, and Andr\'as Stipsicz for helpful discussions. KH
thanks Francis Bonahon and Toshitake Kohno for helpful discussions
on TQFT.  KH also thanks Takashi Tsuboi and the University of Tokyo
for their hospitality; much of the writing of this paper was done
during his five-month stay in Tokyo in the summer of 2007.

\end{document}